\documentclass[10pt]{amsart}
\usepackage{rotating}
\usepackage{mathcomp,amscd}
\usepackage{amssymb}
\usepackage{mathtools}
\usepackage[all]{xy}
\usepackage{subfigure}
\usepackage[margin=1cm]{caption}
\usepackage{wasysym}
\usepackage{color}
\definecolor{dg}{rgb}{0,0.67,0}
\definecolor{gr}{rgb}{.67,0.67,0.67}
\usepackage{amsfonts}
\usepackage{url}
\usepackage{graphicx}
\usepackage{amssymb}
\usepackage{epstopdf}
\usepackage{array}
\usepackage{multirow}

\newcommand{\Br}{\mathrm{Br}}
\newcommand{\br}{\mathrm{HBr}}

\newcommand{\GZ}{T}

\newcommand{\Out}{\mathrm{Out}}

\newcommand{\Hur}{{\mbox{\sc Hur}}}

\newcommand{\AHur}{{\sf Hur}}
\newcommand{\APGL}{{\sf PGL_2}}

\newcommand{\Conf}{{\mbox{\sc Conf}}}
\newcommand{\AConf}{{\sf Conf}}

\usepackage[OT2,T1]{fontenc}
\DeclareSymbolFont{cyrletters}{OT2}{wncyr}{m}{n}
\DeclareFontFamily{OT1}{rsfs}{}
\DeclareFontEncoding{OT2}{}{} 
  
     \DeclareFontShape{OT1}{rsfs}{n}{it}{<-> rsfs10}{}
\DeclareMathAlphabet{\mathscr}{OT1}{rsfs}{n}{it}

\newcommand{\C}{{\Bbb C}}
\newcommand{\Z}{{\Bbb Z}}
\newcommand{\R}{{\Bbb R}}

\newcommand{\cG}{\mathcal{G}}

\newcommand{\cF}{\mathcal{F}}

\newcommand{\F}{{\Bbb F}}

\newcommand{\bbP}{{\Bbb P}}
\newcommand{\AP}{{\sf P}}  
\newcommand{\SP}{{\mbox{\sc P}}}
\newcommand{\AX}{{\sf X}}

\newcommand{\SX}{{\mbox{\sc X}}}
\newcommand{\AY}{{\sf Y}}

\newcommand{\AU}{{\sf U}}
\newcommand{\SU}{{\mbox{\sc U}}}

\newcommand{\Q}{{\Bbb Q}}
 
\newcommand{\Gal}{\mbox{Gal}}
\newcommand{\Aut}{\mbox{Aut}}

\newcommand{\cP}{\mathcal{P}}

\newcommand{\disc}{\mbox{disc}}

\numberwithin{equation}{section}
\numberwithin{table}{section}
\numberwithin{figure}{section}
\newtheorem{Theorem}{Theorem}[section]
\newtheorem*{Theorem*}{Theorem}
\newtheorem*{corollary*}{Corollary} 
\newtheorem{Corollary}{Corollary}[section] 

\newtheorem*{Proposition*}{Proposition}
 
\newtheorem{Definition}[Theorem]{Definition}
\newtheorem{Conjecture}[Theorem]{Conjecture}

\newcommand{\cmmt}[1]{}
 \allowdisplaybreaks
\title{Hurwitz-Belyi maps}
\author{David P. Roberts}
\address{Division of Science and Mathematics, University of Minnesota Morris; 
Morris, Minnesota, 56267, USA}
\email{roberts@morris.umn.edu}
\urladdr{http://cda.morris.umn.edu/~roberts}
\begin{document}
\begin{abstract}
The study of the moduli of covers of the projective line leads to 
the theory of  Hurwitz varieties covering 
configuration varieties.    Certain one-dimensional
slices of these coverings are particularly interesting Belyi maps.    
We present systematic examples of such ``Hurwitz-Belyi maps.''
Our examples illustrate a wide variety of theoretical phenomena and computational
techniques.   
 \end{abstract}

\maketitle
\tableofcontents

\sloppy

\section{Introduction} 
The theory of Belyi maps sits at an attractive intersection in mathematics 
where group theory, algebraic geometry, and number theory all play 
fundamental roles.  In this paper we first introduce a simply-indexed class
of particularly interesting Belyi maps which arise in solutions of Hurwitz
 moduli problems. 
Our main focus is then the systematic computation 
of sample Hurwitz-Belyi maps and the explicit exhibition of their remarkable properties. 
We expect that our exploratory work here will support future more
theoretical studies.    
We conclude this paper by speculating that as degrees become 
large, Hurwitz-Belyi maps become extreme outliers among all
Belyi maps.   The rest of the introduction amplifies
on this first paragraph.

\subsection{Belyi maps} 
\label{Belyimaps}
In the classical theory of smooth projective complex algebraic curves,
ramified covering maps  from a given curve $\AY$ to the projective line $\AP^1$ 
play a prominent role.   If $\AY$ is connected with genus $g$, then 
any degree $n$ map $F : \AY \rightarrow \AP^1$ has  $2n+2g-2$ 
critical points in $\AY$, counting multiplicities.  For generic $F$, 
these critical points $y_i$ are all distinct and moreover the
critical values $F(y_i)$ are also also distinct.   
A {\em Belyi map} by definition is a map 
$\AY \rightarrow \AP^1$
having all critical values in $\{0,1,\infty\}$.  One should 
think of Belyi maps as the maps which are as far from generic
as possible, with their critical values being moreover
normalized to a standard position.  

 The term {\em Belyi map} has become standard in acknowledgement of 
  the
 importance of a theorem proved by Belyi \cite{Bel79} in the late 1970s.  
 This theorem says that a curve $\AY$ is the domain of 
 a Belyi map if and only if $\AY$ is defined over $\overline{\Q}$, the 
 algebraic closure of $\Q$ in $\C$.    
About Belyi's theorem, Grothendieck wrote \cite[p.15]{Gro97}, ``never, without
a doubt, was a such a deep and disconcerting result 
proved in so few lines!''   He went on to describe the 
``essential message'' of Belyi's theorem as a ``profound identity'' between
different subfields of mathematics.   The recent paper \cite{SV14} provides a 
computationally-focused survey of Belyi maps, 
with many references.

\subsection{An example}   The main focus of this
paper is the explicit construction of Belyi maps with
certain extreme properties.  A map from \cite{RobCheb} arises from
outside the main context of this paper but still exhibits 
these extremes.  It serves as a useful initial example: 
\begin{eqnarray}
\nonumber  \pi : \AP^1 & \rightarrow & \AP^1, \\
\label{U89} x & \mapsto &    \frac{ (x+2)^9 x^{18}  \left(x^2-2\right)^{18} (x-2)}{(x+1)^{16} \left(x^3-3 x+1\right)^{16}}.
\end{eqnarray}
The degree is $64$ and the $126$ critical points 
are easily identified as follows.   From the numerator $A(x)$,
one has the critical points $-2$, $0$,  $\sqrt{2}$, $-\sqrt{-2}$,
with total multiplicity  $8+17+17+17=59$ and critical 
value $0$.  From the denominator $C(x)$, one has critical 
points $-1$, $x_2$, $x_3$, $x_4$ with total 
multiplicity $15+15+15+15=60$ and critical value $\infty$.
Since both $A(x)$ and $C(x)$ are monic, one has 
$\pi(\infty)=1$.  The exact coefficients in \eqref{U89} are chosen so that the
degree of $A(x)-C(x)$ is only
$56$.  This means
that $\infty$ is a critical point of multiplicity $63-56=7$. 
As $59 + 60 + 7 = 126$, there can be no critical values
outside $\{0,1,\infty\}$ and so  
$\pi$ is indeed a Belyi map.  

In general, a degree $m$ Belyi map $\pi$ has a monodromy group $M_\pi \subseteq S_m$, 
a number field $F_\pi \subset \C$ of definition, and a finite 
set $\cP_\pi \subset \{2,3,5,\dots\}$ of bad primes. 
We call $\pi$ {\em full} if $M_\pi \in \{A_m,S_m\}$.  
Our example $\pi$ is full because $M_\pi = S_{64}$.  
It is defined over $F_\pi=\Q$ because all the coefficients in 
\eqref{U89} are in $\Q$.   It has bad reduction set $\cP_\pi = \{2,3\}$ 
because numerator and denominator have 
a common factor in $\F_p[x]$ exactly for $p \in \{2,3\}$.  
In the sequel, we almost always drop the subscript $\pi$,
as it is clear from context.  

To orient the reader, we remark that the great bulk of the explicit
literature on Belyi maps concerns maps which are not full.  
Much of this literature, for example \cite[Chapter II]{MM99}, focuses on Belyi maps with $M$ 
a finite simple group different from $A_m$.  On the other hand, seeking
Belyi maps defined over $\Q$ is a common focus in
the literature.  Similarly, preferring maps with
small bad reduction sets $\cP$ is a common viewpoint.  

\subsection{An inverse problem}
\label{inverseproblem}  To provide a framework for 
our computations, we pose the following inverse problem:
{\em given a finite set of primes $\cP$ and a degree $m$,
find all full degree $m$ Belyi maps $\pi$ defined over $\Q$ 
with bad reduction set within $\cP$.}    The finite 
set of full Belyi maps in a given degree $m$ is 
parameterized in an elementary way by group-theoretic
data.  So, in principle at least,  this problem is simply asking to extract those for which
$F_\pi = \Q$ and $\cP_\pi \subseteq \cP$.  

While the  Belyi map \eqref{U89} may look rather 
ordinary, it is already unusual for full Belyi
maps to be defined over $\Q$.  It seems to 
be extremely rare that their bad reduction set
is so  small.    In fact, we know
of no full Belyi maps defined over $\Q$ with $m \geq 4$ and 
 $|\cP_\pi| \leq 1$.   For $|\cP_\pi|=2$ we know 
 of only a very sparse collection of such maps \cite{MR05}, \cite{RobCheb},
 as discussed further in our last section here.
The largest degree of these with both primes less than seventeen is $m=64$, 
 coming from \eqref{U89}.

\subsection{Hurwitz-Belyi maps}   Suppose now that 
$\cP$ contains the set $\cP_T$ of primes dividing
the order of a finite nonabelian simple group $T$. 
The theoretical setting for this paper is a systematic method of constructing
Belyi maps of arbitrarily large degree defined over $\Q$ and
 ramified within $\cP$.  
This method is
to extract Belyi maps $\AX \rightarrow \AP^1$ from solutions to 
Hurwitz moduli problems involving 
maps $\AY \rightarrow \AP^1$ having $r \geq 4$ critical values
and  monodromy
group suitably close to $T$.
From the full monodromy theorem of \cite{RV15}, we expect 
that these Hurwitz-Belyi maps are typically full.

 \subsection{Contents of this paper} 
 Our viewpoint is that Hurwitz-Belyi maps form a remarkable
 class of mathematical objects, and are worth studying
 in all their aspects.   This paper focuses on presenting explicit defining
 equations for systematic collections of Hurwitz-Belyi maps, 
 and exhibits a number of theoretical structures in the process. 
 The defining  equations are obtained by two 
 complementary methods.  What we call the {\em standard method} centers
 on algebraic computations directly with the
  $r$-point  Hurwitz source.  The {\em braid-triple method} 
  is a novelty of this paper.  It
 uses the $r$-point Hurwitz source only to 
 give necessary braid group information;  its remaining
 computations are then the same ones
 used to compute general Belyi maps.   
 
 We focus primarily on the case $r=4$ 
 which is the easiest case for computations for a given $T$. 
 This case was studied in some generality  by
 Lando and Zvonkin in \cite[\S5.5]{LZ04} under the term {\em megamap}.  
 In the last two sections, we shift the focus 
 to $r \geq 5$, which is necessary to obtain the
 very large degrees $m$ we are most interested in.
 The standard method is
  insensitive to genera of covering
 curves $\AX$, and so we could
 easily present examples of quite high
 genus.  However, to give 
 a uniform tidiness to our final equations,
 we present defining equations only
 in the case of genus zero.   
 Thus the reader will find many
 explicit rational functions 
 in $\Q(x)$ with properties 
 similar to those of our initial
 example \eqref{U89}.  All these 
 rational functions and related information are
 available in the {\em Mathematica} file \verb+HMB.mma+ file on the author's homepage.

 Section~\ref{Belyi} reviews the theory of Belyi maps.
 Section~\ref{HurwitzBelyi1} reviews the theory of 
 Hurwitz maps and explains how carefully 
 chosen one-dimensional slices are
 Hurwitz-Belyi maps.  
   Of the many Belyi 
maps appearing
  in Section~\ref{Belyi}, two
are unexpectedly defined over $\Q$.   These maps 
 each appear again in Section~\ref{HurwitzBelyi2}, with now their 
 rationality obvious from the beginning via the Hurwitz
 theory. 
 
 Section~\ref{Cubical} illustrates the phenomenon that sometimes 
 infinitely many  Hurwitz-Belyi maps form a ``clan'' in 
 that they can be described  simultaneously via formulas uniform in certain 
 parameters.  We illustrate this phenomenon
 with a three-parameter clan extracted from a four-parameter
 clan studied by Couveignes \cite{Cou97}.  
 We find defining equations via the standard method.  
 The simple groups $T$ involved here are $A_n$.
 Interestingly, the bad reduction
 sets $\cP$ are substantially smaller than $\cP_T$. 
 However for any given $\cP$, the clan gives only 
 finitely many Belyi maps ramified within $\cP$.

 Section~\ref{4vs3} introduces the alternative
 braid-triple method for finding defining equations.  
 We give general formulas for the preliminary
 braid computations in the setting $r=4$.
 Passing from braid information to defining equations can then be
 much more computationally demanding than
 in our initial examples, and we find equations 
 mainly by  $p$-adic techniques.  
 Section~\ref{separation} then presents three examples for which both methods  work,
 with these examples having the added interest that
 lifting invariants force $\AX$ to be disconnected.
 In each case, $\AX$ in fact has two components, each
 of which is full over the base projective line.  
 
 Sections~\ref{23p}, \ref{235}, and \ref{237} consider
a  systematic collection
 of Hurwitz-Belyi maps, with all 
 final equations computed by 
 the braid-triple method.  
 They focus on the cases where 
 $|\cP_T| \leq 3$.  By the classification
 of finite simple groups, the possible
 $\cP_T$ have the form $\{2,3,p\}$ with
 $p \in \{5,7,13,17\}$.  Section~\ref{23p}
 sets up our framework and 
 presents one example each for 
 $p=13$ and $p=17$.  Sections~\ref{235} 
 and \ref{237} then give many examples for 
 $p=5$ and $p=7$ respectively.   

Section~\ref{fivepoint} takes first computational 
steps into the setting $r \geq 5$.  
Working just with $T = A_5$ 
and $r=5$, we summarize braid computations
which easily prove the existence of 
full Hurwitz-Belyi maps  with
bad reduction set $\{2,3,5\}$ and 
degrees into the thousands.   
We use the standard method to find
equations of two such covers related to $T=A_6$,  
one in degree $96$ and the other of 
degree $192$. 

Section~\ref{conj} concludes by tying the considerations
of this paper very tightly to those of \cite{RV15} and \cite{RobHNF}.
It conjectures a direct analog for Belyi maps
of the main conjecture there
for number fields.  The Belyi map conjecture
responds to the above inverse problem
in the case that $\cP$ contains the 
set of primes dividing the order of some
finite nonabelian simple group.  In particular, it says
that there then should be full Belyi maps
defined over $\Q$ and ramified within
$\cP$ of arbitrarily large degree.

\subsection{Notation}
\label{notation} Despite the arithmetic nature 
of our subject, we work as much as possible
over $\C$.  We use a sans serif font for 
complex spaces, as in $\AY$, $\AX$, $\AP^1$
above, or $\AHur_h$ and $\AConf_\nu$ below.  
Projective lines enter our considerations
in many ways.  When useful, we distinguish 
projective lines by the coordinates
we are using on them, as in 
$\AP^1_y$, $\AP^1_x$, $\AP^1_t$, $\AP^1_v$,
$\AP^1_w$, or $\AP^1_j$.  
More general curves also enter in more than one way.
As above, we mostly reserve $\pi : \AX \rightarrow \AP^1$
for Belyi maps, and use
$F : \AY \rightarrow \AP^1$ for general maps.   

The phenomenon that allows us to work 
mainly over $\C$ is that to a great extent
geometry determines arithmetic.    
Thus an effort
to find a function $\pi(x) \in \C(x)$ 
giving a full Belyi map 
$\pi: \AP^1_x \rightarrow \AP^1$
involves choices of
normalization.   Typically, one can make 
these choices in a  geometrically
natural way, and then the coefficients of 
$\pi(x)$ automatically span 
the field of definition.  When this field
is $\Q$, and the normalization
is sufficiently canonical, the
primes of bad reduction 
can be similarly read off.  
In the rare cases that 
it is really necessary
to use the language of
algebraic varieties
over $\Q$, 
we 
use a different font,
as in 
$\Hur_h \rightarrow \Conf_\nu$
or $\SX \rightarrow \bbP^1$.

Partitions play a prominent role throughout this paper, and we use
several standard notations interchangeably, choosing the 
one that is most convenient for the current context.  Thus 
$(5,3,2,2)$,  $5322$, and  $532^2$ all represent the
same partition of $12$.  Our notation for groups closely 
follows Atlas \cite{ATLAS} notation, with the main difference
being a greater emphasis on partitions.  As a smaller deviation,
we help distinguish classes in a simple group $T$ from classes in
$T.2-T$ by using small letters for the former.    For 
example, we denote the conjugacy classes 
$1A$, $2A$, $3A$, $5A$ and $5B$ of \cite{ATLAS} in $A_5$ by 
$11111$, $221$, $311$, $5a$ and $5b$.  

We will commonly present a Belyi map $\pi : \AX \rightarrow \AP^1_v$ not
as a rational function $v=A(x)/C(x)$ but rather via
the corresponding polynomial equation $A(x)-v C(x)=0$.  
This trivial change in perspective has several advantages, 
one being that it lets one see the three-point
property and the primes of bad reduction simultaneously
via discriminants.  For example, the discriminant 
of $A(x)-v C(x)$ in our first example \eqref{U89}  is $2^{256} 3^{126} v^{59} (v-1)^7$.

\subsection{Acknowledgements} 
Conversations with Stefan Kr\"amer, Kay Magaard, Hartmut Monien, 
Sam Schiavone, Akshay Venkatesh, and John Voight were helpful 
in the preparation of this paper.    
The Simons foundation
partially supported this work through grant \#209472.

\section{Two Belyi maps unexpectedly defined over $\Q$} 
\label{Belyi}
This section presents thirty-one Belyi maps as explicit rational 
functions in $\C(y)$, two of which are unexpectedly in $\Q(y)$.  
Via these examples, it provides a quick summary, adapted to this paper's needs,
of the general theory of 
Belyi maps.     We will revisit the
two rational maps from a different point of view 
in Section~\ref{HurwitzBelyi2}.   Our three-point
computations here are providing models for later
$r$-point computations.  Accordingly, we use
the letter $y$ as a primary variable.  
\label{zeroparam}
 
\subsection{Partition triples} 
\label{partitiontriples}
Let $n$ be a positive integer.  Let $\Lambda = (\lambda_0,\lambda_1,\lambda_\infty)$
be a triple of partitions of $n$, with the $\lambda_\tau$ 
having all together $n+2-2g$ parts, with $g \in \Z_{\geq 0}$.
  The two examples 
pursued in this section are
\begin{align}
\label{twotrips1}
 \Lambda' & = (322, \, 421, \, 511), &
\Lambda'' & = (642, \, 2222211, \, 5322).
\end{align}
So the degrees of the examples are $n=7$ and $n=12$, and both
have $g=0$. 

Consider Belyi maps $F : \AY \rightarrow \AP^1$
with ramification numbers of the points in $F^{-1}(\tau)$ forming
the partition $\lambda_\tau$, for each $\tau \in \{0,1,\infty\}$. 
Up to isomorphism, there are only finitely many
such maps.   For some of these maps, $\AY$ may 
be disconnected, and we are not 
interested here in these degenerate cases. 
Accordingly, let $\AX$ be the set of isomorphism
classes of such Belyi maps with $\AY$ connected.   
One wants to explicitly identify $\AX$, and simultaneously
get an algebraic expression for each corresponding
Belyi map $F_x : \AY_x \rightarrow \AP^1$.    The Riemann-Hurwitz
formula says that all these $\AY_x$ have genus $g$.  

Computations can be put into a standard form when $g=0$ and
the partitions $\lambda_0$, $\lambda_1$, and $\lambda_\infty$
have in total at least three singletons.  Then one can pick an ordered
triple of singletons and coordinatize $\AY$ by choosing
$y$ to take the values $0$, $1$, and $\infty$ in order at the
three corresponding points.   In our two examples, we do this via
\begin{align}
\label{twotrips2}
\Lambda'_* & = (3_022, \; 4_12_x1, \; 5_\infty11), &
\Lambda''_* & = (6_04_12_x, \; 2222211, \; 5_\infty 322).
\end{align}
Also we have chosen a fourth point in each case and subscripted it
by $x$.   This choice gives a canonical map from $\AX$ into $\C$,
as will be illustrated in our two examples.   When the
map corresponding to such a marked triple $\Lambda_*$ is injective,
as it almost always seems to be, we say that $\Lambda_*$ is
a {\em cleanly marked} genus zero triple.   Equations~\eqref{hur} and \eqref{nonhur} together 
give an example where injectivity fails. 

When $g=0$ and there is at least one singleton $a$, computations
can be done very similarly. All the explicit examples of this paper 
are in this setting.  When $g=0$ and 
there are no singletons, one often has to take
extra steps, but the essence of the method remains
very similar.  
When $g>0$, computations are still possible,  but they are very much
more complicated.

\subsection{The triple $\Lambda'$ and its associated $4=3+1$ splitting} 
\label{deg7}  
The subscripted triple $\Lambda'_*$ in 
\eqref{twotrips2} requires us to consider rational functions
\[
F(y) =  \frac{1+c+d}{(1+a+b)^2} \cdot \frac{y^3 (y^2 + a y + b)^2}{y^2 + c y + d}
\]
and focus on the equation
\begin{equation}
\label{diff3}
5 y^4  + 3 (a+2c) y^3  +   (4 a c+b+7 d) y^2 + (5 a d+2 b c) y +3 b d \,  = \, 5 (y-1)^3 (y-x).
\end{equation}
The left side is a factor of the numerator of $F'(y)$ and thus its roots are critical points.  
The right side gives the required locations and multiplicities of these critical points. 

Equating coefficients of $y$ in \eqref{diff3} and using also $F(x)=1$ 
gives five equations in five unknowns.  There are four solutions,
indexed by the roots of 
\begin{equation}
\label{split3plus1}
f_{\Lambda_*'}(x) = (x+2) \left(16 x^3-248 x^2-77 x-6\right).
\end{equation}
In general from a cleanly marked genus zero triple $\Lambda_*$, one gets a 
separable {\em moduli
polynomial} $f_{\Lambda_*}(x)$.  The {\em moduli algebra}
\[
K_{\Lambda} = \Q[x]/f_{\Lambda_*}(x)
\]
depends, as indicated by the notation, only on $\Lambda$ and 
not on the marking.   It is well-defined in the general case 
when the genus is arbitrary, even though we are not giving
a procedure here to find a particular polynomial.    

While the computation just presented is typical, the final result is not.
We give three independent conceptual 
explanations for the factorization in \eqref{split3plus1}, two in \S\ref{deg7HB} and one at the end of 
\S\ref{cubic}.   For context, 
 the splitting of the moduli polynomial is one of just
 four unexplained splittings  on the fourteen-page table of moduli algebras 
 in \cite{Mal94}.   While here the 
 degree $7$ partition triple yields a moduli algebra splitting as $3+1$,
 in the other examples the degrees are 8, 9, and 9, and
 the moduli algebras split as $7+1$, $8+1$, and $8+1$.

\subsection{Dessins}
\label{dessins}
A Belyi map $F: \AY \rightarrow \AP_t^1$ can be visualized by
its {\em dessin} as follows.  
      Consider the interval $[0,1]$ in $\AP_t^1$ as the bipartite graph $\bullet \!\!\! - \!\!\!\! - \!\circ$.  
      Then
$\AY_{[0,1]} := F^{-1}([0,1])$ inherits the structure of a bipartite graph.  
This bipartite graph, considered always as inside the ambient real surface $\AY$, is 
the dessin associated to $F$.    A key property is
that $F$ is completely determined by the topology of the 
dessin.    

\begin{figure}[htb]
\includegraphics[width=5.4in]{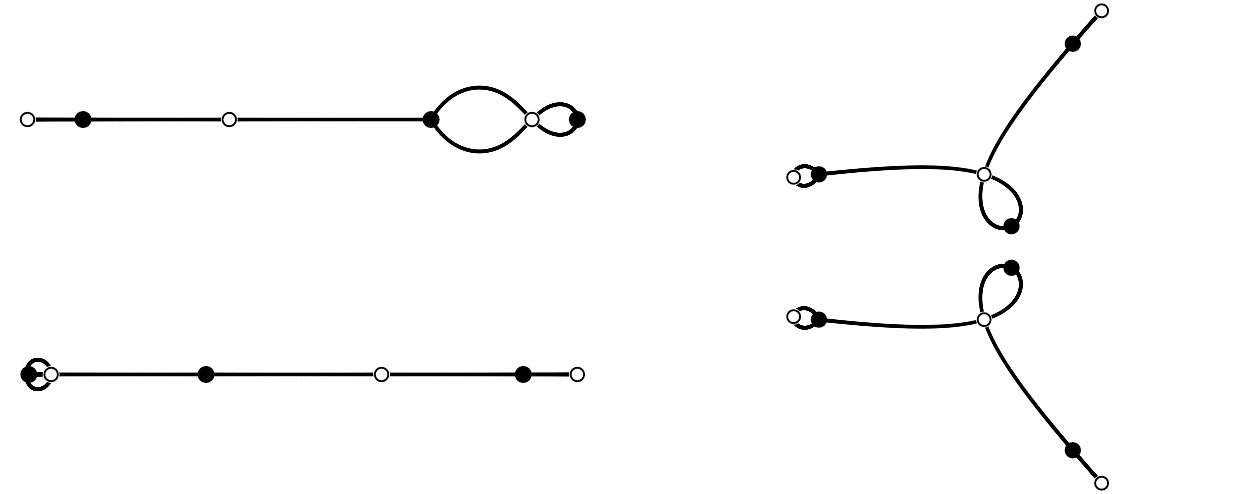}
\caption{\label{deg7four}  Dessins $\AY_{x_i,[0,1]} \subset \AP^1_y$ corresponding to the points of 
$\AX_{\Lambda'} = \{x_1,x_2,x_3,x_4\}$ 
with $\Lambda' =  (322, 421, 511)$ }
\end{figure}

Returning to the example of the previous subsection, the roots indexing the four solutions are
\begin{align*}
x_1 & = -2, &  x_2 & \approx 0.153 - 0.018 i, \\
x_3 & \approx 15.86, & x_4 & \approx 0.153 + 0.018i.
\end{align*}
The complete first solution is 
\begin{equation}
\label{deg7y}
F_1(y) = -\frac{y^3 \left(y^2+2 y-5\right)^2}{4 (2 y-1) (3 y-4)}.
\end{equation}
The coefficients of the other $F_i$ are cubic irrationalities.     
The four corresponding dessins in $\AY_i = \AP^1_y$ 
are drawn in Figure~\ref{deg7four}.  
The scales of the four dessins in terms of the common $y$-coordinate 
are quite different. 
Always the black triple point is at $0$ and the
white quadruple point is at $1$. The white double point
is then at $x_i$.  

\subsection{Monodromy} The dessins visually capture the group
theory which is central to the theory of Belyi 
maps but has not been mentioned so far.    Given 
a degree $n$ Belyi map $F : \AY \rightarrow \AP^1$,
consider the set $\AY_\star$ of the edges of the
dessin.   Let $g_0$ and $g_1$ be the operators 
on $\AY_\star$ given by 
rotating minimally counterclockwise about black
and white vertices respectively.   

   The choice of $[0,1]$ as the base graph is
 asymmetric with respect to the three critical values $0$, $1$, and $\infty$.
Orbits of $g_0$ and $g_1$ correspond to 
black vertices and white vertices respectively.  In 
our first example, the original partitions $\lambda_0 = 322$ and $\lambda_1 = 421$ 
can be recovered from 
each of the four dessins from the valencies of these vertices.  
On the other hand, the orbits of  $g_\infty = g_1^{-1} g_0^{-1}$ correspond to 
faces.   The valence of a face is by definition half the
number of edges encountered as one traverses its
 boundary.  Thus $\lambda_\infty = 511$ is recovered from 
 each of the four dessins in Figure~\ref{deg7four}, with the outer face always
 having valence five and the two bounded faces valence one.  
 
 Let $\mathcal{Y}_*$ be the set of ordered triples $(g_0,g_1,g_\infty)$ 
 in $S_n$ such that
 \begin{itemize}
 \item  $g_0$, $g_1$, and $g_\infty$ respectively
 have cycle type $\lambda_0$, $\lambda_1$, and $\lambda_\infty$,
 \item  $g_0g_1g_\infty=1$,
 \item  $\langle g_0, g_1 \rangle$ is a transitive subgroup of $S_n$.
 \end{itemize}
 Then $S_n$ acts on $\mathcal{Y}_*$ by simultaneous conjugation,
 and the quotient is canonically identified with $\AY_\star$.
    
 For each of the thirty-one dessins of this section, the monodromy group $\langle g_0, g_1 \rangle$ is all of
 $S_n$.  Indeed the only transitive subgroup of $S_7$ having the three cycle
 types of $\Lambda'$ is $S_7$, and the only transitive subgroup of $S_{12}$ having 
 the three cycle types of $\Lambda''$ is $S_{12}$.   
 
 \subsection{Galois action} 
 Let $\Gal(\overline{\Q}/\Q)$ be absolute Galois group of $\Q$.  
 The ``profound identity'' mentioned in the introduction
 centers on the fact that $\Gal(\overline{\Q}/\Q)$ acts
 naturally on the set $\AX$ of Belyi maps belonging 
 to any given $\Lambda$.  In the favorable cleanly marked
 situation set up in \S\ref{partitiontriples}, one has $X \subset \overline{\Q}$ 
 and the action on $\AX$ is the restriction
 of the standard action on $\overline{\Q}$. 
  
 A broad problem is to describe various 
 ways in which $\Gal(\overline{\Q}/\Q)$  may be forced
 to have more than one orbit.    Suppose 
 $x, x' \in \AX$ respectively give rise to 
 monodromy groups  $\langle g_0, g_1 \rangle$ and
 $\langle g'_0, g'_1 \rangle$.  If these 
 monodromy groups are not conjugate in 
 $S_n$ then certainly $x$ and $x'$ are in different 
 Galois orbits.   Malle's paper \cite{Mal94} repeatedly
 illustrates the next most common source
 of decompositions, namely
 symmetries with respect to 
 certain base-change operators 
 $\AP_t^1 \rightarrow \AP_t^1$.  The two splittings
 in this section do not come from either
 of these simple sources.  

\subsection{The triple $\Lambda''$ and its associated $24=23+1$ splitting}  
\label{deg12} Here we summarize the situation for $\Lambda''$.  Again 
the computation is completely typical, but the result is atypical.  The
clean marking on $\Lambda''$ identifies $\AX_{\Lambda''}$ with
the roots of 
 {\small
\begin{align*}
& (5 x+4) \cdot  \\ &  (48828125 x^{23}+283203125 x^{22}-4345703125
    x^{21}-21400390625 x^{20}+134842187500 x^{19} \\& +461968375000
    x^{18}-1670830050000 x^{17}-2095451850000 x^{16}+7249113240000
    x^{15} \\ & +6576215456000 x^{14}-23053309281280 x^{13}-10284915779584
    x^{12}+50191042453504 x^{11} \\ & +9449308979200 x^{10}-74715419574272
    x^9+5031544553472 x^8+71884253429760 x^7 \\ & -35243151065088 x^6-41613745192960
    x^5+29347637362688 x^4+14541349978112 x^3 \\ & +1765701320704 x^2+100126425088
    x+2684354560).
 \end{align*}
 }%
 
  \begin{figure}[htb]
\includegraphics[width=4.6in]{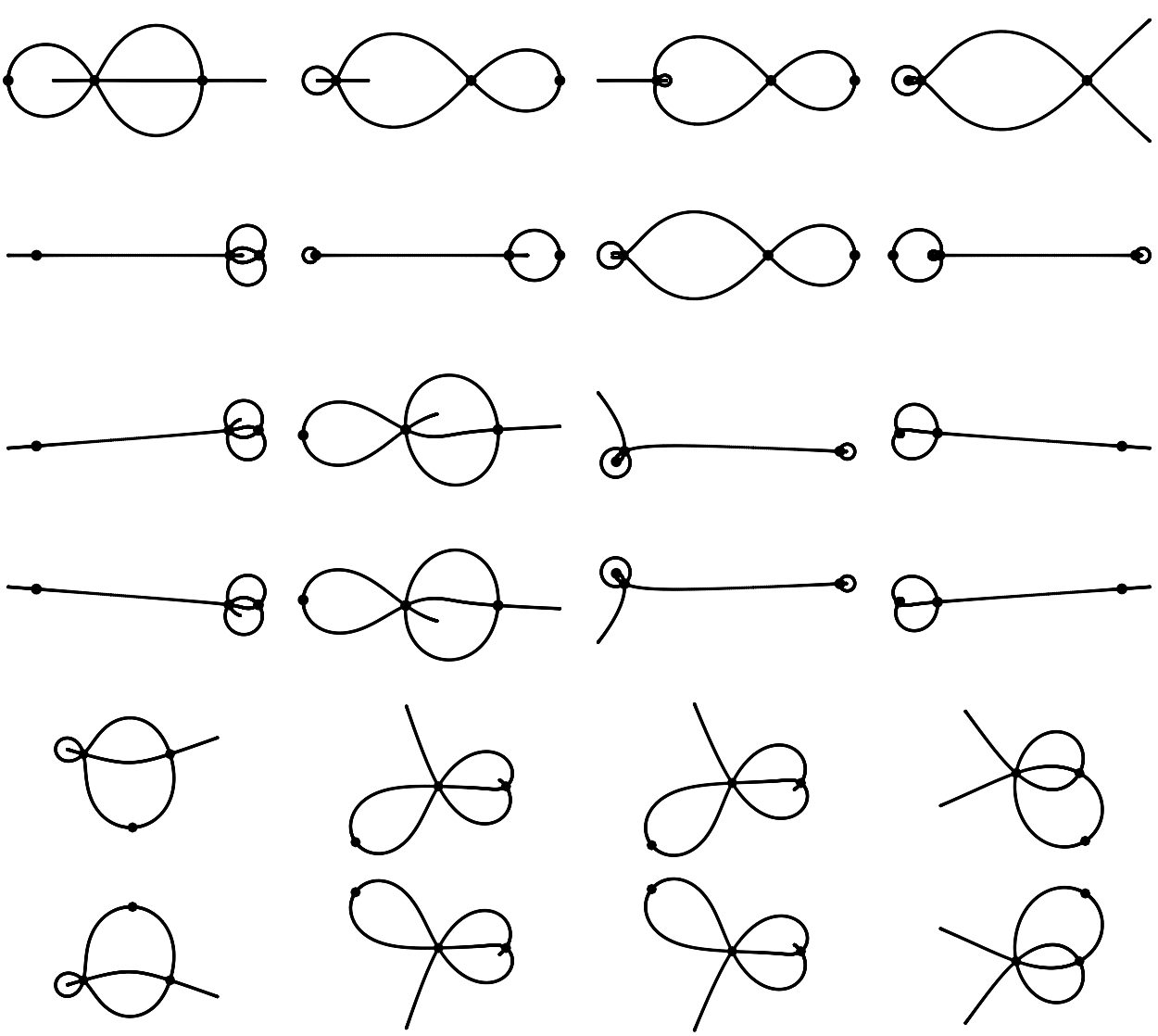}
\caption{\label{picts24}  Dessins in $\AY_{x,[0,1]} \subset \AP^1_y$ corresponding to the twenty-four points $x \in \AX_{\Lambda''}$ 
with $\Lambda'' = (642,2222211,5322)$ }
\end{figure}

 The twenty-four associated dessins are drawn in Figure~\ref{picts24}.
The cover $F_{-4/5} : \AP^1_y \rightarrow \AP^1_t$ is given by 
 \begin{equation}
 \label{deg12y} 
 t= \frac{5^5 y^6 (y-1)^4  (5 y+4)^2}{2^4 3^3 (2 y+1)^3 \left(5 y^2-6 y+2\right)^2}.
 \end{equation} 
 This splitting of one cover away from the other twenty-three covers is
explained in \S\ref{ambiguity}.

In choosing conventions for using dessins to represent covers, one often has to 
choose between competing virtues, such as symmetry
versus simplicity.  
Figure~\ref{picts24} represents the 
standard choice when $\lambda_1$ has the form $2^a 1^b$: one
draws the white vertices just as regular points, because they 
are not necessary for recovering the cover.  With this convention
there are just three highlighted points in each of the
dessins in Figure~\ref{picts24}:  black dots
of valence $6$, $4$, $2$ at $y=0$, $1$, $x$.   
The rational cover, with $x=-4/5$, appears in the upper left.     

\subsection{Bounds on bad reduction} 
\label{bounds} Let $n$ be a positive
integer and let $\Lambda = (\lambda_0,\lambda_1,\lambda_\infty)$ 
be a triple of partitions of $n$ as above.   Let $\cP^{\rm loc}$ be the 
set of primes dividing a part of one of the $\lambda_i$.  Let
$\cP^{\rm glob}$ be the set of primes less than or equal to 
$n$.   In our two examples $\cP^{\rm loc} = \{2,3,5\}$ and
$\cP^{\rm glob}$ is larger, by $\{7\}$ and $\{7,11\}$ respectively.   

Let $K_\Lambda$ be the moduli algebra associated to $\Lambda$.  
Let $D_\Lambda$ be its discriminant, i.e.\ the product of the
discriminants of the factor fields.  In our two examples,
$D_{\Lambda'} =  -2^3 \, 3 \; 5^3 \, 7$ and $D_{\Lambda''} = 2^{38} \, 3^{25} \, 5^{18} \, 7^{6}$.
Let $\cP_\Lambda$ be the set of primes dividing $D_{\Lambda}$.  Then 
one always has $\cP_\Lambda \subseteq \cP^{\rm glob}$.  Of 
course if $K_\Lambda = \Q$, then one has 
$\cP_{\Lambda} = \emptyset$.  Our experience is that 
once $[K_\Lambda:\Q]$ has moderately large degree,
$\cP_\Lambda$ is quite likely to be all or almost all of $\cP^{\rm glob}$,
as in the two examples.  

Suppose now that $\pi : \AY \rightarrow \AP^1$ is 
a Belyi map defined over $\Q$.   Then its set $\cP$ of bad
primes satisfies 
\begin{equation}
\label{reductionbounds}  \cP^{\rm loc} \subseteq \cP \subseteq \cP^{\rm glob}.
\end{equation}
For our two examples, $\cP$ coincides with its lower bound 
$\{2,3,5\}$.  The conceptual explanations of the 
splitting  given in Section~\ref{HurwitzBelyi1} also explain why
the remaining one or two primes in $\cP^{\rm glob}$ are
primes of good reduction.

\section{Hurwitz maps, Belyi pencils, and Hurwitz-Belyi maps}  
\label{HurwitzBelyi1}
    In \S\ref{Hurwitzmaps} we very briefly review the formalism of dealing with moduli
 of maps $\AY \rightarrow \AP_t^1$ with $r$ critical values.
 A key role is played by Hurwitz covering maps 
 $\pi_h : \AHur^*_h \rightarrow \AConf_\nu$.  
    In \S\ref{BP1} we introduce the concept of a Belyi pencil 
    $u : \AP^1-\{0,1,\infty\} \rightarrow \AConf_\nu$ 
    and in \S\ref{BP2} we give three important examples in
    $r=4$.  Finally \S\ref{HB} combines the notion of 
    Hurwitz map and Belyi pencil in a straightforward
    way to obtain the general notion of a Hurwitz-Belyi map $\pi_{h,u}$.  
 
\subsection{ Hurwitz maps}
\label{Hurwitzmaps}
Consider a general 
degree $n$ map $F : \AY \rightarrow \AP_t^1$ as 
in \S\ref{Belyimaps}.  Three fundamental invariants
are 
\begin{itemize}
\item Its global monodromy group $G \subseteq S_n$. 
\item The list $C = (C_1,\dots,C_k)$ 
 of distinct conjugacy classes of $G$ arising as non-identity local monodromy 
 transformations.
\item The corresponding list $(D_1,\dots,D_k)$ of disjoint finite subsets
$D_i \subset \AP^1$ over which these classes arise.
\end{itemize}
To obtain a single discrete invariant, we write $\nu = (\nu_1,\dots,\nu_k)$ with
 $\nu_i = |D_i|$.  The triple $h = (G,C,\nu)$ is then a {\em Hurwitz parameter}
 in the sense of \cite[\S2]{RV15} or \cite[\S3]{RobHNF}.
 
 A Hurwitz parameter $h$ determines
 a {\em Hurwitz moduli space} $\AHur_h$ whose points $x$ index 
 maps $F_x : \AY_x \rightarrow \AP_t^1$ of type $h$.   
 The Hurwitz space covers the {\em configuration space} $\AConf_\nu$ of 
 all possible divisor tuples $(D_1,\dots,D_k)$ of type $\nu$.   
 The common dimension of $\AHur_h$ and $\AConf_\nu$ is 
 $r = \sum_i {\nu_i}$, the number of critical values of any
 $F_x$.   
  
 Sections~2-4 of \cite{RV15} and Section~3 of \cite{RobHNF} provide
 background on Hurwitz maps, some main points
 being as follows.   
  There is a group-theoretic formula for a mass $\overline{m}$ which is
  an upper bound and often agrees with the degree $m$ of 
  $\pi_h : \AHur_h \rightarrow \AConf_\nu$.  
 If all the $C_i$ are rational classes,
 and under weaker hypotheses as well, 
 then the covering of complex varieties
 descends canonically to a covering of
 varieties defined over the rationals, 
 $\pi_h : \Hur_h \rightarrow \Conf_\nu$.   
 The set $\cP_h$ of primes at which this map has bad 
 reduction is contained in the set $\cP_G$ of primes 
 dividing $|G|$.  
 
 On the computational side, \cite{RobHNF} provides many examples
 of explicit computations of Hurwitz covers.   Because the map $\pi_h$ is 
 equivariant with respect to $\APGL$ actions, we normalize to take representatives
 of orbits and thereby replace $\AHur_h \rightarrow \AConf_\nu$ by 
 a similar cover with three fewer dimensions.  
 Our computations in the previous section for
 $h = (S_7,(322,421,511),(1,1,1))$ and 
 $h = (S_{12},(642,2222211,5322),(1,1,1))$ illustrate
 the case $r=3$.  Computations in the cases $r \geq 4$ proceed 
 quite similarly.  The next section gives some simple 
 examples and a  collection of more complicated examples is given 
 in \S\ref{Cubical}.
 
 Let $\Out(G,C)$ be the subgroup of $\Out(G)$ which
 fixes all classes $C_i$ in $C$.  Then $\Out(G,C)$
 acts freely on $\AHur_h$.  For any subgroup $Q \subseteq \Out(G,C)$,
 we let $\AHur_h^Q$ be the quotient $\AHur_h/Q$ and 
 let $\pi_h^Q : \AHur_h^Q \rightarrow \AConf_\nu$ be
 the corresponding covering map.   One can expect that
 $Q$ will usually play a very elementary role.  For 
 example, it can be somewhat subtle to get the exact degree
 $m$ of a map $\pi_h$, but then the degree of $\pi_h^Q$ is
 just $m/|Q|$.  
 
 We already have the simpler notation $\AHur_h$ for $\AHur_h^{\{e\}}$.  
 Similarly, we use $*$ as a superscript to represent the entire group
 $\Out(G,C)$.    In the literature,
 $\AHur_h$ is often called an inner Hurwitz space
 while $\AHur_h^*$ is an outer Hurwitz space. 
 In the entire sequel of this paper, the only $Q$ that
 we will consider are these two extreme cases.  
 It is important for us to descend to the $*$-level to
 obtain fullness. 
   
\subsection{Belyi pencils} 
\label{BP1} For any $\nu$ as above, the variety $\Conf_\nu$ naturally
comes from a scheme over $\Z$.  Thus for any commutative
ring $R$, we can consider the set $\Conf_\nu(R)$.  
 Section~8 of \cite{RV15} and then the sequel paper \cite{RobHNF} 
 considered $R=\Q$ and its subrings 
 $\Z[1/\cP] = \Z[\{1/p\}_{p \in \cP}]$.  From fibers
 $\AHur_{h,u} \in \Hur_h(\overline{\Q})$ above
 points in $\Conf_\nu(\Z[1/\cP])$ one
 gets interesting number fields, the 
 Hurwitz number fields of the title of \cite{RobHNF}.
 
 Our focus here is similar, but more 
 geometric.  A Belyi pencil $u$ is an 
 algebraic map
\begin{equation}
\label{BelyiPencil}
u : \AP^1_v -\{0,1,\infty\} \rightarrow \AConf_\nu,
\end{equation}
with image not contained in a single $\APGL$ orbit.  
One can think of $v$ as a time-like variable here.  The Belyi pencil 
$u$ then can be understood as giving $r$  points in
$\AP^1_t$, typically moving with $v$.  There are $\nu_i$ points of 
color $i$; points are indistinguishable except for color,
and they always stay distinct except for 
collisions at $v \in \{0,1,\infty\}$. 

To make the similarity clear, for $R$ a ring  
let 
 \[
R\langle v \rangle = R[v,\frac{1}{v(v-1)}].
\]
Then a Belyi pencil can be understood as a point in 
$\Conf_\nu(\C\langle v \rangle)$. 
We say that the Belyi pencil $u$ is {\em rational} if it is in 
$\Conf_\nu(\Q(\langle v \rangle))$.  For rational pencils, one has
a natural bad reduction set $\cP_u$.  It is the smallest set 
$\cP$ with $u \in \Conf_\nu(\Z[1/\cP]\langle v \rangle)$. 

As an example, consider the eight-tuple 
\begin{eqnarray}
\nonumber && \left(
     \left(t^6-8 v t^3 +9 v t^2 -2 v^2\right), \left(t^6-3 t^5+10 v t^3 -15 v t^2 +9 v t -2 v^2\right) , \right. \\
\label{complicatedbelyi} 
 &&  \left.   \left(t^6-6 v t^5 +15 v t^4 -20 v t^3 +6 v^2 t^2 +9 v t^2 -6  v^2 t+v^2\right), \right. \\ 
\nonumber &&  \left. 
     \left(t^4-2 t^3+2 v t -v\right) ,\left(t^4-4 v t +3 v\right) ,\left(2 t^3-3 t^2+v\right) , (t),\{\infty\}
\right).
\end{eqnarray}
The product of the seven irreducible polynomials presented
has leading coefficient $2$ and discriminant $2^{161} 3^{266}  v^{125} (v-1)^{125}$.  Thus 
$u : \AP^1_v-\{0,1,\infty\}  \rightarrow \AConf_{6,6,6,4,4,3,1,1}$ is a Belyi pencil
with bad reduction set $\{2,3\}$.   All the other examples
considered in this paper are substantially simpler.
We include \eqref{complicatedbelyi} to indicate
that Belyi pencils themselves may be quite complicated. 

\subsection{Belyi pencils for $r=4$}
\label{BP2}
Three Belyi pencils play 
a special role in the case $r=4$, and we denote them by $u_{1,1,1,1}$, $u_{2,1,1}$, and $u_{3,1}$.
  For $u_{1,1,1,1}$,
we keep our standard  variable $v$.  To make $u_{2,1,1}$ and $u_{3,1}$ stand
out when they appear in the sequel, we switch
the time-like variable $v$  to respectively $w$ and $j$ for them.  
These special Belyi pencils are then given by
 \begin{align}
 \label{threebasicBelyi}
   \left( \{v\},\{0\},\{1\},\{\infty\} \right), &&
   \left(D_w,\{0\},\{\infty\} \right), &&
   \left( D_j,\{\infty\} \right).
   \end{align}
 Here the divisors $D_w$ and $D_j$ are the root-sets of 
$t^2 + t + \frac{1}{4(1-w)}$ and $ 4(1-j) t^3 + 27 j t + 27 j$ respectively.
So the three Belyi pencils are all rational, and their 
bad reduction sets are respectively $\{\}$, $\{2\}$,
and $\{2,3\}$. 

The images of these Belyi pencils are curves
  \begin{align*}
  \AU_{1,1,1,1} & \subset \AConf_{1,1,1,1},  &
  \AU_{2,1,1} & \subset \AConf_{2,1,1}, & 
  \AU_{3,1} & \subset \AConf_{3,1}.
  \end{align*}
The three curves are familiar as coarse moduli spaces of elliptic curves.  
Here $\AU_{1,1,1,1} = \AY(2)$ parametrizes elliptic curves with a basis of $2$-torsion, 
$\AU_{2,1,1} = \AY_0(2)$ parametrizes elliptic curves with a $2$-torsion point, 
and $\AU_{3,1}=\AY(1)$ is the $j$-line parametrizing elliptic curves.    
The formulas 
\begin{align}
\label{basechange}
w & = \frac{(2v-1)^2}{9}, & 
j = \frac{(3w+1)^3}{(9w-1)^2} =  \frac{2^2(v^2-v+1)^3}{3^3 v^2 (v-1)^2}
\end{align}
give natural maps between these three bases: $\AP^1_v \rightarrow \AP^1_w \rightarrow \AP^1_j
$.

The cases $\nu = (2,2)$ and $\nu = (4)$ are complicated by the presence
of extra automorphisms.  Any configuration $(D_1,D_2) \in \AConf_{2,2}$ 
is in the $\APGL$ orbit of a configuration of the special form $(\{0,\infty\},\{a,1/a\})$.
This latter configuration is stabilized by the automorphism $t \mapsto 1/t$.  
Similarly a configuration $(D_1) \in \AConf_4$ has a Klein four group
of automorphisms.   To treat these two cases, the best approach seems to be
modify the last two pencils in \eqref{threebasicBelyi} to 
$(D_w,\{0,\infty\})$ and $(D_j \cup \{\infty\})$.  Outside of a quick example for $\nu = (4)$ near \eqref{serretheorem}, 
we do not pursue
any explicit examples with such $\nu$ in this paper.

\subsection{Hurwitz-Belyi maps}
\label{HB}
We can now define the objects in our title.  
\begin{Definition}  Let $h = (G,C,\nu)$ be a Hurwitz parameter, 
let $Q$ be a subgroup of $\Out(G,C)$, and let
$\pi^Q_h : \AHur^Q_h \rightarrow \AConf_\nu$ be the corresponding 
Hurwitz map.
Let $u : \AP^1_v \rightarrow  \AConf_\nu$ be a Belyi
pencil.  Let 
\begin{equation}
\pi^Q_{h,u} : \AX^Q_{h,u} \rightarrow \AP_v^1 
\end{equation}
be the Belyi map obtained by pulling back the Hurwitz 
map via the Belyi pencil and canonically completing.
A Belyi map obtainable by this construction
is a {\em Hurwitz-Belyi map}.
\end{Definition}

A terminological explanation is in order. The term {\em Hurwitz-Belyi map} is meant to be parallel
to the term {\em Hurwitz number algebra} that figures 
prominently in \cite{RobHNF}.   Both are constructed
from a Hurwitz parameter and a specialization point,
with $\pi^Q_{h,u}$ having $u \in \Conf_\nu(\C(t))$ and 
$K^Q_{h,u}$ having $u \in \Conf_\nu(\Q)$.     In general, for a Belyi map
$\AX \rightarrow \AP^1$, one has a decomposition
 $\AX = \sqcup_i \AX^i$ into connected components, and one is typically interested
 in the individual Belyi maps $\AX^i \rightarrow \AP^1$.
 So too in a Hurwitz number algebra $K^Q_{h,u} = \prod_i K^{Q,i}_{h,u}$ one
 is typically interested in the factor 
 {\em Hurwitz number fields} $K^{Q,i}_{h,u}$.

A notational convention will be useful as follows.  When $r=4$ and $u$ is
one of the three maps \eqref{threebasicBelyi}, 
then we are essentially not specializing
as we are taking a set of representatives
for the $\APGL$ orbits on $\AConf_\nu$.  
We allow ourselves to drop the $u$ in this
situation, writing e.g.\  $\pi_h : \AX_h \rightarrow \AP^1_j$.

Rationality and bad reduction are both essential to this paper.  If $h$ and $u$ are both defined over $\Q$, then so is $\pi_{h,u}$.  
If $h$ has bad reduction set $\cP_h$ and $u$ has bad reduction
set $\cP_u$ then the bad reduction set of $\pi_{h,u}$ is
within $\cP_h \cup \cP_u$.   All the examples we pursue in this paper, 
with the exception of the example of \S\ref{sixpoint}, satisfy $\cP_u \subseteq \cP_h$.

\section{The two rational Belyi maps as Hurwitz-Belyi maps}  
\label{HurwitzBelyi2}
    This section presents some first examples in the setting $r=4$, 
    and in particular 
interprets the two rational Belyi maps of Section~\ref{Belyi}
 as Hurwitz-Belyi maps.   
 
\subsection{A degree $7$ Hurwitz-Belyi map: 
 computation and dessins}
\label{dessins4} \label{deg7HB}
To realize the Belyi map \eqref{deg7y} as a Hurwitz-Belyi map,
we start from the Hurwitz parameter
\begin{equation}
\label{deg7hur}
h = (S_6,(2_x1111, 3_03, 3_1111,4_\infty 11),(1,1,1,1)).
\end{equation}
Here and in the sequel we often present
Hurwitz parameters with subscripts which indicate 
our normalization, without being as formal as we were in 
Section~\ref{Belyi}.  The marked Hurwitz parameter \eqref{deg7hur}
tells us to consider rational functions 
of the form 
\begin{equation}
\label{g7}
F(y) = \frac{(y^2-a)^3 (b+c+1)}{(1-a)^3 \left(y^2+b y+c\right)}
\end{equation}
and the equation
\begin{equation}
\label{diff7}
4 y^3 + 5 b y^2 + 2 a y + a b = 4 (y-1)^2 (y-x).
\end{equation}
The left side of \eqref{diff7} is a factor of the numerator of $F'(y)$ and thus
its roots are critical points.  The right side gives the required locations and multiplicities
of these critical points.    

 \begin{figure}[bht]
\includegraphics[width=4.4in]{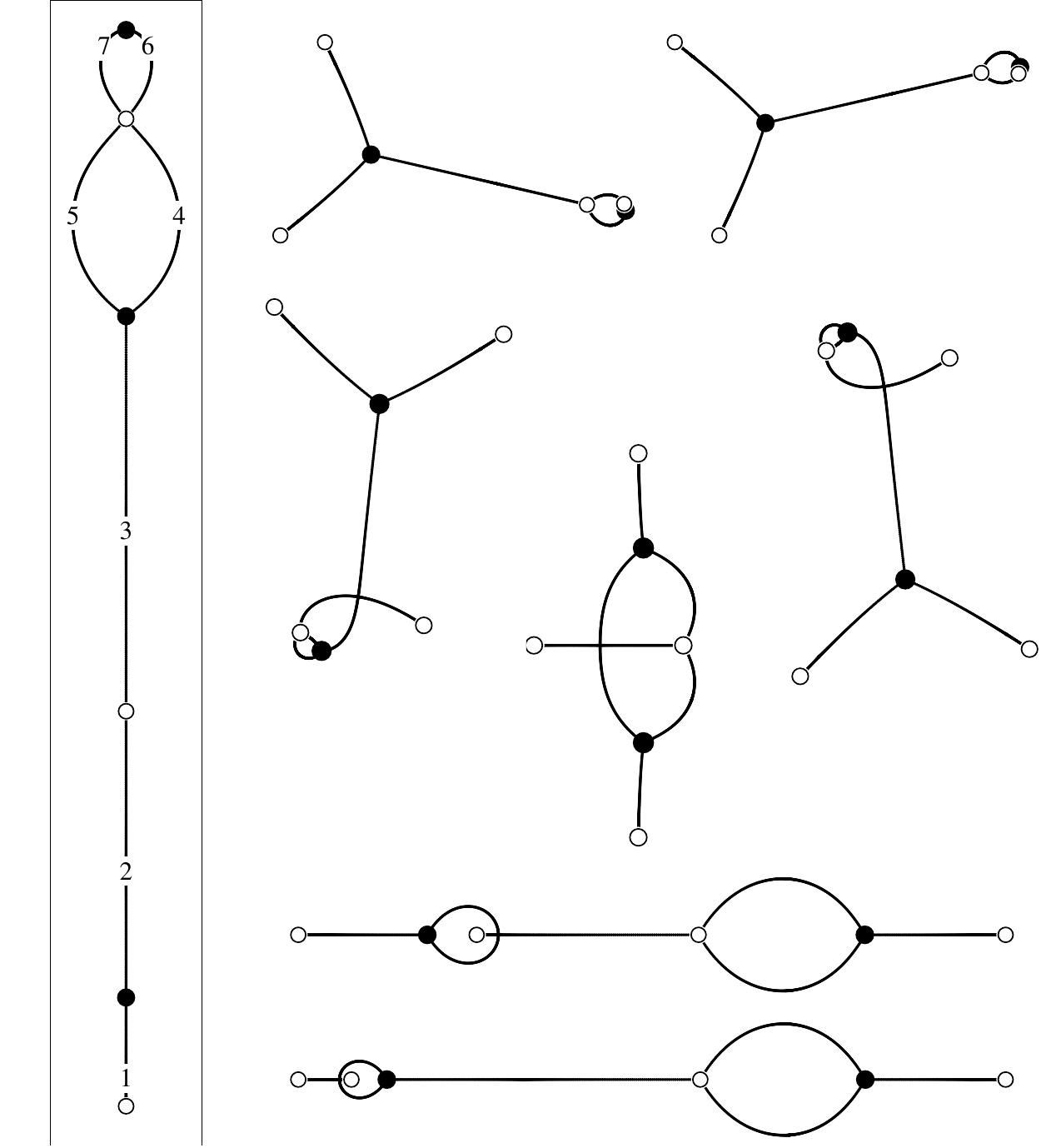}
\caption{\label{deg7nice} Left: The dessin $\AX_{h,[0,1]} \subset \AP^1_x$ belonging to $h = (S_6, \, (21111,33,3111,411), \, (1,1,1,1))$, with real axis 
pointing up; the integers $i$ are at preimages $x_i \in \AX_{h,[0,1]}$ of $1/2 \in \AP^1_v$.  
  Right:  the dessins $\AY_{h,x_i,[0,1]} \subset \AP^1_y$
 for the seven $i$.} 
\end{figure}

Equating coefficients of powers of $y$ in \eqref{diff7} and solving,
we get
\begin{align}
a & = \frac{5 x}{x+2}, &
b & = -\frac{4}{5} (x+2), &
c & = \frac{4 x^2+5 x+4}{3(x+2)}.
\end{align}
Summarizing, we have realized  
$\AX_h$ as the complex projective $x$-line and identified each $\AY_{h,x}$
with the complex projective $y$-line $\AP^1_y$ so that the
covering maps $F_{h,x} : \AP^1_y \rightarrow \AP_t^{1}$ 
become
\begin{equation}
\label{deg7source}
F_{h,x}(y) = \frac{\left((x +2) y^2-5 x\right)^3}{4 (2 x-1)
 \left(-15 (x+2) y^2+12 (x+2)^2 y -5 (4 x^2+5 x+4)
    \right)
    }.
\end{equation}
Since the fourth critical value of $F_{h,x}$ is just $F_{h,x}(x)$, we 
have also coordinatized the 
 Hurwitz-Belyi map   $\pi_h: \AX_h \rightarrow \AP^1_v$ to 
\begin{equation}
\label{deg7x}
\pi_h(x) =  -\frac{x^3 \left(x^2+2 x-5\right)^2}{4 (2 x-1) (3 x-4)}.
\end{equation}
  Said more explicitly, to pass from the right side of \eqref{deg7source}
 to  the right side \eqref{deg7x}, one substitutes $x$ for $y$ and 
 cancels $x^2+2x-5$ from top and bottom.

 The rational function \eqref{deg7x} appeared
 already as \eqref{deg7y}, with its dessin
 printed in the upper left of Figure~\ref{deg7four}.
 This connection is our first explanation of why 
 \eqref{split3plus1} splits.  It also explains
 why the bad reduction set of the rational
 Belyi map is in $\{2,3,5\}$.

 Figure~\ref{deg7nice} presents the current situation pictorially, with 
 $v = 1/2$ chosen as a base point.   The elements of 
 $\AX_{h,1/2} = \pi_h^{-1}(1/2)$ are labeled in the box at the left, where 
 the real axis runs from bottom to top for a better overall
 picture.  For each $x \in \AX_{h,1/2}$, a corresponding
  dessin $\AY_{x,[0,1]}$ is drawn to its right.  Like
  the standard dessins of \S\ref{dessins}, these dessins
  have black vertices, white vertices, and faces.  
  However they each also have five vertices of a fourth type
  which we are not marking, corresponding to the five parts of the
  partition $21111$.  The valence of this type of vertex
  with ramification number $e$ is $2e$, so only
  the extra vertex coming from the critical
  point with $e=2$ is visible on Figure~\ref{deg7nice}.  The
  action of the braid group to be discussed in \S\ref{braidbackground} can
  be calculated geometrically from these dessins.
  
\subsection{Cross-parameter agreement}  
\label{crossparam} An interesting
phenomenon that we will see repeatedly is 
{\em cross-parameter agreement}, discussed also in \cite[\S3]{RobHNF}.  
This phenomenon
occurs when two different Hurwitz parameters
give rise to isomorphic Hurwitz
covers.    Already this phenomenon occurs
 for our septic Belyi map, which we realize in a second 
way as a Hurwitz-Belyi map as follows.   

   For the normalized Hurwitz parameter 
\[
\hat{h}= (S_5, ( 2_z111, 221, 3_111 , 3_\infty 2_0), (1,1,1,1)),
\]
 the computation is easier than it was for  the Hurwitz parameter $h$ of \eqref{deg7hur}. 
 An initial form of $\hat{F}(y)$, analogous to 
 \eqref{g7}, is 
 \begin{equation}
 \hat{F}(y)  =  \frac{(y-c) \left(y^2+a y+b\right)^2}{(1-c) y^2 (a+b+1)^2}.
 \end{equation}
 Analogously to \eqref{deg7source}, the covering maps $\AP^1_{y} \rightarrow \AP^1_t$ 
 are
 \begin{equation}
\hat{F}_z(y)  =  \frac{(4 y z+2 y-z) \left(-2 y^2 z-y^2+6 y z^2+14 y z+6 y+12 z^2+12
    z+3\right)^2}{4 y^2 (3 z+2)^5}.
 \end{equation}
 Analogously to \eqref{deg7} the Hurwitz-Belyi map $\AP^1_z \rightarrow \AP^1_v$ is
 \begin{equation}
 \label{deg7alt}
 \hat{\pi}(z) = \hat{F}_z(z) = 
 \frac{\left(4 z^2+z\right) \left(4 z^3+25 z^2+18 z+3\right)^2}{4 z^2 (3
    z+2)^5}.
 \end{equation}
 The map \eqref{deg7alt} agrees with the map \eqref{deg7}
  via the substitution 
  $z = (1-2x)/(3x-4)$.
In terms of Figure~\ref{deg7nice}, 
the dessin at the left remains exactly the same, up to change of coordinates.
 In contrast, the seven sextic
 dessins at the right would be each replaced by a corresponding quintic dessin.

\subsection{Two degree $12$ Hurwitz-Belyi maps}
\label{ambiguity}  In our labeling of conjugacy classes mentioned in \S\ref{notation},
we are accounting for the fact that the five cycles
in $A_5$ fall into two classes, with representatives $(1,2,3,4,5) \in 5a$ and
$(1,2,3,4,5)^2 = (1,3,5,2,4) \in 5b$.  Consider two Hurwitz parameters, as in the left column:
\begin{align*}
h_{aa} & =  (A_5, \,  (5a, \;  311, 221),  \; (2,1,1)), &  (\beta_0,\beta_1,\beta_\infty) & = (5331,222222,5322), \\
h_{ab} & = (A_5, \, (5a,\; 5b, \; 311, 221),(1,1,1,1)), &  (\beta_0,\beta_1,\beta_\infty) & = (642,2222211,5322).
\end{align*}
Applying the outer involution of $A_5$ turns $h_{aa}$ and $h_{ab}$ respectively
into similar Hurwitz parameters $h_{bb}$ and $h_{ba}$, and so it would be redundant to 
explicitly consider these latter two.     

    It is hard to computationally distinguish $5a$ from $5b$.  We will deal with this 
 problem by treating $h_{aa}$ and $h_{bb}$ simultaneously.
Thus we working formally with 
\begin{align*}
h = (S_5, (5, 311, 221),(2,1,1)),
\end{align*}
ignoring that the classes do not generate $S_5$.  A second problem is that there
are only eight parts altogether in the partitions $5$, $5$, $311$, and $221$, so 
the covering curves $\AY$ have genus two.  

    To circumvent the genus two problem, we use the braid-triple method, 
as described later in Section~\ref{4vs3}. The mass formula \cite[\S3.5]{RobHNF} applies to $h$, with only the two abelian characters 
of $S_5$ contributing. It says that
the corresponding cover 
 $\AX_h \rightarrow \AP_w^1$ has degree
\[
\frac{|C_5|^2 |C_{311}| |C_{221}|}{|A_5|^2} = \frac{24^2 \cdot 20 \cdot 15}{120 \cdot 60} = 24. 
\]
A braid group computation of the type described in \S\ref{stepone}
 says $\AX_h$ has two components, each of degree $12$.  
The braid partition triples are given in the right column above.   
The $\beta_\tau$ then enter the formalism of 
Section~\ref{Belyi} as the $\lambda_\tau$ there.  Conveniently.
each cover sought has genus zero, and so 
the covers are easily computed.  
The resulting polynomials are
\begin{eqnarray*}
f_{12aa}(w,x) & \!\! = \!\! &   x^5\left(9 x^2-21 x+16\right)^3 (x+3) - 2^8 w (x-1)^3 \left(9 x^2-12  x+8\right)^2, \\
f_{12ab}(w,x) & \!\! = \!\! & 5^5 (x-1)^4 x^6 (5 x+4)^2  -2^4 3^3  w(2 x+1)^3 \left(5 x^2-6 x+2\right)^2.
\end{eqnarray*}
Up to the simultaneous letter change $y \leftrightarrow x$, $t \leftrightarrow w$,  
the equation $f_{12ab}(w,x)=0$ defines the exact same map as \eqref{deg12y}.
The current context explains why this map is defined over $\Q$ and
has bad reduction at exactly $\{2,3,5\}$.   

 \begin{figure}[htb]
\includegraphics[width=4.6in]{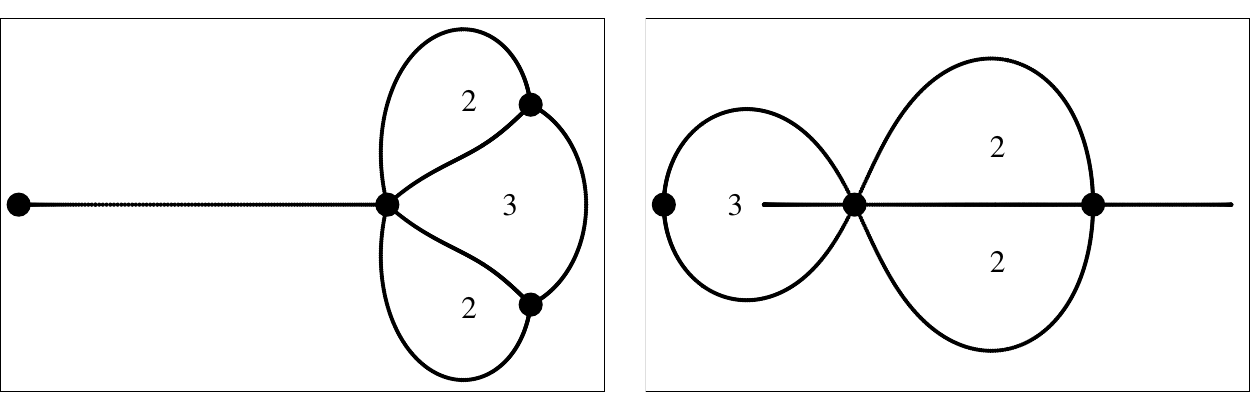}
\caption{\label{two12}  Rational degree $12$ dessins coming from Hurwitz parameters involving irrational classes.
$h_{12aa}$ yielding partition triple $(5331,222222,5322)$ on the left and $h_{12ab}$ yielding $(642,2222211,5322)$ on the right}
\end{figure}

Dessins corresponding to $f_{12aa}$ and $f_{12ab}$ are drawn in Figure~\ref{two12}.  
The two dessins 
present an interesting contrast:  the dessin 
on the left of Figure~\ref{two12} is the unique dessin with its partition triple, while the 
dessin on the right is one of the $24$ locally equivalent dessins drawn in Figure~\ref{picts24}.

\subsection{An $M_{12}$ specialization} 
\label{m12} Specializing Hurwitz-Belyi maps yields
interesting number fields.  Except for this subsubsection, we are saying nothing about
this application, because we discuss specialization of Hurwitz covers 
quite thoroughly in \cite{RV15},  
\cite{Rob15}, and  \cite{RobHNF},  However the cover given by $f_{12ab}(w,x) = 0$ yields 
a particularly interesting number field and so we discuss specialization 
for this cover. 

In general, let $f(w,x) = 0$ be a polynomial defining a degree $m$ Belyi map $\AX \rightarrow \AP_w^1$ with
ramification partition above $1$ having the form $2^a 1^{m-2a}$.  Suppose the Belyi map
has bad reduction within $\{2,3,5\}$.  Then specializing $w$ to any rational number 
in the $183$-element set $\SU_{2,1,1}(\Z[1/30])-\{1\}$ of \cite[\S4]{Rob15} gives a number algebra ramified within $\{2,3,5\}$.  

The principles of \cite[\S4]{RobHNF} give different ways that we
expect specialization to behave generically.  Like all the covers of this paper,
the one given by $f_{12ab}(w,x)=0$ conforms very well to these principles.  
Principle A has no exceptions: the $183$ algebras $\Q[x]/f_{12ab}(w,x)$ are
all non-isomorphic.   Principle B has two exceptions: $w=-625/216$ and $w=1/4$
each give Galois groups strictly smaller than $A_{12}$.  

The first exceptional specialization can be compared with complete tables of number fields \cite{JR14}.  
The polynomial $f_{12ab}(-625/216,x)$ factors into two sextics each
with Galois group $S_6$.  For one, the splitting field is the third least ramified
all of $S_6$ Galois fields, having root discriminant  
$2^{31/12} 5^{23/20} \approx 38.15$.   

For the second exceptional specialization, one can only compare 
with incomplete lists.  The polynomial $f_{12ab}(1/4,x)$,
or equivalently
\[
x^{12}-24 x^{10}+180 x^8-60 x^6-2520 x^5+4320 x^4-2520 x^3+864 x^2-216,
\]
has the Mathieu group $M_{12}$ as its Galois group.  The splitting 
field has root discriminant $2^{3/2} 3^{3/2} 5^{23/20} \approx 93.55$.
Comparing with the discussion in \cite[\S6.2]{Rob16}, one sees that this
is currently the second least ramified of known  $M_{12}$ Galois fields, 
being just slightly greater than the current minimum 
$2^{25/12} 3^{10/11} 5^{13/10} \approx 93.23$.

\section{The semicubical clan}
\label{Cubical}
    An interesting phenomenon is that certain infinite collections of 
 Hurwitz-Belyi maps can be studied uniformly by means of 
 parameters.  We illustrate this phenomenon by 
 means of a three-parameter clan in the setting $r=4$.
 This clan serves as a general model, because it 
 exhibits behavior we have seen in many other clans. 
 We often think of a single Hurwitz-Belyi map as a family of number fields,
 as just illustrated by \S\ref{m12}.    
 Accordingly, we are using the word ``clan'' to indicate a collection larger than
 a family.
  
    This section has a different tone than the other sections,
  as we are studying infinitely many maps at once.  We treat enough
  topics to give a sense of how we expect clans  to behave in general.    In the
  last subsection we explain how the clan approach responds
  interestingly but not optimally to our inverse problem.  
  In the remaining sections, we will therefore return to studying
  Hurwitz-Belyi maps one by one.  
 
 \subsection{Couveignes' cubical clan} 
 Our clan is very closely related to a four-parameter 
 clan studied by Couveignes \cite{Cou97}. 
    To place our results
 in their natural context, we first define
  Couveignes' clan, using notation adapted to our context, and present
 some of his results.   
 
 For Couveignes' clan, take $a$, $b$, $c$, and $d$ distinct positive 
 integers and set $n=a+b+c+d$.  The Hurwitz
 parameters are
 \[
 h(a,b,c,d) = (S_n, \;  (2 \, 1^{n-2}, \; a \, b \, c \, d,  \;  3 \, 1^{n-3}, \;  n), \; (1, \, 1, \, 1, \, 1)).
 \]
  The degree of the Hurwitz-Belyi map 
  $\pi_{a,b,c,d} : \AX_{a,b,c,d} \rightarrow \AP^1_v$ is $m = 6n$.  The
  braid partition triple $(\beta_0,\beta_1,\beta_\infty)$ is deducible 
  from the figure in \cite[\S5.2]{Cou97}:
 \begin{eqnarray}
\nonumber \beta_0 &= &  (a+b)^2 (a+c)^2 (a+d)^2 (b+c)^2 (b+d)^2 (c+d)^2, \\
\label{Cbraid}  \beta_1 & = & 4^6 1^{m-24}, \\
 \nonumber  \beta_\infty & = &(a+b+c)^2 (a+b+d)^2 (a+c+d)^2 (b+c+d)^2.
 \end{eqnarray}
 The total number of parts is $12+(6+m-24)+8=m+2$,
 and so $\AX_{a,b,c,d}$ has genus zero.  
 
\begin{figure}[bht]
\includegraphics[width=3in]{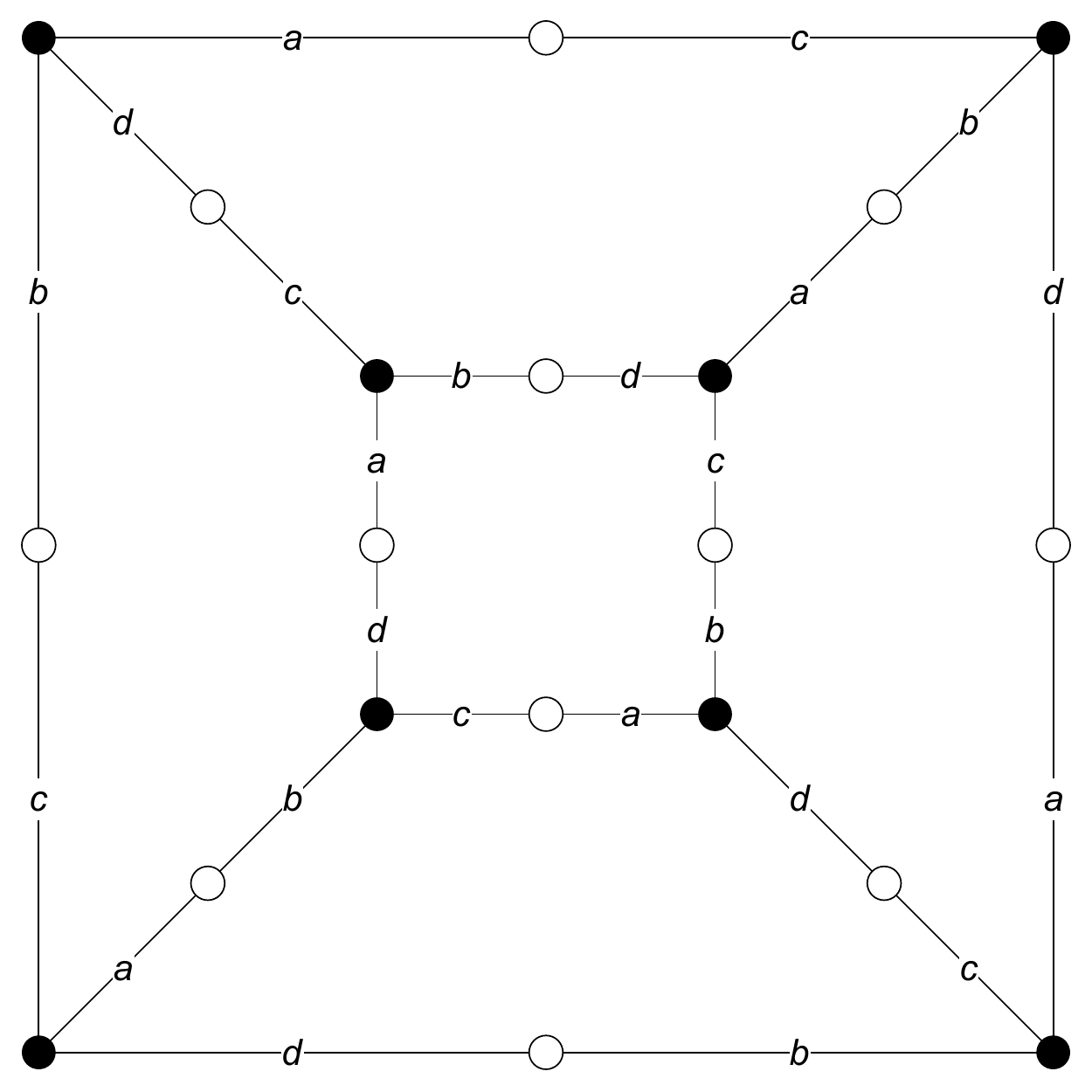}
\caption{\label{CouveignesCube} Schematic indication of Couveignes' dessin
with parameters $(a,b,c,d)$ based on the combinatorics of a cube.  The
actual dessin is obtained by replacing each 
$\bullet  \! \! - \!\!\!-  u \! - \! \! \!  - \! \!  \circ$
by $u$ parallel edges.}
\end{figure}

 In the current context, it is better to modify the visualization conventions of \S\ref{dessins}, to exploit
 that the number of parts in $\beta_0$ and $\beta_\infty$ is small and independent of the parameters.  Accordingly,
 we now view the interval $[-\infty,0]$ in the projective line $\AP^1_v$  as the
 simple bipartite graph $\bullet \!\! - \!\!\! - \!\!\! - \!\! \circ$.  
 The dessin $\pi^{-1}_{a,b,c,d}([-\infty,0]) \subset \AX_{a,b,c,d}$, capturing
 Couveignes' determination \cite[\S5.2 and \S9]{Cou97} of the permutation triple $(b_0,b_1,b_\infty)$ underlying
 $(\beta_0,\beta_1,\beta_\infty)$,
  is indicated
 schematically by Figure~\ref{CouveignesCube}.   Note 
 that our visualization is dual to that of Couveignes,
 as our dessin is formatted on a cube rather than
 an octahedron.

The $\Q$-curve $\SX$ underlying $\AX = \SX(\C)$ is naturally given 
 in $\SP^3$ by the following
symmetric equations \cite[\S5.1]{Cou97}:
\begin{eqnarray*}
a x_1 + b x_2 + c x_3 + d x_4 & = & 0, \\
a x_1^2 + b x_2^2 + c x_3^2 + d x_4^2 & = & 0. 
\end{eqnarray*}
The second equation has no solution besides $(0,0,0,0)$ and so $\SX(\R)$ is empty.  
This non-splitting of $\SX$ over $\R$, which forces non-splitting over
$\Q_p$ for an odd number of $p$, is one of the main focal points 
of \cite{Cou97}.

 Note finally that our requirement that 
$a$, $b$, $c$, and $d$ are all distinct is just so that the
above considerations fit immediately into our formalism.  
One actually has natural covers $\SX_{a,b,c,d} \rightarrow \SP^1$
of degree $6n$
even when this requirement is dropped.  However these covers have extra symmetries and can never be full,
as illustrated by the rotation $\iota$ discussed in the next subsection.       

\subsection{The semicubical clan and explicit equations}
Couveignes does not give explicit equations for the map 
$\SX_{a,b,c,d} \rightarrow \SP_v^1$.  Such equations
could not be given in our simple standard form because, as just
discussed, $\SX_{a,b,c,d}$ is not isomorphic to $\SP^1$ over $\Q$.  
The fact that all multiplicities  are even in the triple \eqref{Cbraid} is
necessary for this somewhat rare obstruction.  

For our semicubical clan, we still require that $a$, $b$, and $d$
are distinct.  But now we essentially set $c=d$ in Couveignes'
situation, so that the degree takes the asymmetric form $n = a + b + 2d$.  
We thus are now considering the Hurwitz parameters
\begin{equation}
\label{semicubical}
h(a,b,d) = (S_n, (2 \, 1^{n-2}, \; a_0 \, b_1 \, d^2, \; 3_x \, 1^{n-3}, \; n_\infty), (1, \, 1, \, 1, \, 1)).
\end{equation}
 The four subscripts are present as usual for the purposes of normalization 
 and coordinatization.  They will enter into our proof of Theorem~\ref{clantheorem}
 below.
 
The fact that two distinguishable points have now become 
indistinguishable implies that $\AX_{a,b,d} = \AX_{a,b,d,d}/\iota$,
where $\iota$ is the rotation interchanging $c$ and $d$ in
Figure~\ref{CouveignesCube}.  The fixed points of this rotation are the upper-left 
and lower-right white vertices, each with valence $c+d=2d$.   Thus 
$\AX_{a,b,d}$ has degree $3 n$ over $\AP^1_v$.    
The new braid partition triple is
\begin{eqnarray}
\nonumber \beta_0 & = & (a+b)_\infty (a+d)^2 (b+d)^2 d^2, \\
\label{Cbraid2} \beta_1 & = & 4^3 1^{3 n-12}, \\
\nonumber \beta_\infty & = & (a+b+d)^2 (a+2d)_0(b+2d)_1. 
\end{eqnarray}
There are three singletons, namely the parts subscripted
$0$, $1$, and $\infty$.  So not only is the 
$\Q$-curve $\SX_{a,b,d}$ split, but also our choice of
subscripts gives it a canonical coordinate.

To compute $\pi_{a,b,d}$ as an explicit rational function,
we follow our standard procedure.   Since this particular
instance of our procedure is done with parameters, we display the
result as a theorem and give some details of the calculation.

\begin{Theorem} \label{clantheorem}   For distinct positive integers $a$, $b$, $d$, 
let $n = a+b+2d$ and 
\begin{eqnarray*}
A & = & -2 n x (a+d)+(a+d) (a+2 d)+n x^2 (n-d), \\
B & = & a (a+d)-2 a n x+n x^2 (-(d-n)), \\
C & = & x^2 (a+b) (n-d)-2 a x (n-d)+a (a+d), \\
D & = & n x^2 (a+b)+a (a+2 d)-2 a n x. 
\end{eqnarray*}
Then the Hurwitz-Belyi map for \eqref{semicubical} is
\begin{equation}
\label{claneq}
\pi_{a,b,d}(x) = \frac{a^a b^b A^{a+d} B^{b+d} D^d}{2^d d^{2d} n^n x^{a+2d} (1-x)^{b+2d} C^{n-d}}.
\end{equation}
\end{Theorem}

\proof The polynomial
\begin{equation}
\label{semicubicalF}
F_x(y) = \frac{ y^a (y-1)^b (y^2 + r y + s)^d}{ x^a (x-1)^b (x^2 + r x + s)^d}
\end{equation}
partially conforms to \eqref{semicubical}, including that $F_x(x) = 1$.   From the
$3_x$ we need also that $F_x'(x)=0$ and $F_x''(x)=0$.  The 
derivative condition is satisfied exactly when 
\[
r=\frac{-n x^3 +(a+2 d) x^2 - (a +b)s  x+a s}{((a+b+d)x^2-(a+d) x)}.
\]
The second derivative condition is satisfied exactly when
\[
s = \frac{n  (a+b+d) x^4-2 n (a+d) x^3+ (a+d) (a+2 d) x^2}{ (a+b) (a+b+d) x^2-2 a  (a+b+d)x+a (a+d)}.
\]
The identification of $r$ and $s$ completely determine the maps
$F_{x} : \AP^1_y \rightarrow \AP^1_t$.

From a linear factor in the numerator of $F'_{x}(y)$, one gets that
the critical point corresponding to the $2$ in the first class $2 \, 1^{n-2}$ in \eqref{semicubical} is 
\[
y_{x} =
 \frac{a s}{n x^2}.
\]
Substantially simplifying $F_x(y_x)$ gives the right side of \eqref{claneq}. 
 \qed
 
 From the form of the normalized partition tuple \eqref{Cbraid2},
one knows {\em a priori} that 
\begin{equation}
\label{sym1}
\pi_{a,b,d}(1-x) = \pi_{b,a,d}(x).
\end{equation}  
 Indeed, one can check that the simultaneous interchange
 $a \leftrightarrow b$, $x \leftrightarrow 1-x$ interchanges
 $A$ and $B$ and fixes $C$ and $D$.  Given this fact, 
 the symmetry \eqref{sym1} is visible in the main formula 
 \eqref{claneq}.  
 
\subsection{Dessins}
\begin{figure}[htb]
\includegraphics[width=4.9in]{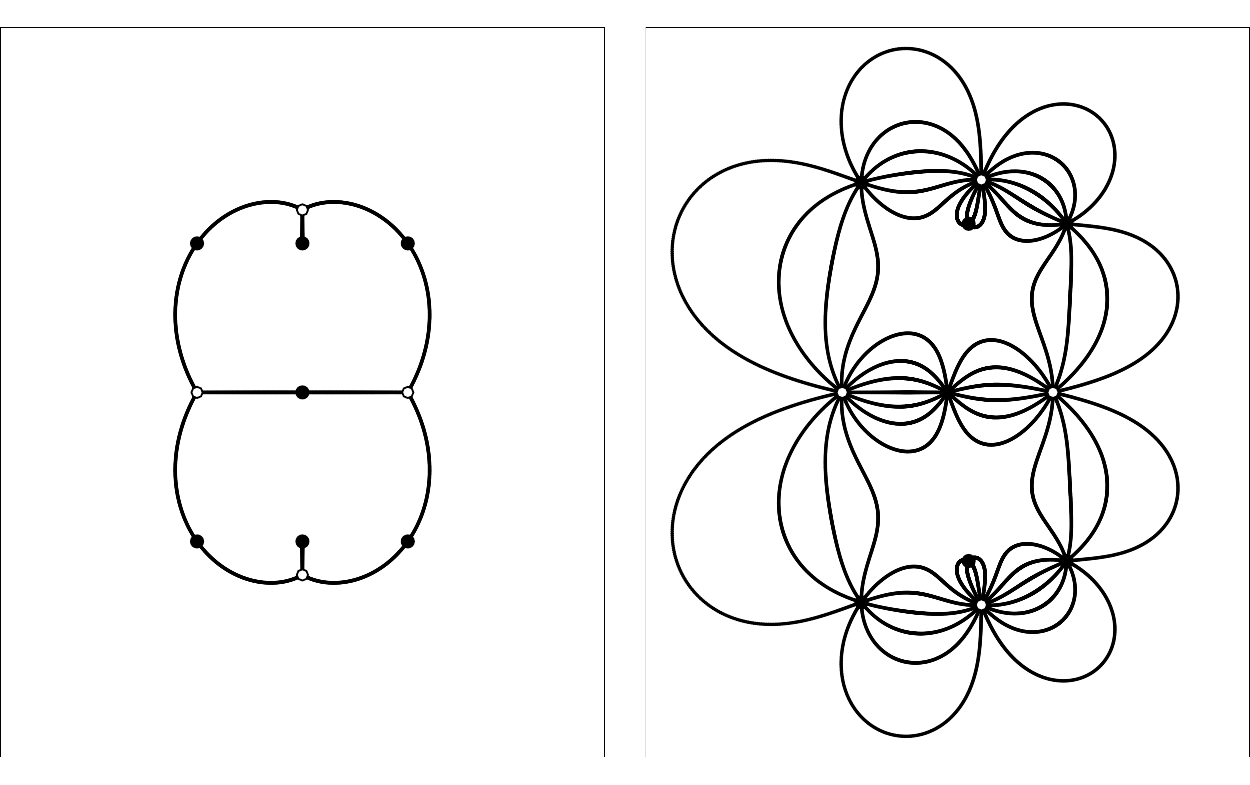}
\caption{\label{couvpos} Left: $\gamma(1,1,1)$.  Right: $\gamma(7,6,4)$. }
\end{figure}

As with Couveignes' cubical clan, the semicubical clan gives covers $\pi_{a,b,d}$, 
even when $a$, $b$, and $d$ are not required to be distinct. 
For example, the simplest case is the dodecic cover
\[
\pi_{1,1,1}(x) = 
 -\frac{\left(6 x^2-8 x+3\right)^2  \left(8 x^2-8 x+3\right) \left(6 x^2-4 x+1\right)^2}{2^8 x^3 (x-1)^3
     \left(3 x^2-3 x+1\right)^3}.
\]
An example representing the main case of distinct parameters is 
\[
\pi_{7,6,4}(x) =  \frac{\left(119
    x^2-154 x+55\right)^{11} \left(13 x^2-14 x+5\right)^4 \left(51 x^2-42 x+11\right)^{10}}{2^{14} x^{15} (x-1)^{14}\left(221 x^2-238
    x+77\right)^{17}}.
\]
Let $\gamma(a,b,d) =  \pi_{a,b,d}^{-1}([-\infty,0]) \subset \AP^1_x$.  The left part of the Figure~\ref{couvpos} is a view on $\gamma(1,1,1)$.
 To obtain the general $\gamma(a,b,d)$ topologically, one replaces each 
 segment of $\gamma(1,1,1)$
 by the appropriate number $a$, $b$, or $d$ of  parallel  segments, so as to 
 create $m= 3n = 3a+3b+6d$ edges in total.   As an example, the right part of the figure draws
 the degree sixty-three dessin $\gamma(7,6,4)$.
 
The view on $\gamma(a,b,d)$ given by \eqref{couvpos} was obtained via the involution $s(x) = x/(2x-1)$ which fixes
$0$ and $1$ and interchanges $\infty$ and $1/2$.    The sequence of white, black, and white vertices on the
real axis are the points $0$, $\infty$, and $1$ in the Riemann sphere.  Thus the point not in the plane
of the paper is the point $x=1/2$.  The involution \eqref{sym1} corresponds to rotating the figure one half-turn
about its central point $\infty$.  

\subsection{Three imprimitive cases} Given a Belyi map $\AX \rightarrow \AP^1$ with
connected $\AX$, a natural first question is whether it has strictly intermediate 
covers.  As preparation for the next subsection, we exhibit three settings
where there is such an intermediate cover
\begin{equation}
\label{intermediate}
\AX_{a,b,d} \stackrel{\delta}{\rightarrow} \AY \stackrel{\epsilon}{\rightarrow} \AP_v^1.
\end{equation}
Another way to view the general question is that connectivity of $\AX$ is exactly
equivalent to transitivity of the monodromy group $M = \langle g_0,g_1 \rangle$.  
The cover $\AX \rightarrow \AP^1$ has imprimitive monodromy group exactly
if there exists a $\AY$ as in \eqref{intermediate} and primitive otherwise.  

\subsubsection*{Case 1.} Let $e = \gcd(a,b,d)$.  If $e>1$ then the explicit formula
\eqref{claneq} says that 
\begin{equation}
\label{homogeneity}
\pi_{a,b,d}(x) = \pi_{a/e,b/e,d/e}(x)^e.
\end{equation}  Thus
one has imprimitivity here, with the cover $\epsilon$ naturally coordinatized
to $y \mapsto y^e$.

\subsubsection*{Case 2.}  Suppose $a=b$.  As a special case of \eqref{sym1}, the cover $\pi_{a,a,d}$ has
the automorphism $x \mapsto 1-x$, corresponding to rotating 
dessins as in Figure~\ref{couvpos}.  To coordinatize $\AY$, we introduce the
function $y = \delta(x) = x(1-x)$.  Then $\epsilon(y)$ works out to 
\begin{equation}
\label{formulaV}
\pi_{a,d}^V(y) =
 \frac{(4 a y-a+4 d y-2 d)^d \left(4 y^2 (2 a+d)^2-4 a y (2 a+3 d)+a (a+2 d) \right)^{a+d}}{d^{2 d} 2^{2 a+3 d} y^{a+2 d} (4 a y-a+2 d y-d)^{2 a+d}}.
\end{equation}
The superscript $V$ indicates that $\pi^V_{a,d}$ comes from $\pi_{a,a,d,d}$ by quotienting by a noncyclic group of 
order four.   The ramification partitions and the induced normalization of $\pi_{a,d}^V$ are
\begin{eqnarray}
\nonumber \alpha_0 & = &  a_\infty \; d \,  (a+d)^2, \\
\label{ramtriple2} \alpha_1 & = &  4 \, 2_4 \, 1^{1.5 n-6},\\
\nonumber \alpha_\infty & = &  (a+2d)_0, \, (2a+d).
\end{eqnarray}
The covers $\pi_{a,d}^V$ and $\pi_{d,a}^V$ are isomorphic, although our
choice of normalization obscures this symmetry.  

\subsubsection*{Case 3.} Suppose $d \in \{a,b\}$.  Via \eqref{sym1}, it suffices to consider the case $b=d$.  
Then while the cover $\pi_{a,d,d}$ does
not have any automorphisms, the original cubical cover $\pi_{a,d,d,d}$ has 
automorphism group $S_3$,  This implies that $\pi_{a,d,d}$ has a subcover
$\epsilon = \pi_{a,d}^S$ of index three.    To coordinatize $\AY$ in this case, we use the
function
\[
y=\delta(x) = \frac{(x-1)^3 (a+d) (a+2 d)}{x \left(x^2 (a+d) (a+2 d)-2 a x (a+2 d)+a (a+d)\right)}.
\]
Then     
\begin{equation}
\label{formulaS}
\pi_{a,d}^S(y) = \frac{(-a)^a (y-1)^{a+d} \left(a^3 y^2 (a+d)+2 a^2 d y (5 a+9 d)+27 d^2 (a+d) (a+2
    d)\right)^d}{2^d d^d n^n y^d}.
 \end{equation}
 So the ramification partitions and the induced normalization of $\pi_{a,d}^S$ are 
 \begin{eqnarray}
\nonumber \alpha_0 & = &  (a+d)_1 \, d^2  \\
\label{ramtriple3} \alpha_1 & = &  4 \, 1^{n-4},\\
\nonumber \alpha_\infty & = &  (a+2d)_\infty  \, d_0.
\end{eqnarray}

\subsection{Primitivity and fullness}  
\label{primandfull} The next theorem says in particular
that all $\pi_{a,b,d}$ not falling into Cases 1-3 of the previous subsection
have primitive monodromy.   

\begin{Theorem} \label{monodromy}
 Let $a$, $b$, $d$ be distinct positive integers with $\gcd(a,b,d)=1$.   
 Let $\pi : \AX \rightarrow \AP^1$ be a Belyi map with the same
 ramification partition triple \eqref{Cbraid2} as $\pi_{a,b,d}$.  Then
 $\pi$ has primitive monodromy.  
\end{Theorem}

\noindent Note that the hypothesis $\gcd(a,b,d)=1$ excludes Case 1 from the previous subsection.  The distinctness
hypothesis excludes Cases 2 and 3.  Since these cases all have imprimitive monodromy,
we will have to use these hypotheses.

\proof The hypothesis $\gcd(a,b,d)=1$ has implications on the parts of $\beta_0$ 
and $\beta_\infty$ in \eqref{Cbraid2}.  For $\beta_0$ it implies $\gcd(a+b,a+d,b+d,d)=1$
while for $\beta_\infty$ it implies $\gcd(a+b+d,a+2d,b+2d) \in \{1,3\}$.  
 As before,
let $n = a+b+2d$ and $m = 3n$. 
The two hypotheses together say that the smallest possible 
degree $m$ is twenty-one, coming from $(a,b,d)=(3,2,1)$.  

Let $\AY$ be a strictly intermediate cover as in \eqref{intermediate}, 
with $\AX_{a,b,d}$ replaced by $\AX$.    Let $e$ be the degree of $\AY \rightarrow \AP^1$ and 
  write its ramification partition triple as $(\alpha_0,\alpha_1,\alpha_\infty)$.    
Because  $\gcd(a+b,a+d,b+d,d)=1$, the cover $\AY$ cannot be totally ramified
over $0$.   Because $\gcd(a+b+d,a+2d,b+2d) \in \{1,3\}$, it can be totally ramified
over infinity only if $e = 3$.    There are then two possibilities, as
$(\alpha_0,\alpha_1,\alpha_\infty)$ could be $((1,1,1),(3),(3))$ or $((2,1),(2,1),(3))$.  
The $\alpha_1=3$ in the first possibility immediately contradicts $\beta_1=(4^3,1^{m-{12}})$.  
The $\alpha_1 = (2,1)$ in the second possibility allows 
two possible forms for $\beta_1$, namely  $(4^3,1^6)$ and $(4^3)$. 
But both of these have degree less than twenty-one.  So $e=3$ is not possible,
and thus $\AY$ cannot be totally ramified over $\infty$ either.

Since $\alpha_0$ and $\alpha_\infty$ both have at least
two parts, the minimally ramified partition $(2,1^{e-2})$
is eliminated as a possibility for $\alpha_1$, by the
Riemann-Hurwitz formula.     
The candidates $(2^2,1^{e-4})$ and $(2^3,1^{e-6})$ for $\alpha_1$
both force  $\AX$ to be a double cover of $\AY$, so that $e=m/2$;
but both are then incompatible with $\beta_1 = (4^3,1^{m-12})$.  
This leaves  $(4,2,1^{m/2-6})$ and $(4,1^{m/3-4})$
 as the only possibilities for $\alpha_1$.    In the first case,
 the two critical values of the double cover 
$\AX \rightarrow \AY$ would have to correspond
to the $2$ in $\alpha_1$ and the image of the singleton 
$a+b$ of $\beta_0$; the parts of $\beta_\infty$ 
would have to be those of $\alpha_\infty$ with multiplicities
doubled; from the form of $\beta_\infty$ in \eqref{Cbraid2} this
forces $a=b$, putting us in Case 2 and contradicting the distinctness hypothesis.  
In the second case, the combined partition $\alpha_0 \alpha_\infty$ would 
have the form $(k_1,k_2,k_3,k_4,k_5)$ and the 
the combined partition $\beta_0 \beta_\infty$ would have 
the form $(3k_1,2k_2,k_2,2k_3,k_3,k_4^3,k_5^3)$ or
$(3k_1,3k_2,k_3^3,k_4^3,k_5^3)$. From \eqref{Cbraid2},
the first possibility occurs exactly when $a=d$ or $b=d$,
putting us into Case 3 and contradicting the 
distinctness hypothesis; the second possibility
cannot occur as it is incompatible with the shape of $\beta_0 \beta_\infty$ in \eqref{Cbraid2}.   We have now eliminated 
all possibilities for $\alpha_1$ and so $\AX \rightarrow \AP^1$ 
has to be primitive.     
 \qed

The two smallest degrees of covers as in Theorem~\ref{monodromy} 
are $m=21$ and $m=24$.  There are nine primitive groups in degree
$21$ and five in degree $24$, all accessible via {\em Magma}'s database
of primitive groups.   None of them, besides $S_{21}$ and $S_{24}$,
contain an element of cycle type $4^3 1^{m-12}$.

In an e-mail to the author on July 5, 2016, Magaard has sketched a 
proof that, for all $m \geq 25$, likewise $4^3 1^{m-12}$ is not a cycle partition for a 
primitive proper subgroup of $S_m$. His proof
appeals to Theorem~1 of \cite{GM98}, which has as essential hypothesis that the $1$'s in
$4^3 1^{m-12}$ contribute more than half the degree.  Special arguments are needed 
to eliminate the other possibilities that parts 1, 2, and 3
of Theorem~1 of \cite{GM98} leave open.   Thus, Theorem~\ref{monodromy} can be strengthened
by replacing {\em primitive} by {\em full}.  

\subsection{Primes of bad reduction}   The dessins $\gamma(a,b,d)$ have four white vertices: 
$0$, $1$, and the roots of $D$.  They have seven black vertices: $\infty$ and the roots of $ABC$.  If 
any of these eleven points agree modulo a prime $p$, then the map $\pi_{a,b,d}$ has
bad reduction at $p$.   To study bad reduction, one therefore has to consider
some special values, discriminants, and resultants.  Table~\ref{ABCDtab} gives the
relevant information, with ``value at $\infty$'' meaning the coefficient of $x^2$ in the 
quadratic polynomial heading the column.
\begin{table}[htb]
{\renewcommand{\arraycolsep}{3pt}
\[
 \begin{array}{c|ccccc}
    &                A & B &  C & D \\
   \hline
   \mbox{Value at $0$} &                (a+d) (a+2 d) & a (a+d)  & a (a+d) & a (a+2 d)\\
   \mbox{Value at $1$} &           b (b+d) & (b+d) (b+2 d)  & b (b+d) & b (b+2 d) \\
   \mbox{Value at $\infty$} & (a + b + d) n & (a + b + d) n  & (a + b) (a + b + d) &(a + b) n \\
   \mbox{Disc.} &                 -4 b d (a+d) n & -4 a d (b+d) n  & -4 a b d (a+b+d) & -8 a b d n \\
   \hline
 \mbox{Res.\ with $A$}&     & 4 d^3 n^2 e & 4 b^2 d^3 e & b^2 d^2 (a+2 d)^2 n^2  \\
 \mbox{Res.\ with $B$}&                   &   & 4 a^2 d^3 e & a^2 d^2 (b+2 d)^2 n^2  \\
 \mbox{Res.\ with $C$}&                   &  &    & a^2 b^2 (a+b)^2 d^2  \\
   \end{array}
 \]
 }
 \caption{\label{ABCDtab} Special values, discriminants, and resultants 
 of the four quadratic polynomials $A$, $B$, $C$, $D$ from Theorem~\eqref{clantheorem}, using the 
 abbreviation $e = (a+d)(b+d)(a+b+d)$.}
 \end{table}
 
 Combining the explicit formula \eqref{claneq}, the general discriminant formula \cite[(7.14)]{MR05},
 and the elementary facts collected in Table~\ref{ABCDtab}, one gets the following
 discriminant formula.  
 \begin{Corollary}
 \label{disccor}
Let $a^a b^b A^{a+d} B^{b+d} D^d - v 2^d d^{2d} n^n x^{a+2d} (1-x)^{b+2d} C^{n-d}$ 
be the polynomial whose vanishing defines $\pi_{a,b,d}$.  Its discriminant is
 \begin{eqnarray*}
\lefteqn{D(a,b,d)=} \\
&& (-1)^{ (a-1) a/2+ (b-1) b/2+d} 2^{n (d+2 n)} \\
  &&  a^{2 n^2 -a^2+2 a n-n} b^{2 n^2 -b^2 - n + 2bn} d^{ (10n^2 - (1 + a + b)(a + b + 3n))/2}  \\
  &&  (a+b)^{(a+b+d-1) n} (a+d)^{a n+a+d n+d+n^2-n} (b+d)^{  b n+b+d n+d+n^2-n} \\
  &&   (a+2 d)^{(a+2 d)^2} (b+2 d)^{(b+2 d)^2}  (a+b+d)^{(a+b+d) (2 (a+b+d)+1)} \\
    && n^{n (3 n+2)}   v^{3 n-7} (v-1)^9. \qed
\end{eqnarray*}
\end{Corollary}
\noindent The discriminants corresponding to $\pi_{a,d}^V$ and $\pi_{a,d}^S$ are
given by similar but slightly simpler formulas.  

\subsection{Allowing negative parameters}  The formula \eqref{claneq} makes sense for arbitrary integer parameters
$(a,b,d)$ satisfying $abdn \neq 0$, although individual factors may switch from numerator to denominator
or vice versa.  As a special case of \eqref{homogeneity}, one has the formula $\pi_{-a,-b,-d}(x) = \pi_{a,b,d}(x)^{-1}$.  
Using this symmetry,  we can
and will restrict attention to the half-space $d \geq 0$.  
 \begin{figure}[htb]
\includegraphics[width=4.9in]{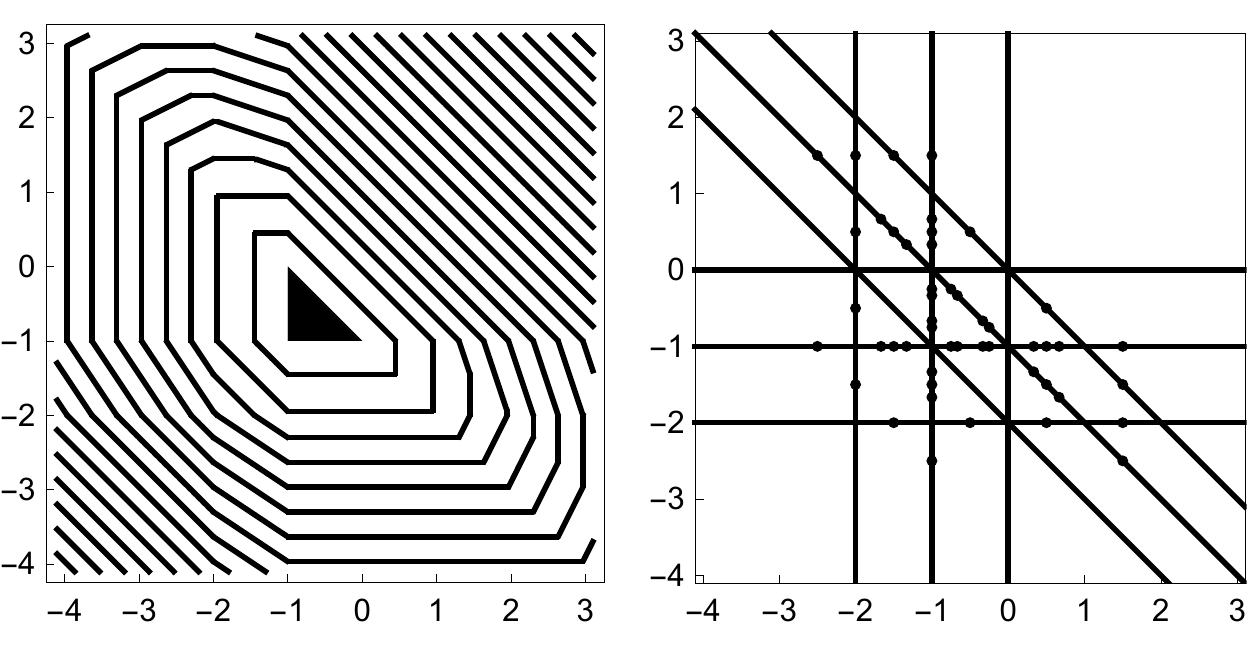}
\caption{\label{disclocus}  Left: the formal degree of $\pi_{\alpha,\beta,1}$, meaning the  quantity $\mbox{degree}(\pi_{\alpha d,\beta d,d})/d$,
drawn in the $\alpha$-$\beta$ plane. Contours range from $4$ (middle triangle) to $24$ (upper right corner). 
Right: The discriminant locus of the semicubical clan, drawn in the same $(\alpha,\beta)$ plane.  Points 
are particular $(a/d,b/d)$-values giving covers appearing in Table~\ref{modalgebrastab}.}
\end{figure}

Assuming none of the quantities in Table~\ref{ABCDtab} vanish, the degree $N(a,b,d)$
of $\pi_{a,b,d}$ is the total of those quantities on the list $2(d-n)$, $2(a+d)$, $2(b+d)$, $2d$, $-a-2d$, $-b-2d$, and
$a+b$ which are positive.  This continuous, piecewise-linear function is homogeneous in the parameters $a$, $b$, and $d$, and
so it can be understood by its restriction to $d=1$ via $N(a,b,d) = d N(a/d,b/d,1)$.  The left half of Figure~\ref{disclocus} is 
a contour plot of $N(\alpha,\beta,1)$.  Thus, for $d \geq 4$ fixed, the minimum degree for $\pi_{a,b,d}$ is $4d$, occurring for all 
$(a,b,d)$ with $(a/d,b/d)$ in the middle black triangle.

If $(a+d)(a+2d)(b+d)(b+2d)(n+d)(n+2d)=0$, then there is a cancellation among at least a pair of factors, and
$\pi_{a,b,d}$ has degree strictly less than $N(a,b,d)$.   If $a b d n=0$ then, taking a limit, $\pi_{a,b,d}$ 
is still naturally defined, and again has degree strictly less than $N(a,b,d)$.   Taking $d=1$, the lines given
by the vanishing of the other nine factors are drawn in the right half of Figure~\ref{disclocus}.    
The complement of these lines has $37$ connected components, called chambers.  The middle chamber is the
interior of the triangle given by $N(a,b,1)=4$.   As indicated by the caption of Figure~\ref{disclocus},
one can also think of the right half of Figure~\ref{disclocus} projectively.   From this viewpoint,
the line at infinity is given by the vanishing of the remaining linear form, i.e.\ by $d=0$.    
Our main reference \cite{Cou97} already illustrates some of this wall-crossing behavior in the context of
distinct $a$, $b$, $c$, and $d$:
the dessins with all parameters positive are described as being a {\em chardon}, while the dessins
with parameters in certain other chambers are described as being a {\em pomme}.  

\subsection{The central chamber and symmetric coordinates}
Each chamber corresponds to a different family of Hurwitz parameters, with corresponding
rational functions being uniformly given by \eqref{semicubicalF}.  To study the middle chamber, switch to 
new parameters $(u,v,w) = (d+a,d+b,d-n)$.  In the new parameters, the quantity $d=u+v+w$ is still
convenient.   The middle chamber is given by the positivity of $u$, $v$, and $w$.   We indicate
the presence of the new parameters by capital letters, changing $h$ to $H$, $\pi$ to $\Pi$, and 
$\gamma$ to $\Gamma$.  As in the case with other clans as well, both parametrization systems
have their virtues.  

The normalized Hurwitz parameter 
\eqref{semicubical} gets replaced by 
\begin{equation}
\label{semicubicalsym}
H(u,v,w) = (S_{2d}, (2 \, 1^{2d-2}, d^2, \; 3_x \, 1^{2d-3}, \; (d-u)_0 \, (d-v)_1 (d-w)_\infty), (1, \, 1, \, 1, \, 1)).
\end{equation}
Simply writing factors with the new parameters and in different places corresponding to the new signs, \eqref{claneq}
becomes
\begin{equation}
\label{claneqsym}
\Pi_{u,v,w}(x) = \frac{(-1)^{d-w} A^u B^v C^w D^d}{2^d d^{2d} (d-u)^{d-u} (d-v)^{d-v} (d-w)^{d-w} x^{d+u}
 (1-x)^{d+v}}.
\end{equation}
The degree is $4d$ and the partition triple \eqref{Cbraid2} changes to
\begin{eqnarray}
\nonumber \beta_0 & = &  u^2 v^2 d^2 w^2, \\
\label{Cbraid2sym} \beta_1 & = & 4^3 1^{4d-12}, \\
\nonumber \beta_\infty & = &  (d+u)_0(d+v)_1(d+w)_\infty.
\end{eqnarray}
The symmetry \eqref{sym1} in the new parameters takes the similar form
$\Pi_{u,v,w}(1-x) = \Pi_{v,u,w}(x)$.  But now the 
symmetry 
\begin{equation}
\label{sym2}
\Pi_{u,v,w}(1/x) = \Pi_{w,v,u}(x)
\end{equation}
is equally visible.  

In terms of a sheared version of Figure~\ref{disclocus} 
in which the central triangle is equilateral, the symmetries
just described generate the $S_3$ consisting of rotations
and flips of this triangle.  One has quadratic reduction as in 
\eqref{formulaV}
whenever two of the parameters are equal.   
One has cubic reduction as in \eqref{formulaS} whenever one 
of the parameters is $2d$.     

 \begin{figure}[htb]
\includegraphics[width=4.9in]{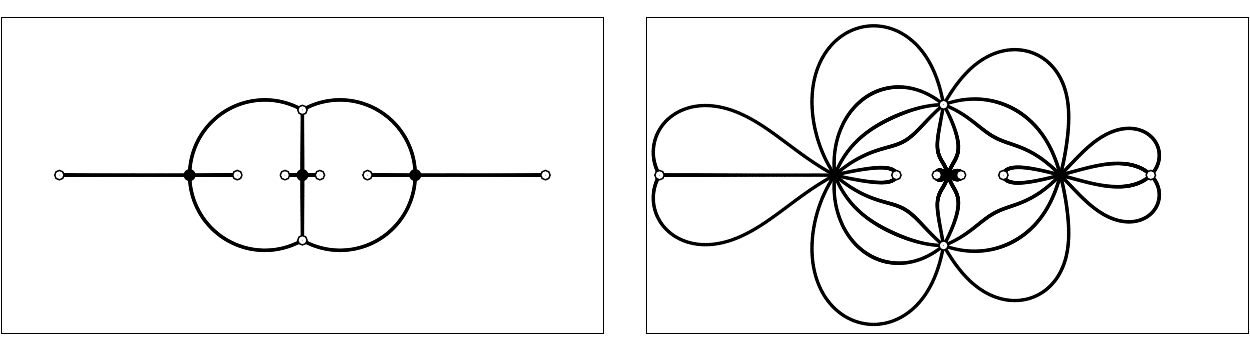}
\caption{\label{couvneg}  Left: $\Gamma(1,1,1)$.  Right: $\Gamma(4,3,2)$. }
\end{figure}

For dessins we take $\Gamma(u,v,w) = \Pi_{u,v,w}^{-1}([-\infty,0])$ with $[\infty,0] = \bullet \!\! - \!\!\! - \!\!\! - \!\! \circ$ 
as before.   Figure~\ref{couvpos} is then direct analog of Figure~\ref{couvpos}.  The black points on the 
real axis from left to right are, as before $0$, $\infty$, $1$.   The unique point above this axis is connected 
to $0$, $\infty$, and $1$ by respectively $u$, $w$, and $v$ edges.  To pass from $\Gamma(1,1,1)$ 
to $\Gamma(u,v,w)$, one replaces each edge by either $u$, $v$, or $w$ parallel edges, illustrated
by the example of $\Gamma(4,3,2)$.   

While symmetric parameters are motivated by the central chamber, they are often better for analysis of the entire
clan.  As an example, we consider an aspect about degenerations over discriminantal lines,
always excluding the intersections of these lines.   To begin, consider the numerator of the
 logarithmic derivative $\Pi'_{u,v,w}(x)/\Pi_{u,v,w}(x)$.  
In conformity with $\beta_1 = 4^3 1^{m-12}$, one gets that this numerator is the cube of a cubic polynomial 
$\Delta(u,v,w,x)$.  An easy computation gives  
 $\Delta(u,v,w,x) =$
\begin{equation}
\label{cubicDelta}
w x^3 (w-d) (w+d) + 3 w (u-d) (w-d) x^2+3 u  (u-d) (w-d) x-u (u-d) (u+d).
\end{equation}
Its discriminant has the completely symmetric form
\[
\disc_x(\Delta(u,v,w,x)) = 2^2 3^3 u v w (u+d)^2 (v+d)^2 (w+d)^2 d^3.
\]
By symmetry, to understand the degenerations at the nine lines visible on the right
half of Figure~\ref{disclocus}, one needs only to understand the degenerations
at the horizontal lines.  From bottom to top, the lines are given by 
$u=-d$, $u=0$, and $u=d$.  Visibly from \eqref{cubicDelta} the polynomials $\Delta(-d,v,w)$, 
$\Delta(0,v,w)$, and $\Delta(d,v,w)$ have $0$ as a root of multiplicity $1$, $2$, 
and $3$.  Continuing this analysis, the Hurwitz-Belyi maps with parameter 
on these lines has $\beta_1$ of the form $4^2 \, 1^{m-8}$, $4 \, 2 \, 1^{m-6}$, and $3 \, 1^{m-1}$.
The only discriminantal line not discussed yet  is the line $d=0$ at infinity.  Here
again by symmetry, we need to consider only the case $(u,v,-m)$ with 
$u$ and $v$ positive satisfying $u+v=m$.   In this case, one has 
$(\beta_0,\beta_1,\beta_\infty) = (u \, v, 2 \, 1^{m-2}, m)$.

\subsection{Moduli algebras}   In Section~\ref{Belyi}, we saw some unexpected splitting of moduli algebras for certain partition
triples $(\lambda_0,\lambda_1,\lambda_\infty)$.   The triples coming from the degenerations \eqref{ramtriple3} and \eqref{ramtriple2}
and the semicubical clan itself \eqref{Cbraid} give us many moduli algebras which have $\Q$ a factor.  
We computationally investigate small degree members of this collection here, as
a further illustration that Hurwitz-Belyi maps are very special among all Belyi maps.

\begin{table}[htb]
\[
{\renewcommand{\arraycolsep}{2.7pt}
\begin{array}{|r |rrr| lll | l | lllllll |}
\hline
m & a &     & d & \lambda_0 & \lambda_1 & \lambda_\infty & \;\;\; \mu & \multicolumn{7}{c|}{D}  \\
\hline
9 &   1 &    & 2  & 3 3  2  1 & 4 \, 2 \, 1^3 & 54 & 18+1 & - & \mathbf{2}^{37} & \mathbf{3}^{26} &  \mathbf{5}^{12} & 7^{2} &  &  \\ 
12 & 1 &    & 3 &  4 4  3  1 & 4 \, 2 \, 1^6 & 75 & 39+1 &  &  \mathbf{2}^{84} &  \mathbf{3}^{41} &  \mathbf{5}^{30} & \mathbf{7}^{26} & 11^{10} &   \\
15 & 1 &    & 4 & 5  5  4  1 & 4 \, 2 \, 1^9 & 96 & 60+1 &&  \mathbf{2}^{105} &  \mathbf{3}^{107} &  \mathbf{5}^{55} & 7^{17} & 11^{10} &\\
15 & 2 &    & 3 & 5 5  3  2 & 4 \, 2 \, 1^9 & 87 & 60+1 & -&  \mathbf{2}^{151} &  \mathbf{3}^{60} &  \mathbf{5}^{55} &\mathbf{7}^{44} && 13^{14} \\
\hline
\multicolumn{12}{c}{\;} \\
\hline
m & u & v & w & \lambda_0 & \lambda_1 & \lambda_\infty & \;\;\; \mu & \multicolumn{7}{c|}{D}  \\
\hline
  7 &  0 & 3 & -1  & 322 & 4 \, 2 \,1 & 511 & 3+1 & - & \mathbf{2}^3 & \mathbf{3} & \mathbf{5}^2 & 7 & &  \\
                   9 &  0 & 2 & 1  & 3321 & 4\,2 \,1^3 & 54 & 18+1 & - & \mathbf{2}^{37} & \mathbf{3}^{26} & \mathbf{5}^{12} & 7^2 && \\
                   10 &  0 & 4 & -5      & 5311 & 4 \,2 \,1^4 & 64 & 28+1 && \mathbf{2}^{56} & \mathbf{3}^{41} & \mathbf{5}^{20} & 7^2 &&  \\
                   10 &  0 & 4 & -1      & 433 & 4 \,2 \,1^4 & 721 & 28+1 && \mathbf{2}^{50} & \mathbf{3}^{34} & 5^{14} & \mathbf{7}^{19} &&  \\
                   10 &   0 & 3 & -5      & 5221 & 4 \,2 \, 1^4 & 73 & 31+1 && \mathbf{2}^{66} & \mathbf{3}^{30} & \mathbf{5}^{23} & \mathbf{7}^8 &&  \\
                   10 &  2 & -1 & -3        & 33211 & 4\,4\,1^2 & 532 & 33+1+1 && \mathbf{2}^{69} & \mathbf{3}^{42} & \mathbf{5}^{23} & 7^{21} && \\
                   11 &    0 & 5 & -2      & 533 & 4 \,2 \,1^5 & 821 & 38+1 && \mathbf{2}^{113} & \mathbf{3}^{45} & \mathbf{5}^{28} & 7^5 & 11^{10} &  \\
                   12 &   0 & 3 & 1      & 4431 & 4 \,2 \,1^6 & 75 & 39+1 && \mathbf{2}^{84} & \mathbf{3}^{41} & \mathbf{5}^{30} & \mathbf{7}^{26} & 11^{10} & \\
                   12 &  2 & 1 & -5  & 552 & 4 \,4 \, 1^4 & 72111 & 41+1 && \mathbf{2}^{74} & \mathbf{3}^{18} & \mathbf{5}^{31} & \mathbf{7}^{31} & 11^{15} & \\
                   12 &   0 & 5 & -6 & 6411 & 4\, 2 \,1^6 & 75 & 42+1 && \mathbf{2}^{86} & \mathbf{3}^{44} & \mathbf{5}^{27} & \mathbf{7}^{30} & 11^{20} &\\
\hline
\end{array}
}
\]
\caption{ \label{modalgebrastab}  Degrees $\mu$ and discriminants $D$ of the moduli algebras 
coming from the partition triple \eqref{ramtriple2} of $\pi_{a,d}^V$ and the partition triple \eqref{Cbraid} of $\Pi_{u,v,w}$.
The primes which are bad for the Hurwitz-Belyi map are in boldface. }
\end{table}

The dessin of $\pi_{1,1}^S$ has the shape $\circ \!\! - \!\! \bullet \!\! - \!\! \circ \! \langle \!\!\!\! \begin{array}{c} \bullet \\ \bullet \end{array}$.  The dessin of $\pi_{a,c}^S$ is obtained by replacing the middle edge by $a$ parallel edges and the remaining three edges by 
$c$ parallel edges.   There are no other dessins that share the partition triple \eqref{ramtriple3}.  So all moduli algebras of \eqref{ramtriple3}
are simply $\Q$.  

The dessin of $\pi_{1,1}^V$ have a more complicated shape and there are many dessins that share the partition triple
\eqref{ramtriple2}.    The four cases with $0 < a <d$ with $a+d \leq 5$ are in the top part of Table~\ref{modalgebrastab}.

For our main case of $\pi_{a,b,d}$, the triple satisfying the conditions of Theorem~\ref{monodromy} giving the lowest degree
$m = 3(a+b+2d)$ is $(a,b,d) = (3,2,1)$ with $m=21$.   This case is beyond our computational reach.  
Allowing $a$ and/or $b$ to be negative, but staying off the discriminantal hyperplanes, the 
lowest degree is $m=16$ from $(-1,-2,4)$.   This case may be easier, but instead
we allow $(a,b,d)$ to be on a discriminantal hyperplane, excluding the 
extreme degenerations $abdn=0$, as they give mass $\mu \leq 2$.   We still require that
the $a$, $b$, and $d$ are distinct without a common factor.    We switch to symmetric 
coordinates, so as to see the $S_3$ symmetries of the previous subsection clearly.  Modulo these symmetries,
 Table~\ref{modalgebrastab} gives all cases with $m \leq 12$. 
The top line is the septic example from Section~\ref{Belyi} yet again. 

The behavior summarized in Table~\ref{modalgebrastab} is similar to the
behavior discussed in \S\ref{Belyi}, but now the $(\lambda_0,\lambda_1,\lambda_\infty)$
have been chosen to ensure splitting.    Calculation shows that the
the large factor of the 
moduli algebra always has Galois group the full symmetric group
on the degree.  In the only case where there is extra splitting, the
two rational Belyi maps are
\begin{eqnarray}
\label{hur} \;\;\;\;\;\;\;\;\; \pi_{1,-1,2}(x)& \!\! \!\! = \!\! \!\! &  
\frac{(8 x-5)^2 \left(8 x^2-24 x+15\right)^3 \left(8 x^2-8
    x+3\right)}{2^{14} (x-1)^3 x^5 (4 x-3)^2}, \\
\label{nonhur} \;\;\;\;\;\;\;  \pi(x) & \!\! \!\! = \!\! \!\! & \frac{-(8 x-5)^2 \left(464 x^2-840
    x+375\right)^3 
  \left(2528 x^2-4400 x+1875\right)}{2^{16} 5^5 (x-1)^3 x^5 (4   x-3)^2}.
 \end{eqnarray}
 Like $\pi_{1,-1,2}(x)$, the unexplained rational factor $\pi(x)$ has bad prime
 set just $\{2,3,5\}$ and monodromy group $S_{10}$.

\subsection{Responsiveness to the inverse problem} 
Let $a$, $b$, $d$ be distinct positive integers without a common factor.
The fullness conclusion of \S\ref{primandfull} and the 
discriminant formula in Corollary~\ref{disccor} combine
to say that the explicit rational Hurwitz-Belyi maps $\pi_{a,b,d}$ 
of Theorem~\ref{clantheorem}
respond interestingly to the inverse problem of \S\ref{inverseproblem}.

More precisely, 
the set $\cP_{a,b,d}$ of bad primes of $\pi_{a,b,d}$ is the set
of primes dividing 
\[
a b d (a+b) (a+d) (b+d) (a+b+d) (a+2d) (b+2d)  (a+b+2d).
\]
This set can only contain primes at most $n=a+b+2d$.  From this 
fact alone, it is substantially smaller than the set 
than the set $\cP^{\rm glob}$ of \S\ref{bounds}, which is the the set of primes 
less than the degree $3n$. 

On the other hand, if one fixes a small set $\cP$, the largest
degree cover coming from the semicubical clan ramified within
$\cP$ is often very small, even if one works with all the
full $\Pi_{u,v,w}$, $\pi^V_{a,d}$, and $\pi^S_{a,d}$.   For example, taking 
$\cP = \{2,3,5\}$, one has relatively large degree covers coming from very degenerate cases
like $\Pi_{1,80,-1}$, with degree $81$.   However if one requires that $\beta_1$ contain a $4$, then
 the largest degree cover seems to be twenty-six.  This degree comes
 from two nondegenerate parameters $\Pi_{-5,-3,9}$ and $\Pi_{-5,-1,9}$, 
 with partition triples $(9^2 \, 4 \, 2 \, 1^2, 4^3 \, 1^{14}, 10 \, 5^2 \, 3^2)$
 and $(9^2 \, 3^2 \, 2, \, 4^3 \, 1^{14}, \; 12 \, 5^2 \, 2 \, 1^2)$ respectively.

We have looked at many clans, some quite different in nature from
the semicubical clan.  All seem to share the property
that the analog to $\cP_{a,b,d}$ is relatively sparse, but nevertheless grows 
with the parameters.  We are more interested
in this paper in fixing a small $\cP$  and
providing examples of full rational covers in degrees
as large as possible.  In this direction, the focus of 
Sections~\ref{23p}-\ref{conj}, 
clans do not seem to be helpful.   The fundamental
problem is that in clans the groups $G$ are $A_n$ or 
$S_n$ and one is increasing $n$.  
What one needs to do instead is fix $G$ and
let the number $r$ of ramifying points increase.

\section{The braid-triple method} 
\label{4vs3}
     Our first Hurwitz-Belyi map \eqref{deg7x} and the entire semicubical
 clan \eqref{claneq} were computed with the standard method.  After  \S\ref{braidbackground} gives background on braids,
 \S\ref{stepone} and \S\ref{steptwo} describe the alternative braid-triple method, already used twice in \S\ref{ambiguity}.   
 The two methods are complementary, as we explain in the short \S\ref{comparison}.
 
     A key step in the braid-triple method is to pass from a Belyi pencil
 $u$
  to a corresponding
 {\em braid triple} $B = (B_0,B_1,B_\infty)$.   In this paper,
 we use the braid-triple method only for four $u$, and
 the corresponding $B$ are given by the simple 
 formulas \eqref{braid1111}, \eqref{braid211}, \eqref{braid31}, and \eqref{braid41}.  
   We do not pursue the 
 general case here; 
  our policy in this paper is to be 
  very brief with 
  respect to braid groups, saying
  just enough to allow the reader to
  replicate our computations of individual covers.   
   
 \subsection{Algebraic background on braid groups} 
 \label{braidbackground}
 The {\em Artin braid group} on $r$ strands is the most widely known of all braid groups, 
 and our summary here follows \cite[\S3]{RV15}.  The group is defined via $r-1$ generators and $\binom{r-1}{2}$ relations:
\begin{equation}
{\renewcommand{\arraycolsep}{3pt}
\label{braiddef}
\Br_r = \left< \sigma_1,\dots,\sigma_{r-1} : 
 \begin{array}{rclrcl} \sigma_i \sigma_j & = & \sigma_j \sigma_i,  & \mbox{ if } |i-j| & > & 1 \\
 \sigma_i \sigma_j \sigma_i & = & \sigma_j \sigma_i \sigma_j, & \mbox{ if } |i-j| & = & 1 \\
  \end{array} \right>.
 }
\end{equation}
The assignment $\sigma_i \mapsto (i, i+1)$ extends to a surjection $\Br_r \twoheadrightarrow S_r$.
For every subgroup of $S_r$ one gets a subgroup of $\Br_r$ by pullback.  
Thus, in particular, one has surjections $\Br_{\nu} \twoheadrightarrow S_{{\nu}}$ for
$S_\nu = S_{\nu_1} \times \cdots \times S_{\nu_r}$.    

Given a finite group $G$, let $\cG_r \subset G^r$ be the set of tuples $(g_1,\dots,g_r)$
with the $g_i$ generating $G$ and satisfying $g_1 \cdots g_r = 1$.  The braid group $\Br_r$ acts on the
right of $\cG_r$ by the braiding rule
 \begin{equation}
\label{bdef}
(\dots, \; g_{i-1}, \;g_i, \; g_{i+1}, \; g_{i+2}, \; \dots)^{\sigma_i} = (\dots,g_{i-1}, \; g_{i+1}, \;  g_i^{g_{i+1}}, \; g_{i+2},\dots).
\end{equation}
The group $G$ acts diagonally on $\cG_r$ by simultaneous conjugation.   The actions of $\Br_\nu$ 
and $G$ commute with one another.  

Let $h = (G,C,\nu)$ be an $r$-point Hurwitz parameter, 
with $C = (C_1,\dots,C_k)$, $\nu = (\nu_1,\dots, \nu_k)$, and $\sum_{i=1}^k \nu_i=r$,
all as usual.  Consider the subset $\cG_h$ of $\cG_r$, 
consisting of tuples with $g_j \in C_i$ for 
$\sum_{a=1}^{i-1} \nu_a <j \leq \sum_{a=1}^i \nu_a$.
This subset is stable under the action of $\Br_\nu \times G$.    
The group $\Br_\nu$ therefore acts on the right of the 
quotient set $\cF_h = \cG_h/G$.   

The actions of $\Br_\nu$ on $\cF_h$ all factor through a certain
quotient ${\br}_\nu$.   Terminology is not important for us here, but for comparison
with the literature we remark that $\br_\nu$ is the quotient
of the standard Hurwitz braid group by its two-element center.  
Choosing a base point $\star$ and certain identifications appropriately, the group
${\br}_\nu$ is identified with the fundamental group 
$\pi_1(\AConf_\nu,\star)$, in a way which makes the
 action of ${\br}_\nu$ on $\cF_h$ agree 
with the action of  $\pi_1(\AConf_\nu,\star)$ on 
the base fiber $\pi_h^{-1}(\star) \subset \AHur_h$. 

Let $Z$ be the center of $G$, this center being trivial for most of our examples.   
Let $\Aut(G,C)$ be the subgroup of $\Aut(G)$ which stabilizes each conjugacy
class $C_i$ in $C$.  Then not only does $G/Z$ act diagonally on $\cG_h$, but
so does the entire overgroup $\Aut(G,C)$.    The group $\Out(G,C)$ introduced in \S\ref{Hurwitzmaps} is the quotient
of $\Aut(G,C)$ by $G/Z$.  Quotienting by the natural action of a subgroup $Q \subseteq \Out(G,C)$, 
gives the base fiber $\cF_h^Q$ corresponding to the cover $\AHur^Q_h \rightarrow \AConf_\nu$.
 
 \subsection{Step one:  computation of braid triples} 
 \label{stepone}
       A Belyi pencil $u : \AP^1_v \rightarrow \AConf_\nu$ determines, up to conjugacy depending
       on choices of base points and a path between them,  an {\em abstract braid triple} $(B_0,B_1,B_\infty)$ of elements of ${\br}_\nu$.   These
  elements have the property that in any Hurwitz-Belyi map $\pi_{h,u} : \AX_{h,u} \rightarrow \AP^1$,
  the images of the $B_\tau$ in their action on $\cF_{h}$ give the global monodromy of the cover.  
  When $\cF_{h}$ is identified with $\{1,\dots,m\}$, we denote the image of $B_\tau$ by $b_\tau \in S_m$ 
  and its cycle partition by $\beta_\tau$.   We call $(b_0,b_1,b_\infty)$ a {\em braid permutation
  triple}.  As we have already done several times before, e.g.\  \eqref{Cbraid2},  we 
  call $(\beta_0,\beta_1,\beta_\infty)$ a {\em braid partition triple}.

For the three $4$-point Belyi pencils introduced in \eqref{BP2}, the abstract braid triples are 
 \begin{align}
\label{braid1111} \;\;\;\;\;\;\;\;\;\;\;\;  u_{1,1,1,1}\!\! : && (B_0,B_1,B_\infty) & =  (\sigma_1^2, \sigma_2^2, \sigma_2^{-2} \sigma_1^{-2}), \;\;\;\;\;\;\;\;\;\;\;\; \\
\label{braid211}  u_{2,1,1} \!\!: && (B_0,B_1,B_\infty) & =   (\sigma_1, \sigma_1^{-1} \sigma_2^{-2},\sigma_2^2), \\
\label{braid31} u_{3,1}\!\! : && (B_0,B_1,B_\infty) & =  (\sigma_1 \sigma_2, \sigma_2^{-1}  \sigma_1^{-2}, \sigma_1).
 \end{align}
The triple for $u_{1,1,1,1}$ is given in \cite[\S5.5.2]{LZ04}.  The other two can be deduced by quadratic and then
cubic base change, using \eqref{basechange}.   An important point is that, in the quotient group $\br_\nu$, 
the $B_1$ for $u_{2,1,1}$ and $u_{3,1}$
have order $2$ and the $B_0$ for $u_{3,1}$ has order $3$.  
 
 It is worth emphasizing the conceptual simplicity of our braid computations.
 They  repeatedly use the generators
$\sigma_i$ of \eqref{braiddef} and their actions on $\cF_h$ from
\eqref{bdef}.  However they do not explicitly use the relations in
 \eqref{braiddef}.  Likewise they do not explicitly
 use the extra relations involved in passing from $\Br_\nu$ to $\br_\nu$.
 Our actual computations are at the level of the permutations $b_\tau \in S_m$ 
 rather than the level of the braid words $B_\tau$.  At the permutation
 level, all these relations automatically hold.  
 
 Computationally, we realize $\cF_h$ via a set of 
 representatives in $\cG_h$ for the conjugation 
 action.  A difficulty is that the set  $\cG_h$ in
 which computations take place is large.  Relatively
 naive use of \eqref{braiddef} and
 \eqref{bdef} suffices for the braid computations
 presented in the next five sections.  To 
 work as easily with larger groups and/or larger degrees,
 a more sophisticated implementation as in 
 \cite{MSV03} would be essential.  

\subsection{Step 2:  passing from a braid triple to an equation}
\label{steptwo}  Having computed a braid partition triple $(\beta_0,\beta_1,\beta_\infty)$
belonging to $\pi_{h,u}$, one can then try to pass from the triple to 
an equation for $\pi_{h,u}$ by algebraic methods.   
We did this in \S\ref{ambiguity} for two partition triples with
degree $m=12$.  For one, as described in \S\ref{deg12},
 the desired $\pi_{h,u}$ 
is just one of $\mu=24$ locally equivalent covers, the one
defined over $\Q$.  

The braid partition triples $(\beta_0,\beta_1,\beta_\infty)$ 
arising in the next five sections have degrees $m$ into 
the low thousands.   The number $\mu$ of Belyi maps with 
a given such braid partition triple is likely to be more than
  $10^{100}$  in some cases.    The numbers 
  $12$ and $24$ are therefore being replaced by very much larger
  $m$ and $\mu$ respectively.  It is completely impractical to follow the purely
 algebraic approach of getting all the maps belonging to the 
 given  $(\beta_0,\beta_1,\beta_\infty)$ and extracting the
 desired rational one.  Instead, there are three more feasible
 techniques for computing only the desired cover.
 
 For almost all the covers in this paper, Step~2 was carried out by a $p$-adic 
 technique for finding covers defined over $\Q$ explained in detail in \cite{Mal88}.  Here one picks a good
 prime $p$ for the cover sought, and first searches
 for a tame cover with the correct $(\beta_0,\beta_1,\beta_\infty)$ 
 defined over $\F_p$.   Commonly, one finds several
 covers, and one cannot yet tell which is the reduction of 
 the cover sought.    One then uniquely lifts all these candidates 
 iteratively to $\Z/p^c$ for some large 
 $c$.  This step requires solving linear equations and is
 easy.  We commonly took $c=50$.  Then one recognizes
 the coefficients of the lifted covers as $p$-adically near 
 rational numbers.  In practice this is easy too, and only one 
 of the initial solutions over $\F_p$ gives small height rational numbers.  One concludes
 by checking that the monodromy of the cover constructed really
 does agree with the braid permutation triple $(b_0,b_1,b_\infty)$.  
 The efficiency of this $p$-adic technique decreases rapidly with
 $p$.   Since the covers we pursue all have a small
 prime of good reduction, typically $5$ or $7$, the technique
 is well adapted to our situation.  
 
 The second and third technique have been recently introduced, and both
 are undergoing further development. They take the permutation triple $(b_0,b_1,b_\infty)$
 rather than the partition  triple $(\beta_0,\beta_1,\beta_\infty)$ as a starting point.  Thus they
 isolate the cover sought immediately, and there is no issue of a large 
 local equivalence class.  The technique of \cite{KMSV14} centers
 on power series while the technique of \cite{Kra15} centers on numerically
 solving partial differential equations.   Schiavone used the programs
 described in \cite{KMSV14} to compute
 \eqref{prime17} here.  Our tables in Sections~\ref{235}-\ref{fivepoint}
 present braid information going well beyond where one can currently
 compute equations, in part to provide 
 targets for these developing computational
 methods. 
 
\subsection{Comparison of the two methods}
\label{comparison}
We used the standard method many time in \cite{RobHNF} 
in the context of constructing covers of surfaces.  In the present
context of curves, the braid-triple method complements the standard method 
as follows.  In the standard method, one can expect the 
difficulty of the computation to increase rapidly with the genus $g_Y$ of $\AY_x$ 
and the degree $n$ of the cover $\AY_x \rightarrow \AP^1_t$.  In
the braid-triple method, these measures of difficulty are replaced
by the genus $g_X$ of the curve $\AX_{h,u}$ and the degree
$m$ of the cover $\AX_{h,u} \rightarrow \AP^1_v$. 
The quantities $(g_Y,n)$ and $(g_X,m)$ are not tightly 
correlated with each other, and in practice each method has
 a large range of parameters 
for which it works well while the other method does not.

\section{Hurwitz-Belyi maps exhibiting spin separation}
\label{separation}
This section presents three Hurwitz-Belyi maps for which we were able to find
a defining equation by both the standard and the braid-triple method.   
 Each map has the added interesting feature that  the covering curve $\AX_h$ has two 
components.  We explain this splitting by means of 
lifting invariants.  Many of the covers in the next
four sections are similarly forced to split via lifting invariants.  

\subsection{Lifting in general} 
\label{liftingingeneral} Decomposition of Hurwitz varieties was studied by Fried and Serre.
Here we give a very brief summary of the longer summary given in 
\cite[\S4]{RobHNF}.    The decompositions come from central extensions $\tilde{G}$ of the
given group $G$.   The term {\em spin separation} is used because many
double covers are induced from the double cover $\mbox{Spin}_n$ of the special
orthogonal group $SO_n$ via an orthogonal representation.  

  Let $h=(G,C,\nu)$ be a Hurwitz parameter.  First, one has the Schur multiplier $H_2(G,\Z)$, always abbreviated in
  this paper as $H_2(G)$.
 Any universal central extension $\tilde{G}$ of  $G$ has the form $H_2(G).G$.   Second, one has a quotient
 $H_2(G,C)$ of the Schur multiplier, with $H_2(G,C).G$ being the largest quotient in which each $C_i$ splits
 completely into $|H_2(G,C)|$ different conjugacy classes.  Third, one has a torsor $H_h = H_2(G,C,\nu)$ over $H_2(G,C)$.   So $|H_h| = |H_2(G,C)|$, 
 but the set $H_h$ does not necessarily have a distinguished point like the group $H_2(G,C)$ does.   
 
  The group $\Out(G,C)$ defined in \S\ref{Hurwitzmaps} acts on the set $H_h$.  For any
  subgroup $Q \subseteq \Out(G,C)$, one has a natural map from the component set $\pi_0(\AHur^Q_h)$ of
  $\AHur_h^Q$ to $H^Q_h$.
   The most common behavior is that these maps $\pi_0(\AHur^Q_h) \rightarrow H_h^Q$ 
   are bijective.  As said already
   in \S\ref{Hurwitzmaps}, our main interest is in $Q = \Out(G,C)$, in which
   case we replace $\Out(G,C)$ by $*$ as a superscript.  
    
   In practice, the key groups $H_2(G,C)$ and $\Out(G,C)$ are extremely small.  
   In the next three subsections $H_2(G,C)$ has order $2$, $3$, and $3$ respectively,
   while $\Out(G,C)$ has order $1$, $2$, and $1$.   We explain lifting 
   in some detail in these subsections and also in \S\ref{235second} where $H_2(G,C)$ and
   $\Out(G,C)$ can be slightly larger.  

\subsection{A degree $25=15+10$ family} 
\label{quadratic} 
 Applying the mass formula
\cite[(3.6)]{RobHNF} to a Hurwitz parameter $h = (G,C,\nu)$  requires the use of 
the character table of $G$.    Common choices for 
$G$ in this paper are $A_5$ and $S_5$.  Table~\ref{chartable}   
gives the character table for these two groups,
as well as their Schur double covers $\tilde{A}_5$,
and $\tilde{S}_5$.  In this subsection, we use this table to illustrate
how mass formula computations for a given
Hurwitz parameter $h$ appear in practice, 
including refinements involving covering groups $\tilde{G}$.  

\begin{table}[htb]
\[
\begin{array}{r |  rrrrr| rrr}
|C_i|          &    1 & 15 & 20 & 12 & 12 & 10 & 30 & 20 \\
 C_i           &   1^5 & 221 & 311 & 5a & 5b & 2111 & 41 & 32 \\
            \hline 
\chi_1 &  1 & 1 & 1 & 1 & 1 & 1 & 1 & 1 \\
\chi_2 & 3 & -1 & 0 & -b & -b' &  0 & 0 & 0 \\
\chi_3 &  3 & -1 & 0 & -b' & -b' &   &  &  \\
\chi_4 & 4 & 0 & 1 & -1 & -1 & 2 & 0 &  -1 \\
\chi_5 & 5 & 1 & -1 & 0 & 0 & 1 & -1 & 1 \\
\hline
\mbox{Orders of }& 1 & 4 & 3 & 5 & 5 & 2 & 8 & 6 \\
\mbox{classes }& 2 &    & 6 & 10 & 10 &   & 8 & 6 \\   
\hline 
\chi_6 & 2 & 0 & -1 & b & b' & 0 & 0 & 0 \\
\chi_7 & 2 & 0 & -1 & b' & b &  &  & \\
\chi_8 & 4 & 0 & 1 & -1 & -1 & 0 & 0 & \sqrt{-3} \\
\chi_9 & 6 & 0 & 0 & 1 & 1 & 0 & \sqrt{-2} & 0 \\
\end{array}
\]
\caption{\label{chartable} Character tables for $A_5$, $S_5$, $\tilde{A}_5$, and
$\tilde{S}_5$, with the abbreviations $b=(-1+\sqrt{5})/2$ and $b' = (-1-\sqrt{5})/2$. }
\end{table}

The character table for $A_5$ is given simply by the upper 
left $5$-by-$5$ block.
The remaining character
tables require the use of Atlas conventions.  
The double cover $\tilde{A}_5$ has the listed nine
characters.  The nine conjugacy classes arise
because all but the class $221$ splits into 
two classes.  We label these classes of $\tilde{A}_5$ according
to whether the order of a representing element
is even ($+$) or odd $(-)$.   Thus $5a$ splits into $5a+$ and
$5a-$.  The printed character values
refer to the class with odd order elements.  
Thus, e.g., $\chi_8(311+) = 1$ but $\chi_8(311-) = -1$.

Only the groups $A_5$ and $\tilde{A}_5$ are used
in the example of this subsection, but $S_5$ and $\tilde{S}_5$
are equally common in the sequel and we explain them
here.  The group $S_5$ has  seven conjugacy classes, the
classes $5a$ and $5b$ having merged to a single class
$5$.  The corresponding seven characters are 
the printed $\chi_1$, $\chi_4$, $\chi_5$, 
the sum $\chi_2+\chi_3$ extended by
zero, and the twists $\chi_1 \epsilon$, $\chi_4 \epsilon$, 
and $\chi_5 \epsilon$.   Here $\epsilon$ is the sign character,
taking value $1$ on $A_5$ and $-1$ on $S_5-A_5$.  
The cover $\tilde{S}_5$ has twelve characters, 
the seven from before and the five new ones
$\chi_6+\chi_7$, $\chi_8$, $\chi_8 \epsilon$, $\chi_9$, 
and $\chi_{9} \epsilon$.

For the example of this subsection, let $h = (A_5,(311,5a),(3,1))$.   Because of $0$'s appearing as character values, only the characters $\chi_1$ and $\chi_4$ appear when
evaluating the mass formula:
\[
\overline{m}_h =   \frac{|C_1|^3 |C_2|}{|G|^2} \sum_{i=1}^5 \frac{\chi_i(C_1)^3 \chi_i(C_2)}{\chi_i(1)^2} =
 \frac{20^3 12}{60^2} \left( \frac{1^3(1)}{1^2} + \frac{1^3 (-1)}{4^2} \right) = \frac{20}{3} \frac{15}{16} = 25.
\]
Because $A_5$ does not have a proper subgroup meeting the both the classes
$311$ and $5a$, the desired degree is just the mass, $m_h = \overline{m}_h =  25$.  

The joint paper \cite{RV15} was originally planned to include this $h$ as an example.  
 The curves $\AY_x$ parameterized have genus one
but Venkatesh nonetheless computed the Belyi map
$\pi_{h} : \AX \rightarrow \AP_j^1$ by the standard method, seeing directly
that $\AX$ breaks into two components, each of genus zero, 
of degree $15$ and $10$ over
$\AP_j^1$.   The present author simultaneously used the braid-triple method, 
using \eqref{braid31} to get the braid partition triples $(3331,22222,541)$ and $(33333,22222221,5433)$.  Both methods 
end at the explicit equations \eqref{common541} and \eqref{common5433}.  

To explain the splitting, consider the Hurwitz parameters 
\begin{align}
\label{15plus10}
h^+ & = (\tilde{A}_5,(311+,5a+),(3,1)), & 
h^- & = (\tilde{A}_5,(311+,5a-),(3,1)).
\end{align}
Let $(g_1,g_2,g_3,g_4) \in \cG_h$.  For $i=1,2,3$, let 
$\tilde{g}_i$ be the unique preimage of
$g_i$ in $311+$.   Then there is a unique lift $\tilde{g}_4$ of
$g_4$ which satisfies 
$\tilde{g}_1 \tilde{g}_2 \tilde{g}_3 \tilde{g_4} = 1$.  
This lift can be in either $5a+$ or $5a-$.  
In this way one gets a map from $\cG_h$ to $H_h = \Z/2$. 
 This invariant does not change under either the braid 
 or conjugation action.  

The mass formula \cite[(3.6)]{RobHNF}, applied to $\tilde{G}$ now, lets one find the degrees of the factors.  
In this simple case, where proper subgroups of $\tilde{G}$ are not involved, one has
\begin{eqnarray*}
m_h^{\pm} & = & \frac{1}{2} \frac{|C_1|^3 |C_2|}{|G|^2} \sum_{i=1}^9 \frac{\chi_i(C_1)^3 \chi_i(C_2)}{\chi_i(1)^2}\\
& = & \frac{1}{2} \frac{20^3 12}{60^2} \left( \frac{1^3(1)}{1^2} + \frac{1^3 (-1)}{4^2} \pm 
\frac{(-1)^3 (b+b')}{2^2} \pm \frac{1^3 (-1)}{4^2} \right) \\
& = & \frac{40}{3} \left( 1 - \frac{1}{16} \pm \frac{1}{4} \pm \frac{-1}{16} \right) = \frac{40}{3} \left(\frac{15}{16} \pm \frac{3}{16} \right) = \frac{5}{2} \left(5 \pm 1\right) = 15,10.
\end{eqnarray*} 
Similar mass computations let one properly identify components
with lifting invariants in general.  We make these identifications,
typically with no further comments, in the next two subsections
and many times in \S\ref{235}-\ref{fivepoint}.

 \subsection{A degree $70=30+40$ family: 
 rational cubic splitting} 
 \label{cubic}  
 Let 
 \[
 h = (A_7, \; (22111, \, 511, \, 322), \; (2, \, 1, \, 1)).
 \]
 The large singletons $5$ and $3$ help keep the standard method 
 within computational feasibility.   By a direct application of this method,
 one sees at the end that the degree $m$ is $70$ and there
 is a splitting into two components of degrees $30$ and $40$.  
 
 In the braid-triple method the order of events is reversed.   
 Mass formula computations says that the desired 
$\AX^*_h \rightarrow \AP^1_w$ has degree $70$.  A braid group 
computation using \eqref{braid211} says that $\AX^*_h$ has two components
 of degrees $30$ and $40$.  The monodromy groups are
 $A_{30}$ and $S_{40}$ respectively, with braid partition triples 
 \begin{eqnarray*}
( \beta_0,\beta_1,\beta_\infty) & = & (7^2 \, 5 \; 3 \; 2^4, \;  2^{14} \, 1^2, \; 6 \;5^2 \, 4 \; 3^3 \,1), \\
( \beta_0,\beta_1,\beta_\infty) & = & (7^2 \, 5^2 \, 4^2 \, 2^3 \, 1^2, \;  2^{20}, \;  5^2 \, 4^3 \, 3^6).
\end{eqnarray*}
As the total number of parts is 32 and 42 respectively, 
the genus is zero in each case.   

The second step in the braid-triple method is challenging,
 since the smallest prime not in $\cP_{A_7}$ is
$11$.  This step is only within feasibility because of the splitting $70=30+40$,
and the fact that one can compute the two components independently.
Explicit equations are
 \begin{eqnarray*}
f_{30}(w,x) &=& 2^2 3^3  \left(7 x^2+14 x+4\right)^7 x^5 (2 x+1)^3 \left(x^2+3 x+1\right)^2 \left(2 x^2+x+2\right)^2 \\ 
    && + w \left(7 x^2+6 x+2\right)^5 (5 x+2)^4 
    \left(14 x^3+39 x^2+18 x+2\right)^3 (x+2),   \\
    &&\\
f_{40}(w,x) & = &  2^2 3^4 \left(5 x^2-12 x+3\right)^7  \left(5 x^2-15 x+12\right)^5  \left(x^2-3 x+6\right)^4   \\ &&
\qquad \!\!\!\! 
 \left(4 x^2-15 x+15\right)^2 x (5 x-9)  \\ &&+ w\left(x^2-3\right)^5 \left(5 x^3-45 x^2+120 x-108\right)^4 \\
&&  \qquad  \!\!\!\!  \left(400 x^6-2700 x^5+7425 x^4-10530 x^3+7830 x^2-2430  x-27\right)^3 \!\! .
 \end{eqnarray*}  
The polynomial discriminants are 
\begin{eqnarray*}
D_{30}(w) & = & -2^{450} 3^{285} 5^{95} 7^{105} w^{22} (w-1)^{14}, \\ 
D_{40}(w) & = & 2^{930} 3^{1254} 5^{230} 7^{105} w^{29} (w-1)^{20}. 
\end{eqnarray*}
Modulo squares these quantities are $-105$ and $7 w$, reflecting the fact 
that monodromy groups and generic Galois groups 
are $(A_{30},S_{30})$ in the first case
and $(S_{40},S_{40})$ in the second.     

 The group $A_7$ has
the unusually large maximal non-split central extension $6.A_7$.  
For both this subsection and the next, only the subextension
$3.A_7$ is relevant because all the classes $C_i$ split in it, while
the class $22111$ is inert in $2.A_7$.  In the notation reviewed in \S\ref{liftingingeneral},
this reduction is expressed by an identification  $H_2(A_7,C) = \Z/3$.
The group $\mbox{Out}(A_7,C)$ is
all of $\mbox{Out}(A_7)=\{\pm 1\}$ because the
classes $22111$, $511$,  and $322$ are
all stabilized by the outer involution.   The action of $\{\pm 1\}$ on
$\Z/3$ is the nontrivial one where $-1$ acts by negation. 
The degree $30$ component
corresponds to the identity class $0 \in \Z/3$ while 
the degree $40$ component corresponds to the orbit
of the two nonidentity classes in $\Z/3$. 

 The promised third conceptual explanation of the degree splitting
 $4= 3+1$ for $h = (A_7,(322, 421, 511), \, (1,1,1))$ from
 \eqref{split3plus1} is in our present context.  All three classes split in
 $3.A_7$ while only the last two split in $2.A_7$. 
 All three classes are stable under outer involution.
 So here again $\Out(A_7,C) = \{\pm 1\}$ acts
 nontrivially on $H_2(A_7,\Z) = \Z/3$.  In this case
 the degree one factor corresponds  to $0 \in \Z/3$ while  
 the degree three factor corresponds to $\{-1,1\} \subset \Z/3$.

\subsection{A degree $42=21+21$ family: 
irrational cubic splitting}  
\label{cubic2}
Lando and Zvorkin \cite[\S5.4]{LZ04} 
investigated splitting of Hurwitz covers
in some generality.  
The unique splitting in their context that
they could not conceptually explain 
comes from the Hurwitz parameter
\[
h = (A_7, \;  (22111, \, 7a), \; (3, \, 1)).
\]
Here one has splitting of the form 
$42=21+21$.   In this subsection 
we complement their study of this example
by both giving an equation and 
explaining the splitting. 

Computing using \eqref{braid31}, one gets that the two components 
have the same braid partition triple, namely
\begin{equation}
\label{braidlz}
(\beta_0, \, \beta_1, \, \beta_\infty) =    (3^7, \; 2^{10} \, 1, \; 6 \, 5 \, 4 \, 3^2).  
\end{equation}
This agreement is in contrast to the previous subsection,
where the two components even had different degrees.  
Lando and Zvonkin speculated (p.~333) that the agreement
is explained by the two components being 
Galois conjugate.  

Indexing the two maps arbitrarily by $\epsilon \in \{+,-\}$,  each map
we seek fits as the right vertical map in a Cartesian square:
\begin{equation}
\begin{array}{ccc}
\tilde{\AX}^\epsilon & \rightarrow & \AX^\epsilon \\
\downarrow & & \downarrow \\
\AP^1_v & \rightarrow & \AP^1_j 
\end{array}.
\end{equation}
Here the bottom map is the degree six $S_3$ cover given in \eqref{basechange},
and so the top map is a degree six $S_3$ cover as well.  

There are $7+11+5=23$ parts in all in \eqref{braidlz}, so that
the genus of each $\AX_\epsilon$ is $0$ by the Riemann-Hurwitz formula.  Lando and Zvonkin 
worked first with the base-changed cover.  Here the braid partition
is $(\tilde{\beta}_0, \, \tilde{\beta}_1, \, \tilde{\beta}_\infty)$, with
each $\tilde{\beta}_\tau = 5333322$.  As $7+7+7=21$,
the genus is $1$.   Jones and Zvonkin \cite{JZ02}
carried out the $S_3$ descent as we are doing here.   

   We find
via explicit equations that the two components are indeed conjugate with respect 
to the two choices
$s = \pm \sqrt{21}$:
\begin{eqnarray*}
\lefteqn{f_{21\pm}(j,x) =} \\
 && 2 \cdot 5 \cdot  \left(2700000 x^7+x^6 (630000 s-3780000)+x^5 (724500-829500
    s) \right.  \\ &&  \qquad \left. +x^4 (2228100-474600 s)  +x^3 (1404725 s-5328225) \right.  \\ && \qquad  \left. +x^2 (7020216-1485456
    s)+x (856800 s-4384800)-252000 s+972000\right)^3 \\ &&  \pm 
    3 \cdot 7^7 j (5239 s-21429) x^5 (10 x-9)^4
   \left(150 x^2+x (40 s-15)-8 s+88\right)^3.
\end{eqnarray*}
 \begin{figure}[htb]
\includegraphics[width=4.6in]{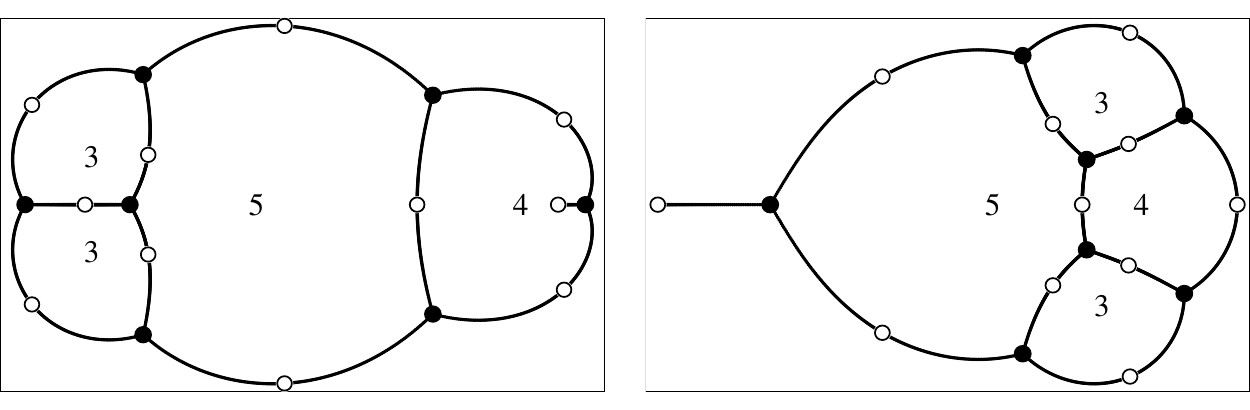}
\caption{\label{twin21}  Conjugate dessins, with $21+$ on left and $21-$ on right}
\end{figure}
Figure~\ref{twin21} draws 
the dessins in $\AP^1_x$ corresponding to Cover $21+$ on the left and its conjugate Cover $21-$ on the right.
The preimages in $\AP^1_x$ of $\infty \in \AP^1_j$ are indicated in the picture by their ramification numbers, 
with the undrawn $\infty \in \AP^1_x$ also being a preimage with ramification number $6$.  
Figure~\ref{twin21} gives the correct analytic shape of Figure~3 of \cite{JZ02}, 
and, after base change, the correct shape of Figure~5.15 of \cite{LZ04}.  

The splitting is induced by the existence of $3.A_7$ as in the previous 
subsection.  Again one has an identification $H_2(A_7,C) = \Z/3$.  
Here however, because $7a$ is not stabilized by the outer involution
of $A_7$, the group $\Out(A_7,C)$ is trivial.   Accordingly
one has a natural function from components of $\AX$ to $\Z/3$.  
One would generally expect all three preimages to have one component
each.  In this case, the preimages of $0$, $1$, and $-1$ are respectively
empty, $\AX^+$ and $\AX^-$. 

\section{Hurwitz-Belyi maps with $|G| = 2^a 3^b p$ and $\nu = (3,1)$} 
 \label{23p} 
    In this section, we set up a framework for studying the Hurwitz-Belyi maps 
coming from a systematic collection of $4$-point Hurwitz parameters $h$.  
Here and in the next two sections, we 
carry out the first part of the braid-triple method for all these $h$, obtaining braid permutation
triples $(b_0,b_1,b_\infty)$ and thus braid partition triples
$(\beta_0, \beta_1, \beta_\infty)$.   In many low degree cases, we 
carry out the second part as well, obtaining a defining
equation for the cover.    

\subsection{Restricting to $|\cP|=3$ and $\nu = (3,1)$} 
\label{23poverview} 
To respond to the inverse problem of \S\ref{inverseproblem}, we consider only 
$h = (G,C,\nu)$ giving covers defined over $\Q$.   
To keep our computational study of manageable size we 
impose two  severe restrictions.  First, we require that 
$G$ be almost simple with exactly three primes 
dividing its order.    Second, we restrict attention to 
the case $\nu = (3,1)$.  There are many more cases within 
computational reach which are excluded because of
these two restrictions.  The rest of this subsection
elaborates on the two restrictions.  

\subsubsection*{Almost simple groups divisible by at most three primes}  
There are exactly eight nonabelian simple groups $T$ for which 
the set $\cP_T$ of primes dividing $|T|$ has size at most $3$.
In all cases, the order has the form 
$2^a 3^b p$ and the classical list is as in Table~\ref{23pgroups}. 
References into this classification literature and the complete list
in the much more complicated case $|\cP_T|=4$ are in  \cite{ZCCL09}.
\begin{table}[htb]
{\renewcommand{\arraycolsep}{3pt}
\[
\begin{array}{rcrclll|rcrclll}
p & \GZ & & \!\!\!\!\! |\GZ| \!\!\!\!\! & & H_2 & A & p & \GZ & & \!\!\!\!\! |\GZ| \!\!\!\!\! & & H_2 & A \\
\hline
5 & A_5 &  60 & = & 2^23^15& 2 & 2 & 7 &  SL_3(2) &168&=&2^3 3^1 7&2&2 \\ 
5 & A_6 &  360 & = & 2^3 3^2 5 & 6 & 2^2 & 7 &  SL_2(8)  &504&=&2^3 3^2 7 & 1 & 3 \\  
5 & W(E_6)^+ & 25920 & = & 2^5 3^4 5 & 2 & 2 & 7 &  SU_3(3) &6048& = &2^5 3^3 7 & 1 & 2  \\
\hline
13 &  SL_3(3) &5616&=&2^4 3^3 13  & 1 & 2 & 17 &  PSL_2(17)& 2448 &=& 2^4 3^2 17 & 2 & 2 \\
\end{array}
\]
}
\caption{\label{23pgroups} The eight simple groups of order $2^a 3^b p$ and related information}
\end{table}

The column $H_2$ gives the order of the Schur multiplier $H_2(T)$.  Non-trivial entries here 
are the source of spin separation as explained in the previous section.    
The column $A$ in Table~\ref{23pgroups} gives the structure of the outer automorphism
group of $T$.   So in every case except $T = A_6$ there are two groups $G$ to consider,
$T$ itself and $\Aut(T) = T.A$.   For $T = A_6$ one has, in Atlas order,
the extensions $S_6$, $PGL_2(9)$, and $M_{10}$, as well as the full extension
$\Aut(A_6) = A_6.2^2$.

\subsubsection*{Attractive features of the case $\nu = (3,1)$} 
The restriction $\nu = (3,1)$ is chosen for several reasons.
First, it makes tables much shorter, and in fact 
Tables~\ref{a5tab}, \ref{a6tab}, \ref{g168}, \ref{g504} are complete.   The
case $\nu = (2,2)$ would have similar length
and the case $\nu=(4)$ would 
be even shorter.  However we stay away from both these 
alternatives as the involutions discussed at the end of \S\ref{BP2} complicate
the situation.    
Second, covers  in a 
given degree $m$ tend to 
have considerably smaller genus for
 $\nu = (3,1)$ than
they do for $\nu = (2,1,1)$ or $(1,1,1,1)$.
  In fact our tables show that for $\nu = (3,1)$,
covers can have genus zero into 
quite high degree.   Third,  in the case $\nu = (3,1)$ and $(\beta_0,\beta_1)=(3^{m/3}, \; 2^{m/2})$, 
Beukers and Montanus \cite{BM08} described a method which allows one to solve the given 
system with $m$ unknowns by first solving an auxiliary system with
approximately $m/3$ unknowns.    This method 
generalizes to the full $(3,1)$ case of $(\beta_0,\beta_1) = (3^{a} \, 1^{m-3a}, \, 2^b \,1^{m-2b})$;
we used it simultaneously with the $p$-adic technique sketched in \S\ref{steptwo}
 to extend the reach of our calculations.
Finally,  as discussed in \S\ref{BP2}, the base $\AP^1_j$ is the familiar
$j$-line.   Transitive degree-$m$ covers $\AX_h \rightarrow \AP^1_j$ 
correspond to index-$m$ subgroups of $PSL_2(\Z)$ and we are in 
a very classical setting.

\subsection{Agreement and indexing} 
As discussed in \S\ref{crossparam},
the interesting phenomenon of cross-parameter agreement says
that different Hurwitz parameters
can give rise to isomorphic coverings.   When the two groups
involve different nonabelian simple groups $T$, 
as in the initial example of \S\ref{crossparam}, 
we use the term {\em cross-group} agreement.  
We note cross-group agreement in our tables mainly by
referencing a common equation.  
Two covers appearing even for $T$ involving different primes
are
\begin{eqnarray}
\label{common31} f_{3,1}(j,x) & = & (x-4) x^3 + 4 j  (2 x+1),  \\
\label{common432} f_{4,3,2}(j,x)  &  = &  (4 x^3-3 x+2)^3 -27 j  x^3 (3 x-2)^2. 
\end{eqnarray}
These covers are capable of arising for $|T|$ of the form $2^* 3^* p$ for various
$p$ because their discriminants are respectively $-2^{12} 3^3 j^2 (j-1)^2$ and
$2^{54} 3^{39} j^6 (j-1)^4$.

The covers \eqref{common31}, \eqref{common432} 
illustrate our convention of indexing by the 
braid partition $\beta_\infty$.   This partition and also $\beta_0 = 3^a 1^{m-3a}$ can be 
read off from the presented polynomial.  The remaining partition 
$\beta_1 = 2^b 1^{m-2b}$ governing the factorization of $f_{\beta_\infty}(1,x)$ is
then determined by the fact that we give polynomials only 
in genus zero cases.

\subsection{A degree $46$ map with bad reduction set $\{2,3,13\}$}
The next two sections focus on Hurwitz-Belyi maps coming from groups of order $2^a 3^b p$ for 
$p \in \{5,7\}$.  Here and in the next subsection, to give
a sense of completeness, we give one map each for $p \in \{13,17\}$.  
From $h = (PGL_3(3),(2B,4B),(3,1))$ we get the braid partition triple 
$(\beta_0,\beta_1,\beta_\infty) = (3^{14} 1^{4}, \; 2^{23}, \;13^2 \, 8 \, 6 \, 3 \, 2 \,1)$.  Our final polynomial is
 \begin{eqnarray*}
\lefteqn{ f_{13^2,8,6,3,2,1}(j,x) = } \\
&&  \left(16 x^4+40 x^3-3 x^2-116 x-8\right)   
 \left(4096 x^{14}+20480 x^{13}-25856 x^{12}- \right. \\ && \qquad 196736 x^{11}+47189 x^{10}+680764 x^9-69384 x^8-1135104 x^7+7638144 x^6 \\ && 
 \qquad \left. -16337408
   x^5+9620480 x^4-2785280 x^3+741376 x^2-16384 x-32768\right)^3 \\ && -2^{13} 3^{12} j (x-2)^3 x^6 (x+4)^2 (2 x-1) \left(3 x^2+2 x-4\right)^{13}.
 \end{eqnarray*}
 The discriminant of this polynomial is $- 2^{2260} 3^{1371} 13^{351} (j-1)^{23} j^{28}$.  
 Modulo squares this discriminant is $-39 (j-1)$.  The factor of $j-1$ is known
 from the outset by the oddness of $\beta_1$ and $\beta_\infty$.  
  
 \subsection{A degree $54$ map with bad reduction set $\{2,3,17\}$}
  The Hurwitz parameters
 \begin{equation}
 h_1  =  (SL_2(17),(17a+,3a-),(3,1)) \mbox{ and }
 h_2  =  (PGL_2(17),(2B,6A),(3,1)) 
 \end{equation}
 each give conjugate braid permutation triples, with common braid partition triple 
 \begin{equation}
 (\beta_0,\beta_1,\beta_\infty) = (3^{17} \, 1^3, \; 2^{27}, \; 9^3 \, 8^2 \,  4^2 \, 2 \, 1).
 \end{equation}
An equation was determined by Schiavone using improvements of the techniques described in \cite{KMSV14}:
 \begin{eqnarray}
 \nonumber \lefteqn{ f_{9^3,8^2,4^2,2,1}(j,x) = } \\
\nonumber && \left(x^3+12 x^2+12 x-8\right) \cdot  \left(x^{17}-52 x^{16}+42136 x^{15}-593008
   x^{14}+10147846 x^{13} \right. \\ 
   \nonumber && \left. \qquad +225862160 x^{12}+1467000268 x^{11}+6342760760 x^{10}+593082769
   x^9 \right. \\
   \label{prime17} && \left. \qquad -1815237116 x^8-5586407260 x^7-258348008 x^6+8975722736 x^5
   \right. \\ 
   \nonumber &&  \qquad \left. -8292246656
   x^4+3424464320 x^3-664160384 x^2+44883968 x-131072\right)^3\\
  \nonumber &&-  2^4 3^9 j x \left(x^2-71
   x+32\right)^4 \left(x^2+2 x-1\right)^8 \left(x^3+18 x^2-48 x-8\right)^9.
  \end{eqnarray}
  It seems that this cross-parameter agreement is one of an infinite family 
  indexed by odd primes as follows. Generalize $h_1$ to $(SL_2(p),(pa+,3a-),(3,1))$.  
  Generalize $h_2$ to $(PGL_2(p),(2B,6A),(3,1))$ when $p \equiv \pm 5 \; (12)$ and
  to $(PSL_2(p),(2b,6a),(3,1))$ when $p \equiv \pm 1 \; (12)$.  Then 
  mass computations confirm that both covers have degree
  \[
  m = \frac{p^2-5}{4} + \left\{ \begin{array}{rl} p & \mbox{ if $p \equiv 1 \; (3)$} \\ -p & \mbox{ if $p \equiv 2 \; (3)$} \end{array} \right. . 
  \]
  Braid computations say that indeed the covers are isomorphic at least for $p \leq 19$.    
  For $m=5$ and $7$ the degrees are $0$ and $18$ respectively, these cases arising in \S\ref{235first}
  and \S\ref{237first}.

\section{Hurwitz-Belyi maps with $|G|=2^a 3^b 5$ and $\nu=(3,1)$}
\label{235} 
   In this section we work in the framework set up in \S\ref{23poverview} and 
   present a systematic collection of 
Hurwitz-Belyi maps having bad reduction
at exactly $\{2,3,5\}$.  

\subsubsection*{Cross-group agreement}  Before getting to the individual groups, we present
equations for covers involved in cross-group agreement.  The covers
\begin{align}
\label{common531} f_{5,3,1}(j,x) & =  5^2 \left(5 x^3-45 x^2+39 x+25\right)^3-2^{14} 3^3 j x^3 (3 x-25),   \\
\label{common532} f_{5,3,2}(j,x) & =  \left(9 x^3+3 x^2-53 x+81\right)^3 (x+9) -2^{14} 3^2 j x^3 (3 x-5)^2,  \\
\label{common543} f_{5,4,3}(j,x) & = \left(4 x^4-24 x^3+24 x^2-48 x+27\right)^3 -2^2 3^3 j x^3 (3  x-4)^5  
\end{align}
appear for all three groups.     The covers
\begin{align}
\label{common541} f_{5,4,1}(j,x) & =  \left(16 x^3-87 x^2+48 x+16\right)^3 (16 x+1) -2^23^{12} j  x^4 (x-5),\\
\label{common5433} f_{5,4,3^2}(j,x) & =  4 \left(256 x^5+640 x^4-440 x^3-3325 x^2-6400
    x-4096\right)^3 &    \\ 
    \nonumber & \qquad  - 3^{12} j x^4 \left(32 x^2+95 x+80\right)^3   
\end{align}
appear for the simple groups $A_5$ and $A_6$.
The covers 
\begin{align}
\label{common51} f_{5,1}(j,x) & = \left(x^2-5\right)^3 -3^3 j  (2 x-5),  \\
\label{common5542} f_{5^2,4,2}(j,x)  & = 2^7 x \left(18 x^5-144 x^4+336 x^3-224 x^2+801 x-162\right)^3 \\
\nonumber & \qquad  - j (2 x-9)^2 \left(36 x^2-52 x-9\right)^5,      
\end{align}
appear for $A_6$ and $W(E_6)^+$.  Several larger degree covers also appear for 
both $A_6$ and $W(E_6)^+$.  A polynomial for the
smallest of these is 
\begin{eqnarray}
\nonumber \lefteqn{f_{10,8^2,6,5,4^2,3}(j,x) =} \\
\label{common108865443} &  &  \left(184528125 x^{16}-984150000 x^{15}+2263545000 x^{14} -2768742000
    x^{13}  \right. \\ 
    \nonumber && \left. 
    +1616849100 x^{12}+181316880 x^{11}-1023304104 x^{10}+721510416
    x^9  \right. \\ 
    \nonumber && \left.  -166620402 x^8-72763728 x^7+59318552 x^6-4952016 x^5-12051828
    x^4 \right. \\ 
    \nonumber && \left. +7406640 x^3-2117016 x^2+314928 x-19683\right)^3  \\
\nonumber && + 2^{20} \, 3^8 \, j \left(9 x^2-10   x+3\right)^8 x^6  (5 x-3)^5 \left(3 x^2-1\right)^4 (3 x-1)^3 .
    \end{eqnarray}

\subsection{The simple group $A_5 \cong SL_2(4) \cong  PSL_2(5)$}
\label{235first}
Tables~\ref{a5tab}, \ref{a6tab}, \ref{g168}, and \ref{g504} have a similar structure, which we explain 
now drawing on Table~\ref{a5tab} where $T=A_5$ for examples.    The top left subtable gives degrees of components of Hurwitz-Belyi maps
$\AX^*_h \rightarrow \AP^1_j$ for $h = (T,(C_1,C_2),(3,1))$.   Here $C_1$ and $C_2$ are distinct conjugacy
classes in $T$.   When the lifting invariant set $H_h^*$ from \S\ref{liftingingeneral} is trivial,  a single number is typically printed.  For $T=A_5$,
this triviality occurs exactly if $221$ is one of the $C_i$, as from Table~\ref{chartable} for $A_5$, only $221$ is inert
in the double cover $\tilde{A}_5$.  When $H_h^*$ is canonically $\Z/2$, typically two numbers are printed; the top 
and bottom numbers respectively give the degrees of $\AX_h^{*+}$  and $\AX_h^{*-}$ over $\AP^1_j$.  

The remaining subtables on the left sides of Tables~\ref{a5tab}, \ref{a6tab}, \ref{g168}, and \ref{g504} similarly
give degrees of components of Hurwitz-Belyi maps, bit now for $h = (T.2,(C_1,C_2),(3,1))$.  
When $H_h=H_h^*$ has $2$ elements but is not canonically $\Z/2$, typically again two numbers are printed in a column.  These
numbers are necessarily the same.  In general, a number is put in italics when the corresponding
component is not defined over $\Q$.  For $G=S_5$, irrationality occurs exactly if $\{C_1,C_2\} = \{41,32\}$, 
a key point being that $41$ and $32$ each split into two irrational classes in $\tilde{S}_5$.   Other
possibilities for $H_h^*$ occur only for $T=A_6$ and will be discussed in \S\ref{235second}.   
In general, if there is splitting beyond that forced by lifting invariants then the corresponding
degree is written as a list of the component degrees separated by commas.   This extra splitting
does not occur on Table~\ref{a5tab} and we expect it to be rare in general.  Indeed for $G=A_5$ 
and any $(C,\nu)$, it never occurs on the level of the entire $r$-dimensional Hurwitz cover 
$\AHur^*_h \rightarrow \AConf_\nu$ \cite{JMS15}.  

\begin{table}[htb]
\[
\begin{array}{|c|c|c|c|}
\hline
\mbox{\small{$\! C_1 \!\! \setminus \!\! C_2 \! $}} & 221 & 311 &  5b \\
\hline
221 & \bullet & 0 & 10a \\
\hline
\multirow{2}{*}{311} & \multirow{2}{*}{12} & \multirow{2}{*}{$\bullet$} & 15 \\
&&& 10b \\
\hline
\multirow{2}{*}{$5a$} & \multirow{2}{*}{4} & 9  & 4 \\
& & 0 & 0 \\
\hline
\multicolumn{4}{c}{\;} \\
\hline
\mbox{\small{$\! C_1 \!\! \setminus \!\! C_2 \! $}}  & 2111 & 41 & 32 \\
\hline
2111 & \bullet & 0 & 0 \\
\hline
\multirow{2}{*}{41} & \multirow{2}{*}{32} &  \multirow{2}{*}{$\bullet$} &  \mathit{36} \\
&&&\mathit{36} \\
\hline
\multirow{2}{*}{32} & \multirow{2}{*}{$10b$} & \mathit{16} &  \multirow{2}{*}{$\bullet$}  \\
&& \mathit{16} & \\
\hline
\end{array}
\;\;\;\;
 \begin{array}{ |  lllllll |  }
 \hline
                 \# &  M &
                  g & \beta_0 & \beta_1 & \beta_\infty & \mbox{Eqn.}    \\
                \hline
                              2+&   A_4 & 0 &  3\; 1 & 2^2 & 3\;1 & \eqref{common31} \\
  1+& A_9 & 0 &  3^3 & 2^41 & 5\; 3\; 1 &  \eqref{common531} \\
                                                                                        1+& S_{10a} & 0 &  3^3 1 & 2^5 & 5\; 3\; 2 &  \eqref{common532}   \\   
                                 2+& S_{10b} & 0&  3^3 1 & 2^5 & 5\; 4\; 1  &  \eqref{common541} \\                           
                               1+&  S_{12} & 0 &  3^4 & 2^5 1^2 & 5\; 4\; 3 &     \eqref{common543}  \\
                       1+& S_{15} & 0 &  3^5 & 3^7 1 & 5\; 4\; 3^2 &  \eqref{common5433} \\ 
                   1& A_{32}& 0&  3^{10} 1^2 & 2^{16} & 10\; 6\; 5\; 4^2 \, 3   &   \eqref{cover1065443} \\
                         \hline
                    \multicolumn{7}{c}{\;} \\
                    \multicolumn{7}{c}{\;} \\
                    \hline
                    \multicolumn{7}{|c|}{\mbox{Two pairs defined over $\Q(\sqrt{6})$}} \\
                    \hline
                    1 & A_{16} & 0 & 3^5 \, 1 & 2^{8} & 6 \; 5 \; 4 \; 1 & \eqref{cover6541} \\
                    1 & A_{36} & 0 & 3^{12}  & 2^{18} & 10 \; 6 \; 5 \; 4^2 \, 3^2 \, 1 &  \\
                    \hline
  \end{array}
\]
\caption{\label{a5tab}  Left: Degrees of components of Hurwitz-Belyi covers with
parameters $(G,(C_1,C_2),(3,1))$ with $G = A_5$ or $S_5$.  Right: further information 
on these covers.}
\end{table}

We are interested primarily in rational covers and we distinguish non-isomorphic 
rational covers of the same degree by identifying labels.   This convention
highlights cross-parameter agreement.  Thus on the left half of Table~\ref{a5tab}
the two 4's and the two $10b$'s each represent isomorphic covers.

The left half of Tables~\ref{a5tab} \ref{a6tab}, \ref{g168}, or  \ref{g504},  as just described, is well thought
of as the Hurwitz half.  The right half can then be considered the Belyi half, as it
makes no reference to its Hurwitz sources beyond the column $\#$.  Here a 
number printed under $\#$ just repeats the number of Hurwitz sources from
the left half; a $+$ sign represents cross-group agreement, as it 
indicates that the cover also arises elsewhere in this paper for a different $T$.  
While our focus is on Hurwitz-Belyi covers defined over $\Q$, when there is
space we include extra lines for Hurwitz-Belyi covers not defined over $\Q$.

Equations for the first six lines of the top right subtable of Table~\ref{a5tab} have already been presented in the context of
cross-group agreement.   An equation for the seventh line is
\begin{align}
\nonumber f_{10,6,5,4^2,3}(j,x) & = \left(x^{10}-38 x^9+591 x^8-4920 x^7+24050
    x^6-71236 x^5+125638 x^4 \right. 
    \\
\label{cover1065443}   & \left. \qquad -124536 x^3+40365 x^2+85050
    x-91125\right)^3 \left(x^2-14 x-5\right) \\
\nonumber    & \;\;\; + 2^{20} 3^3 j x^6 (x-5)^5  \left(x^2-4    x+5\right)^4 (x-9)^3 . 
\end{align}
Note that the four-point covers $\AY_x \rightarrow \AP^1_t$ corresponding
to the seventh line have genus one, and so \eqref{cover1065443} would be hard to compute by the standard
method.   The tables of this and the next section give many examples where
$g_Y$ is large but $g_X=0$.   As we are systematically using the braid-triple method, 
$g_Y$ is irrelevant and the tables present $g = g_X$.

Tables~\ref{a5tab}, \ref{a6tab}, \ref{g168}, and \ref{g504} exclude the case $C_1=C_2$ to stay 
in the context of $\nu=(3,1)$.  The excluded cases $(G,C_1,(4))$ are interesting too
and we mention one of them.  For $h = (A_5,(311),(4))$, the cover $\AX_h^{*+}$ 
is given by $f_{5,3,1}(j,x)$ from \eqref{common531} while $\AX^{*-}_h$ is empty. 
This $h$ is our first of three illustrations of a general theorem of Serre \cite{Ser90} as follows.  Consider
Hurwitz parameters 
\begin{equation*}
h = (A_n,(e_1 1^{n-e_1},\dots,e_k 1^{n-e_k}),(\nu_1,\dots,\nu_k))
\end{equation*}
 with all $e_i$ odd,
so that one has a lifting invariant and thus an equation  $\AX^*_h = \AX_h^{*+} \coprod \AX_h^{*-}$.
Suppose $\sum \nu_i (e_i-1) = 2n-2$ so that the genus $g_Y$ is $0$.   Then the general theorem
says,
\begin{equation}
\label{serretheorem}
\mbox{
{\em If $\prod e_i^{\nu_i} \stackrel{8}{\equiv} \left\{ \begin{array}{l} \pm 1 \\ \pm 3 \end{array} \right.$ then 
$\AX^*_h = \left\{
\begin{array}{l}
\AX^{*+}_{h}  \\
\AX^{*-}_{h}
\end{array}
\right.$.}}
\end{equation}
Table~\ref{a5tab} shows that $\AX_h^{*-}$ is empty for   
$(5a,311)$ and $(5a,5b)$ as well,
even though $g_Y>0$ and so Serre's theorem does not apply in these cases.

The left half of Figure~\ref{g168} refers to two pairs of irrational covers.  Letting $s = \pm \sqrt{6}$,
equations for the smaller degree pair are
\begin{eqnarray}
\nonumber \lefteqn{f_{6,5,4,1}(j,x) =} \\
\nonumber && (3 x+s-3) \left(-225 x^5+1305 x^4 s+4005 x^4-8932 x^3 s-22662 x^3+6594 x^2 s \right. \\
\label{cover6541}  && \qquad \left. +16254 x^2-28476
    x s-69741 x+11673 s+28593\right)^3 \\
\nonumber && - 12288 j x^5 (5 x-9)^4 (53236 s+130401) (-15 x+76 s+186).
\end{eqnarray}

\subsection{The simple group $A_6 \cong Sp_4(2)' \cong  PSL_2(9)$}
\label{235second}
  In terms of both its Schur multiplier
$H_2 \cong \Z/6$ and its outer automorphism group $A \cong (\Z/2)^2$, the group
$T=A_6$ is the most complicated group on Table~\ref{23pgroups}.  
Table~\ref{sexticclasses} gives information on conjugacy classes.

\begin{table}[htb]
\[
\begin{array}{|cccccc|ccccc|}
\hline
\multicolumn{3}{|r}{H_2(A_6)} &= & \multicolumn{2}{l|}{6}&
\multicolumn{2}{|r}{H_2(S_6)} &= & \multicolumn{2}{l|}{2} \\
\multicolumn{3}{|r}{\Out(A_6)} &= & \multicolumn{2}{l|}{2^2}&
\multicolumn{2}{|r}{\Out(S_6)} &= & \multicolumn{2}{l|}{2} \\
\hline
2211 & 3111 & 33 & 42 & 5a & 5b & 21111 & 222 & 411 & 6 & 321 \\
3 & 2 & 2 & 6 & 6 & 6 & 1 & 1 & 1 & 2 & 2 \\
\hline
\end{array}
\]

\[
\begin{array}{|ccccc|ccc|}
\hline
\multicolumn{3}{|r}{H_2(PGL_2(9))} &= & \multicolumn{1}{l|}{3}&
\multicolumn{1}{|r}{H_2(M_{10})} &= & \multicolumn{1}{l|}{3} \\
\multicolumn{3}{|r}{\Out(PGL_2(9))} &= & \multicolumn{1}{l|}{2}&
\multicolumn{1}{|r}{\Out(M_{10})} &= & \multicolumn{1}{l|}{2} \\
\hline
2222 & 811A & 811B & \, 10A \,  & \, 10B \, & 4411 & 82C & 82D \\
1 & 2 & 2 & 2 & 2 & 3 & 3 & 3 \\
\hline
\end{array}
\]
\caption{\label{sexticclasses} Information on conjugacy classes in $A_6$ and conjugacy 
classes not in $A_6$ of its three
overgroups $S_6 = A_6.2_1$, $PGL_2(9) = A_6.2_2$, and $M_{10} = A_6.2_3$.  The
last row gives the number of classes in $\tilde{G}$ mapping to the
given class in $G$.}
\end{table}

\begin{table}[htb]
\[
{\renewcommand{\arraycolsep}{0.7pt}
\begin{array}{|c|c|c|c|c|}
\hline
\mbox{\small{$\! C_1 \!\! \setminus \!\! C_2 \! $}} & 2211 & 33  & 42 &  51b \\
\hline
2211 & \bullet & 0 &  12 & 10a; \! 15 \\
\hline
\multirow{2}{*}{$3111$}  & \multirow{2}{*}{12} &0 &  \multirow{2}{*}{16} &  10a \\
&&6&&0 \\
\hline
\multirow{2}{*}{$42$} & \multirow{2}{*}{$9;48$} 
& \multirow{2}{*}{$108$} & \multirow{2}{*}{$\bullet$}  & \mathit{60}; \! \mathit{40} \\
&&&&\mathit{60}; \! \mathit{40} \\
\hline
\multirow{2}{*}{$51a$} & \multirow{2}{*}{$60;0$} &  45 & \mathit{24};\mathit{40} & 9;45 \\
&&54 &\mathit{24};\mathit{40}&0;42\\
\hline
\multicolumn{5}{c}{\;} \\[.5ex]
\hline
\mbox{\small{$\! C_1 \!\! \setminus \!\! C_2 \! $}} & 21111 & 411  & 321 &  6 \\
\hline
21111 & \bullet & 0 & 0 & 0 \\
\hline
411 & 10a & \bullet & 252 &  \\
\hline
\multirow{2}{*}{$321$} & \multirow{2}{*}{$84$} & \multirow{2}{*}{$544$} & \multirow{2}{*}{$\bullet$} & \mathit{396} \\
&&&& \mathit{396} \\
\hline
6 & 88 &  &  & \bullet  \\
\hline
\multicolumn{5}{c}{\;} \\[.5ex]
\cline{1-4}
\mbox{\small{$\! C_1 \!\! \setminus \!\! C_2 \! $}} & 22222 & 811B  & 10B \\
\cline{1-4}
22222 & \bullet & 16 & 10b \\
\cline{1-4}
\multirow{2}{*}{$811A$} & \multirow{2}{*}{$96a$}  & 96b & \mathit{100}\\
&&108 & \mathit{100}\\
\cline{1-4}
\multirow{2}{*}{$10A$} & \multirow{2}{*}{$4,42$} & \mathit{64}& 4,42  \\
                                     &                                  &\mathit{64}& 54 \\
\cline{1-4}
\multicolumn{5}{c}{\;} \\[.5ex]
\hline
\mbox{\small{$\! C_1 \!\! \setminus \!\! C_2 \! $}} & \multicolumn{2}{c|}{4411} & \multicolumn{2}{c|}{82D}  \\
\hline
4411 & \multicolumn{2}{c|}{\bullet} & \multicolumn{2}{c|}{656; \! \mathit{672}; \! \mathit{672}}   \\
\hline
82C &   \multicolumn{2}{c|}{164; \! \mathit{168}; \! \mathit{168}}          & \multicolumn{2}{c|}{66;\mathit{90};\mathit{90}}  \\
\hline
\end{array}
\;\;\;
 \begin{array}{ | l  lllll  l  | }
 \hline
                 \# & 
                 M &  g \;\; & \beta_0 & \beta_1 & \beta_\infty & \mbox{Eqn.}   \\
                 \hline
                           2+& A_4 &      0  & 3 \; 1 & 2^2 & 3 \; 1 &   \eqref{common31} \\
                           1+&A_6 & 0 & 3^2 & 2^2 1^2 & 5 \;1 & \eqref{common51}  \\
                            2+&A_{9} & 0 &  3^3 & 2^4 \; 1 & 5 \; 3 \; 1 & \eqref{common531} \\                      
                            3+&S_{10a} & 0 &  3^3 1 & 2^5 & 5 \; 3 \; 2 &  \eqref{common532} \\
                             1+& S_{10b} & 0  & 3^3 \, 1 & 2^5 & 5 \; 4 \; 1 &      \eqref{common541} \\
                             2+&S_{12} & 0 & 3^4 & 2^{5} \, 1^2 & 5 \; 4 \; 3 &  \eqref{common543}  \\
                    1+&S_{15} & 0  & 3^5 & 2^7 1 & 5 \; 4 \; 3^2 & \eqref{common5433}  \\
   2+& A_{16} & 0 & 3^5 \, 1 & 2^8 & 5^2 \, 4 \; 2 &    \eqref{common5542} \\
        3&S_{42} & 0&  3^{14} & 2^{21} & 10\; 8^2 \,  4^2 \,  3^2 \,  1^2&   \\
          2&A_{45} &
                    0&  3^{15} & 2^{22} \,  1 & 8^2 \, 5^3 \,  4^2 \, 3^2  &    \\
                      1+& A_{48} & 0 & 3^{16} & 2^{22} \, 1^4& 10 \;  8^2 \,  6 \, 5 \, 4^2 \,  3 & \eqref{common108865443}  \\
                             2 &S_{54} & 0& 3^{18} & 2^{27} & 10\; 8^2 \,  5^2 \, 4^2 \, 3^3 \, 1 & \\
                                              1& S_{60} & 0&  3^{20} & 2^{29} \, 1^2 &  10 \; 8^2 \, 5^4 \, 4^2 \, 3^2 &  \\
                                           1&  S_{66} & 0 & 3^{21} \, 1^3 & 2^{33} & 10^2 \, 8^3 \, 6^3 \, 2^1 \, 1^2 & \\
                  1+& A_{84} & 0 &  3^{27} \,  1^3 & 2^{42} & 10^3 \, 8^2 \, 6^3 \,  5\; 4^3 \,  2\; 1  &    \\
                    1&A_{88} & 0 & 3^{28} \,  1^4 & 2^{44} & 10^3 \,  8^4 \, 6\; 5^2 \, 3^3 \, 1 &   \\
                                    1& A_{96a} & 1 &  3^{32} & 2^{48} & 10^2 \, 8^3 \, 6^3 \, 5^4 \, 4^2 \, 3^2 &       \\
                                    1& G_{96b} & 1 &  3^{32} & 2^{48} & 10^2 \, 8^3 \, 6^3 \, 5^4 \, 4^2 \, 3^2 &       \\
                     2+ &A_{108} & 0& 3^{36} & 2^{54} & 10^2 \, 8^4 \,  6^2 \, 5^2 \,  4^4 \,  3^2  &  \\
                        1& A_{164} & 0 &  3^{53} \, 1 & 2^{82} & 10^6 \, 8^5 \, 6^5 \, 4^8 \, 1^2 &    \\
                       1& A_{252}  &  0  & 3^{84} & 2^{126}  & \multicolumn{2}{l|}{10^{4} \,  8^9  \, 6^6 \,  5^{10} \, 4^8 \,  3^5 \, 1}    \\
                   1& A_{544} & 3 & 3^{181} & 2^{272} &  \multicolumn{2}{l|}{10^{18} \,  8^{18} \, 6^{15} \, 5^{10} \,  4^{13} \,  3^8 \,  1^4}  \\        
                  1& A_{656} & 7    &  3^{217} \, 1^5 & 2^{328} & 10^{26} \, 8^{20} \, 6^{25} \, 4^{21} \, 1^2 &     \\
                                  \hline
                                    \multicolumn{7}{c}{\;} \\[-.7ex]
                    \hline 
                    \multicolumn{7}{|c|}{\mbox{A pair defined over $\Q(\sqrt{10})$}} \\
                    \hline
                    1 & A_{24} & 0 & 3^8  & 2^{12} & 8 \; 5 \; 4 \; 3^2 \, 1 & \\
                    \hline
                  \end{array}
                  }
\]
\caption{\label{a6tab}  Left: Degrees of components of Hurwitz-Belyi covers with
parameters $(G,(C_1,C_2),(3,1))$ with $G = A_6$, $S_6$, $PGL_2(9)$, or $M_{10}$.
  Right: further information 
on these covers.}
\end{table}

Conventions about the $(C_1,C_2)$ entry in the left half of Table~\ref{a6tab} have been given in \S\ref{235first}
whenever $|H_2(G,C)| \in \{1,2\}$.  The remaining possibilities are as follows.   Three entries
in a single row separated by semicolons means $|H_h|=3$ and $\Out(G,C)$ acts trivially on $H_h$, so that
$|H_h^*|=3$ as well.  This possibility arises three times, always in the form $(\mbox{a};b;b)$.  
By the typeface convention of \S\ref{235first}, this means that the degree $\mbox{a}$ component is rational 
and the degree $b$ components are conjugate.   Two entries in a single row separated
by semicolons means $|H_h| = 3$ but $\Out(G,C)$ acts nontrivially on $H_h$,
so that $|H^*_h| = 2$.   This possibility also arises three times, always in the form $(\mbox{c};\mbox{d})$.
Here both components are rational, as indeed in these three cases $\mbox{c} \neq \mbox{d}$.  
Instances of these two situations were described already in \S\ref{cubic2} and \S\ref{cubic}, where degrees were $(\mbox{a};b;b) = (0;21;21)$ and $(\mbox{c},\mbox{d})=(30;40)$ respectively.    In these
situations, one generally expects $\mbox{a} \approx b$ and $\mbox{c} \approx \mbox{d}/2$. 

In the case $(5a,5b)$, one has $H_h \cong \Z/6$ and
$\Out(G,C)$ has order two.  The non-trivial element of $\Out(G,C)$ acts by 
negation, so that $H_h^*$ has order four.  The natural action of $\Gal(\overline{\Q}/\Q)$ on 
$H_h$ is trivial,
and one would generally expect four rational components.  In this case,
the natural map $\pi_0(\AX^*_h) \rightarrow H_h^*$ is injective but not
surjective, and $\AX^*_h$ has only three components.    
The cases $(42,51a)$ and $(51b,42)$ are similar to $(5a,5b)$ but now
all components are defined over $\Q(\sqrt{10})$.  The two degree $24$ 
components have their dessin drawn in the website associated to \cite{BM08}.

A blank in the $(C_1,C_2)$ slot  
means that 
covers belonging to this slot are isomorphic to those of $(C^\alpha_1,C^\alpha_2)$ for
some $\alpha$ in $\Out(G)-\Out(G,C)$.  For example the $(411,6)$ slot is 
left blank because the cover is the same as that represented by the $(411,321)$ slot.  
It is this non-triviality of $\Out(G)-\Out(G,C)$ that makes some of the covers involving $51a$ and/or $51b$ 
rational, even though the classes $51a$ and $51b$ are conjugate to each other.  
Among the further  things to note on Table~\ref{a6tab} are two isomorphic unforced decompositions
of the form $46=42+4$.   Also the cover $96b$ is unexpectedly nonfull.  Finally,
a second instance of Serre's theorem \eqref{serretheorem} is $(C_1,C_2) = (3111,51a)$, so
that  $\AX^{*-}_h$ is forced to be empty.  

\subsection{The simple group $W(E_6)^+ \cong PSp_4(3) \cong PSU_4(2)$} \label{235third}
 The group $W(E_6) = W(E_6)^+.2$ has twenty-five conjugacy classes.
As for all Weyl groups, all the classes are rational. 
\begin{table}[htb]
  \[
 \begin{array}{ | ll |  lllll  l  | }
 \hline
                 C_1 & C_2 & M &  g  & \beta_0 & \beta_1 & \beta_\infty &  \mbox{Eqn.}    \\
                \hline
\multicolumn{8}{c}{\;} \\ 
 \multicolumn{8}{c}{\mbox{Full low degree covers}} \\      
 \hline
                         3d & 4a  & A_6 & 0 & 3^2 & 2^2 \, 1^2 & 5 \; 1 & \eqref{common51} \\
                        3C & 9A & A_9 & 0 &  3^3 & 2^4 \, 1 & 5 \; 3 \; 1   &  \eqref{common531} \\
                    &&   S_9 & 0      &   3^3      &  2^3 \, 1^3         &  5 \; 4      & \eqref{cover54}              \\
                      &&    S_{10}      &   0      &      3^3 \, 1  &    2^5    &  5 \, 3 \, 2     &        \eqref{common532}      \\
&&   S_{12}     & 0  &   3^4      &    2^5 \, 1^2     &   5 \; 4 \; 3     &              \eqref{common543}      \\
6a&4a&   A_{16}     & 0   &   3^5 \, 1         &   2^8      &   5^2 \, 4 \, 2      &         \eqref{common5542}           \\
6a&2b&   A_{24}      &    0     &    3^8     &  2^{10} \, 1^4       &  9 \; 6 \; 5\; 4     &        \eqref{cover9654}      \\
               3c& 9a& S_{27} & 0 &  3^{9} & 2^{13} \, 1 & 5^4 \, 4 \; 3   &  \eqref{cover555543} \\
                           4a&2b& S_{28} & 0 &  3^9 \, 1 &  2^{13} \, 1^2 & 10 \; 9 \; 5 \; 2^2 & \eqref{cover109522}   \\
                              \hline 
\multicolumn{8}{c}{\;} \\ 
 \multicolumn{8}{c}{\mbox{Agreement with covers coming from $A_6$}} \\ 
 \hline
                       4a&3c 1& A_{48} & 0 & 3^{16} &  2^{22} \,1^4 & 10 \; 8^2 \; 6 \; 5 \; 4^2 \, 3 & \eqref{common108865443}  \\
                       6e  &6a        & A_{84} & 0 & 3^{27} \, 1^3 & 2^{42} & 10^3 \, 8^2 \, 6^3 \, 5 \; 4^3 \, 2 \, 1 &  \\
                       4a  &6a       & A_{108} & 0 &  3^{36} & 2^{54} & 10^2 \, 8^3 \, 6^3 \, 5^4 \, 4^4 \, 3^3 \, 1 &   \\
\hline 
\multicolumn{8}{c}{\;} \\ 
 \multicolumn{8}{c}{\mbox{Large degree genus zero examples of spin separation}} \\                                                                                                               
\hline
 3d&5a& A_{165} &  0 & 3^{55} & 2^{80} 1^5 & 12^4 \, 9^4 \, 6^7 \, 5^2 \, 4^2 \, 3^6 \, 2^1 \, 1 &    \\
 3d&5a& S_{225} &  0 & 3^{75} & 2^{109} \, 1^{7} & 12^4 \; 9^7 \, 6^7 \, 5^7 \, 4^6 \, 3^3 \, 2^2 &   \\                                                        
\hline
3d&9b&S_{189} & 0 & 3^{63} & 2^{93} 1^3 & 12^4 \, 9^6 \, 6^7 \, 5^4 \, 4^4 \, 3^9 &    \\
3d&9b &S_{234} & 0 & 3^{78} & 2^{117} & 12^4 \, 9^6 \,  6^7 \, 5^7 \,  4^7 \, 3^8 \,  2 \; 1&  \\
\hline
                                      \multicolumn{8}{c}{\;} \\                                     
                                      \multicolumn{8}{c}{\mbox{An example where all allowed bases appear in $\beta_\infty$}} \\
                                      \hline
                               4D&6G&S_{3186} & 82 & 3^{1062}  & 2^{1593} &
                                      \multicolumn{2}{l|}{  24^6 \, 18^{27} \, 12^{43} \, 10^{78} \, 9^{27} \, 8^{54} }   \\    
                               && & &   & &
                                      \multicolumn{2}{l|}{  \qquad 6^{24} \, 5^{48} \, 
                                       4^{28} \, 3^{23} \, 2^9 \, 1^2 }  \\    
                                      \hline               
 \end{array}
 \]  
\caption{\label{w6tab} Invariants of some covers with $G = W(E_6)^+$ or $W(E_6)$}
\end{table}
Ten classes are in $W(E_6)-W(E_6)^+$ and ten
classes stay rational classes in $W(E_6)^+$.   The remaining five conjugacy classes of $W(E_6)$, namely
$3ab$, $6ab$, $6cd$, $9ab$, and $12ab$, split into two classes in $W(E_6)^+$.  
  If we were presenting complete tables for $\nu = (1,1,1,1)$, there would thousands of lines.    Even complete
  tables for $(3,1)$ would have hundreds of lines.  Accordingly, Table~\ref{w6tab} presents 
  just some of Hurwitz-Belyi covers in a self-explanatory format.   
  
  One of the new covers has the remarkably small degree nine:
\begin{equation}
\label{cover54} f_{5,4}(j,x) = 5^2 \left(10 x^3+15 x^2+48 x-100\right)^3 + 3^{15} j  x^4.
\end{equation}
The other three new covers are
\begin{eqnarray}
\nonumber \lefteqn{f_{9,6,5,4}(j,x)} \\
 \nonumber \;\;\;\;\; & = &  \left(9 x^8-72 x^7+180 x^6-104 x^5-26 x^4-568 x^3 +1620 x^2-1944
    x+729\right)^3 \\ 
\label{cover9654}  &&  + 2^{16} j (x-3)^4 x^6 (2 x-3)^5, \\
\nonumber &&\\
\nonumber  \lefteqn{f_{5^4,4,3}(j,x) }  \\ 
\nonumber &=& (1024 x^9-13824 x^8+81360 x^7-272928 x^6+585144 x^5-879336
    x^4 \\ 
\nonumber    && \qquad \qquad+1012365 x^3 -896832 x^2+516096 x-131072)^3 \\ 
\label{cover555543}    && - 54 j (72 x^4-508 x^3+1350 x^2-1629 x+768)^5 x^3, \\
\nonumber &&\\
\nonumber  \lefteqn{f_{10,9,5,2^2}(j,x)} \\
\nonumber &=&  \left(3125 x^9-9375 x^8+7500 x^7-6500 x^6+9150 x^5-4410 x^4-2484
    x^3 \right. \\ && \nonumber  \qquad \left. -2916 x^2-2187 x+6561\right)^3 (x-3) \\
\label{cover109522}     && + 2^{22} 3^3 j x^9 (5 x-6)^5 \left(3 x^2+2 x+3\right)^2.
\end{eqnarray}
In the entire table for $T=W(E_6)^+$, there are only twelve integers which can appear
as parts for $\beta_\infty$.  The last line of Table~\ref{w6tab} gives the smallest degree
cover where all these integers actually appear.

\section{Hurwitz-Belyi maps with $|G|=2^a 3^b 7$ and $\nu=(3,1)$}
\label{237}
This section is very parallel in structure to the previous one, and presents 
a systematic collection of Hurwitz-Belyi maps having bad reduction
at exactly $\{2,3,7\}$.

\subsubsection*{Cross-group agreement}  Again we present equations for covers involved in cross-group
agreement before getting to the individual groups.  Now we have
only two:  
\begin{eqnarray}
\label{common43} f_{4,3}(j,x) & = & 4 (x-12) \left(9 x^2-20 x-27\right)^3 + 3 \cdot 7^7 j  x^3, \\
\label{common74331} f_{7,4,3^2,1}(j,x) \!\! &=& \!\!  \left(9 x^6-126 x^4+252 x^3-63 x^2-252 x+196\right)^3   \\ 
 \nonumber &&  +  2^6 j (3 x-2)^4  \left(3 x^2-9 x+7\right)^3 (3 x + 14) . 
\end{eqnarray}
The cover $f_{7,4,3^2,1}(j,x)$ was first found by Malle \cite{Mal00} in connection with 
the group $PGL_2(7)$.

\subsection{The simple group $PSL_2(7) \cong SL_3(2)$}
\label{237first}
Equations for the first three Hurwitz-Belyi maps have been given already.  For the fourth, an equation is
\begin{table}[htb]
\[
{\renewcommand{\arraycolsep}{3pt}
\begin{array}{|c|c|c|c|c|}
\hline
\mbox{\small{$\! C_1 \!\! \setminus \!\! C_2 \! $}} & 2a & 3a & 4a & 7b \\
\hline
2a &  \bullet & 0 & 0  & 7 \\
\hline
\multirow{2}{*}{$3a$} & \multirow{2}{*}{$60$} & \multirow{2}{*}{$\bullet$} & \mathit{64} &  70 \\
&&& \mathit{64} & 63 \\
\hline
\multirow{2}{*}{$4a$}  &  \multirow{2}{*}{$18b$}   & \mathit{36} & \multirow{2}{*}{$\bullet$} & \mathit{28} \\
&& \mathit{36} & & \mathit{28} \\
\hline
\multirow{2}{*}{$7a$}  & \multirow{2}{*}{$4$} & 0 & \mathit{16}   & 0 \\
&& 18a & \mathit{16} & 4 \\
\hline
\multicolumn{5}{c}{\;} \\
\cline{1-3}
\mbox{\small{$\! C_1 \!\! \setminus \!\! C_2 \! $}} & 2B & 6A \\
\cline{1-3}
2B & \bullet  & 18a \\
\cline{1-3}
6A & 70 & \bullet \\
\cline{1-3}    
\end{array}
\;\;\;\;
 \begin{array}{  |  l lllll  l  | }
 \hline
                 \# & M &  g &  \beta_0 & \beta_1 & \beta_\infty & \mbox{Eqn.}    \\
                \hline
                  1+ & A_4 & 0 & 3\; 1 & 2^2 & 3 \;1 &  \eqref{common31}   \\      
                  1+  & S_7 & 0 & 3^2 1 & 2^3 1 & 4\; 3 & \eqref{common43}  \\
                       2+      & S_{18a} & 0 &   3^6 & 2^9 & 7 \; 4 \; 3^2 \; 1 &         \eqref{common74331}    \\
                     1       & S_{18b} & 0 & 3^6 & 2^9 & 7 \; 6 \; 3 \; 1^2 & \eqref{cover76311} \\
                                                                 1 & S_{60} & 1 & 3^{20} & 2^{29} \, 1^2 & 14\; 8^2 \,  7\; 6^2 \, 4^2 3\; &    \\                                
                    1   & A_{63} & 1 &  3^{21} & 2^{32}1 & 14\; 8^2 \,  7\; 6^2 \,  4^2 \, 3^2 &      \\
                                2   & S_{70} & 0 &  3^{23} \, 1 & 2^{35} & 14\; 8^2 \, 7\; 6^2 \, 4^2 \, 3^4 \, 1 &         \\
                                    \hline
                                    \multicolumn{7}{c}{\;} \\
                    \hline
                    \multicolumn{7}{|c|}{\mbox{Two pairs defined over $\Q(\sqrt{2})$}} \\
                    \hline
                    1 & A_{16} & 0 & 3^5 \, 1 & 2^8 & 7 \; 4 \; 3 \; 2 & \eqref{cover7432}  \\
                    1 & A_{28} & 0 & 3^9 \, 1 & 2^{14} & 8 \; 7 \;6 \; 3^2 \, 1 &  \\
                    \hline
 \end{array}
 }
 \]
 \caption{\label{g168}  Left: Degrees of components of Hurwitz-Belyi covers with
parameters $(G,(C_1,C_2),(3,1))$ and $G \in \{PSL_2(7),PGL_2(7)\}$.   Right: further information 
on the covers.}  
 \end{table}

  \begin{eqnarray}
  \nonumber
  f_{7,6,3,1^2}(j,x) & = & \left(9 x^6-102 x^5+295 x^4-212 x^3+39 x^2+90 x+9\right)^3 \\
  \label{cover76311} && \qquad - 2^{14} j  x^6 (2 x-3)^3 \left(9 x^2-66 x-7\right).
  \end{eqnarray}

 The left half of Table~\ref{g168} refers to four pairs of irrational covers.  Letting $s = \pm \sqrt{2}$,
 equations for the smallest degree pair are
  \begin{eqnarray}
  \nonumber \lefteqn{f_{7,4,3,2}(j,x) = } \\
\nonumber  && (-7 x+19 s+27) \left(-49 x^5+x^4 (217 s-63)+x^3 (332 s-478)+ \right. \\ 
\nonumber &&  \left. \qquad x^2 (154 s-658)+  x (196 s+147)+441   s+637\right)^3 \\
\label{cover7432}    && -216 j x^4 (35123 s+49688) (-2 x+s-4)^2 (7 s-4 x)^3.
 \end{eqnarray}

 \subsection{The simple group $SL_2(8)$} 
 The group $T = SL_2(8)$ has outer automorphism group $A$ of order $3$.  
 All the corresponding Hurwitz parameters $h$ satisfying the conditions of \S\ref{23poverview}
 have $G = T$, as those of the form $(T.3,(C_1,C_2),(3,1))$ have 
 at least $\Q(\sqrt{-3})$ in their field of definition and hence break the
 rationality restriction in \S\ref{23poverview}.

  \begin{table}[htb] 
  \[
  {\renewcommand{\arraycolsep}{2.8pt}
\begin{array}{|c|c|c|c|c|}
\hline
\mbox{\small{$\! C_1 \!\! \setminus \!\! C_2 \! $}} & 2a & 3a & 7c & \;\; 9c \;\; \\
\hline
2a & \bullet & 18a & 70 & 54 \\
\hline
3a & 16 & \bullet & 49 & 33 \\
\hline
7a & 84 & \!\! 18b,63 \!\!  & 4,84 & \mathit{81} \\
\hline
7b&      &       & 7,90 &      \\
\hline
9a & 4,36 & 4,36 & \mathit{49} & 33 \\
\hline
9b &          &       &                    & 33 \\ 
\hline
\end{array}
\;\;\;
 \begin{array}{ | l  lllll  l  | }
 \hline
                 \# &  M &  g &  \beta_0 & \beta_1 & \beta_\infty & \mbox{Eqn.}    \\
                 \hline
                  3+& A_4 &   0 &  3 \; 1 & 2^2 & 3 \; 1 &       \eqref{common31}   \\
                       1+&    S_7   & 0 &  3^2 1 & 2^3 1 & 4 \; 3 &      \eqref{common43} \\   
                  1& A_{16} & 0 &  3^4 1^4 & 2^8 & 9 \; 7 &   \eqref{cover97} \\
                     1& S_{18a} & 0 &  3^5 1^3 &  2^9 & 9 \; 7 \; 2 &   \eqref{cover972}   \\
                     1+&  S_{18b} & 0 &   3^6 & 2^9 & 7 \; 4 \; 3^2 \, 1 &         \eqref{common74331}    \\
      3& A_{33} & 0 &  3^{11} & 2^{16} \; 1 & 9 \; 7^3 \; 1^3   & \eqref{cover9777111}  \\     
      2& A_{36} & 1 &  3^{12} & 2^{18} & 9 \; 7^3 \; 3^2   &    \\
     1 & A_{49} & 1 &  3^{16} \, 1 & 2^{24} \, 1 & 9^3 \, 7^3 \, 1 &      \\
       1 & S_{54} & 0 &  3^{18} & 2^{27} & 9^2 \, 7^3 \, 3^3 \, 2^3     &    \\
       1 &  S_{63} & 0 &  3^{20} \, 1^3 & 2^{31} \, 1 & 9^3 \; 7^3 \; 4^3 \; 3  &   \\
          1& S_{70} & 0 &  3^{23} 1 & 2^{35} & 9^3 \, 7^4 \, 3^3 \, 2^3 &     \\
         2 & A_{84} & 1 &  3^{28} & 2^{42} & 9^3 \, 7^5 \, 4^4 \, 3^2 &        \\
      1 & G_{90} & 0 &  3^{30} & 2^{45} & 9^3 \, 7^6 \, 4^3 \, 3^2 \, 1^3 & \eqref{cover9774311}     \\  
\hline   
\end{array}
}
\]
\caption{\label{g504} Left: Degrees of components of Hurwitz-Belyi covers with
parameters $(SL_2(8),(C_1,C_2),(3,1))$.  Right: further information 
on these covers. }
\end{table}

Since the Schur multiplier of $SL_2(8)$ is trivial, there is no spin separation.  However Table~\ref{g504} 
exhibits so many Galois degeneracies that it seems likely that at least some of them are forced 
by deeper reasons.  We describe some of these degeneracies here.   
Our conventions follow the Atlas:  if $g \in 7a$, then $g^2 \in 7b$ and $g^4 \in 7c$; 
similarly, if $g \in 9a$, then $g^2 \in 9b$ and $g^4 \in 9c$.  

\begin{table}
\[
\begin{array}{|ccc|cccc|}
\hline
m & |\langle g_i \rangle| & M & g_1\in 7b & g_2 \in 7b & g_3 \in 7b & g_4 \in 7c \\
\hline
7 &  504 &  S_7 &                  (359467182) & (287516439) & (236478159) & (127698453) \\
90 & 504 & G_{90}       &            (318954762) & (978436512) & (978436512) & (127698453) \\
\hline
\multicolumn{7}{c}{\;} \\
\hline
m & |\langle g_i \rangle| & M & g_1\in 7a & g_2 \in 7a & g_3 \in 7a & g_4 \in 7c \\
\hline
4 & 504 & A_4 & (132674598) & (863972514) & (832465917) & (124835697) \\
84 & 504 & A_{84}  & (396482715) &  (685427319) &  (853716942) & (127698453) \\
9 & 56 & SL_2(8) &  (3 2 9 5 1 8 7 4 6) &  (1 2 4 8 3 5 6 9 7) & 
 (8 2 7 1 5 3 9 6 4) &(1 2 7 6 9 8 4 5 3) \\
9 & 56 & SL_2(8) & (1 3 6 9 2 7 4 8 5) & (1 2 4 8 3 5 6 9 7) &
 (1 8 5 2 9 4 7 3 6)&  (1 2 7 6 9 8 4 5 3) \\
1/7 & 7 & S_1 & (136927485) & (136927485) & (136927485) &  (149375286) \\
\hline
\end{array}
\]
\caption{\label{g504extra} Top: Representatives for braid orbits of $\br_{3,1}$ on $\cF_h$, for
$h = (SL_2(8),(7a,7c),(3,1))$.  Bottom: Representatives for braid orbits of $\br_{3,1}$ on
$\cF_h$ for $h = (SL_2(8),(7b,7c),(3,1))$, followed representatives of three degenerate orbits}
\end{table}

For the case $h = (SL_2(8),(7b,7c),(3,1))$, the mass and degree from the mass formula \cite[(3.6)]{RV15} are
 $\overline{m}= m= 97$.  A braid calculation gives two degeneracies; first there are two orbits, of size $7$ and $90$ respectively.
Second, the monodromy group for the degree $90$ orbit is imprimitive, with image inside the
wreath produce $S_3 \wr S_{30}$.  Representatives in $\cG_h$ of the two braid orbits on $\cF_h$
are given in the top part of Table~\ref{g504extra}.  

The the case $h = (SL_2(8),(7a,7c),(3,1))$, the mass is $\overline{m} = 106 \frac{1}{7}$ and the
degree is $m=88$.  The degree decomposes, $m=4+84$, and the degenerate piece decomposes
as well, $\overline{m} - m = 9+9+\frac{1}{7}$.  The two components with $\langle g_i \rangle = 56$
have the same monodromy group $SL_2(8)$, with the rigid braid partition triple 
$(3^3, \; 2^4 \, 1, \; 7 \; 1^2)$.        Representatives in $\cG_4$ are given in the bottom part of
Table~\ref{g504extra} for all five orbits.  Note that the representative of the
orbit with mass $1/7$ has the very simple form $(g,g,g,g^4)$.    All the braid computations in this paper involve $r$-tuples
of permutations like the ones exhibited in Table~\ref{g504extra}.  

Given how differently behaved the last two parameters are, one might expect that 
the parameters $(SL_2(8),(9a,9c),(3,1))$ and $(SL_2(8),(9b,9c),(3,1))$ would be differently
behaved as well.  However here the Galois degeneracy is in the other direction:
not only is $\overline{m} = m = 33$ in each case, but the two degree 33 covers
are isomorphic.  Moreover, these covers are also isomorphic to the cover arising
from $(SL_2(8), (3a,9c), (3,1))$.  

It is because of the $3$-element outer automorphism group that the covers considered
above are all rational, despite the fact that $7a$, $7b$, $7c$ 
and $9a$, $9b$, $9c$ are defined only over the cyclic cubic fields
$\Q(\cos(2 \pi/7))$ and $\Q(\cos(2 \pi/9))$ respectively.  In contrast, the
three-element group $\Out(SL_2(8))$ is not large enough
to make the covers indexed by  $(SL_2(8),(9a,7c),(3,1))$ and $(SL_2(8),(7a,9c),(3,1))$ rational.
They are each defined over a cyclic cubic field ramified at both $7$ and $9$.  As
reported by Table~\ref{g504}, their degrees are $49$ and $81$ respectively.  
Like most of the covers in the upper right of Table~\ref{g504}, they are full
of genus zero.

Equations for three covers coming only from $SL_2(8)$ are
\begin{align}
\label{cover97} f_{9,7}(j,x) & =  \left(441 x^4+1764
    x^3+702 x^2-140 x+49\right)^3   \\
\nonumber    & \qquad \qquad \left(343 x^4+2940 x^3+6594 x^2-468 x+63\right)  -2^{42} j x^7, \\
\label{cover972}  f_{9,7,2}(j,x) & =    \left(7^4 x^5-441 x^4-3366 x^3+2430
    x^2-3^7 x+3^7\right)^3   \\
\nonumber    & \qquad  \qquad  \left(49 x^2+6 x+9\right) (x+3)      -2^{30} 3^9 j x^9 (x-1)^2,  \\
\label{cover9777111} f_{9,7^3,1^3}(j,x) & = \left(16 x^{11}+256 x^{10}+1312 x^9+2208 x^8  \right. \\
\nonumber & \left.  \qquad  \qquad  -1248 x^7-6720 x^6 -1512
    x^5+5652 x^4  \right. \\ 
    \nonumber & \left.  \qquad  \qquad  -6147 x^3-3912 x^2+11712 x-1536\right)^3 \\
\nonumber    & \qquad + 108 j (x-1) (x+2) (x+8) \left(8 x^3+15 x^2-9 x-8\right)^7  .
\end{align}
The cover $f_{9,7,2}(j,x)$ was found by Hallouin \cite{Hal05}.   For $h = (SL_2(8), (7a, 7b), (3,1)),$
an equation for the degree thirty intermediate cover is
\begin{eqnarray}
\nonumber \lefteqn{f_{9,7^2,4,3,1^2}(j,x) = } \\
\nonumber && \left(11664 x^{10}+31104 x^9-38880 x^8-276960 x^7-458528 x^6-245952
    x^5 \right. \\ 
 \label{cover9774311}   && \qquad \qquad \left. +244440 x^4+549396 x^3+475389 x^2+225504 x+46656\right)^3 \\
\nonumber     && -2^2 3^2 7^7 j   \left(8 x^2+15 x+9\right)^7 x^4 (x-3) \left(3 x^2+6 x+4\right). 
\end{eqnarray}

\subsection{The simple group $G_2(2)' \cong PSU_3(3)$}  In parallel with the \S\ref{235third},
 the third simple group of order $2^a 3^b 7$
 is substantially larger than the first and second group.  Again 
 we present only some sample Hurwitz-Belyi maps, following
 the format used in \S\ref{235third}.  
 
 The first block on Table~\ref{g2} represents cases where the Hurwitz-Belyi map has
 degree $1$ and hence is uninteresting in the present context.  
 These three rigid cases are closely related and are studied in detail in 
 \cite{RobG2}, starting from Proposition~3.1 there.  These three cases serve as a reminder that non-trivial
 Hurwitz-Belyi maps measure a failure of rigidity.  
 
\begin{table}[htb] 
  \[
 \begin{array}{ | ll |  lllll  l | }
 \hline
                C_1 & C_2 & M &  g & \beta_0 & \beta_1 & \beta_\infty &
                \mbox{Eqn.}    \\
                \hline
                 \multicolumn{8}{c}{\;} \\                                     
                                      \multicolumn{8}{c}{\mbox{Degree one covers corresponding to rigid Hurwitz parameters}} \\
                                      \hline
                                      3a & 4a   & S_1 & 0 & 1 & 1 & 1 & x-j \\
                                      4a & 4b   & S_1 & 0 & 1 & 1 & 1 & x-j \\
                                        4a & 2a  & S_1 & 0 & 1 & 1 & 1 &x-j  \\
                                      \hline
                                        \multicolumn{8}{c}{\;} \\                                     
                                      \multicolumn{8}{c}{\mbox{Genus zero covers of small degree}} \\
                                      \hline
                   4a&    6a    & G_4 
                                       & 0 &  3^2 & 2^3 & 4 \; 1^2 & \eqref{common31}  \\
                            4c & 2a & S_7 & 0 &  3^2 \, 1 & 2^3 \, 1 & 4 \, 3 &    \eqref{common43}   \\
                                     4a & 3b &
                                    G_9 & 0 &  3^3 & 2^4 \, 1 & 4 \, 3 \, 2 &   \eqref{common432}   \\
                                    4c & 3a & S_{18} & 0 &  3^6 & 2^9 & 7 \, 4 \, 3^2 \, 1 & \eqref{common74331} \\
                                            4D & 2B & A_{24} & 0 &  3^7 \, 1^3 & 2^{12} & 8 \, 7 \, 6 \, 3 & \eqref{cover8763}    \\
                                                  2B & 4D & A_{40} & 0  & 3^{12} \, 1^4 & 2^{20} & 12 \; 8^2 \, 7 \; 3 \; 2 & \eqref{cover1288732}  \\
                                      \hline
                                      \multicolumn{8}{c}{\;} \\                                     
                                      \multicolumn{8}{c}{\mbox{An unforced splitting to two full covers}} \\
                                      \hline
                                      6a&4b  &S_{135}&  1 & 3^{43} \, 1^6 & 2^{65} \, 1^5 & 14^2 \, 12^4 \, 8^2 \, 6^6 \, 4 \; 3 &   \\
                                      6a&4b  & S_{180} & 3 & 3^{60} &  2^{87} \, 1^6 & 14^2 \, 12^4 \, 8^7 \, 7^2 \, 6^4 \, 3^2 \, 2^2 & \\
                                      \hline
                                      \multicolumn{8}{c}{\;} \\                                     
                                      \multicolumn{8}{c}{\mbox{An example where all allowed bases appear in $\beta_\infty$}} \\
                                      \hline
                 8b&2a&S_{750} & 25 & 3^{248} \, 1^6 & 2^{375} & 
                 {24^5 \, 16^{10} \, 14^{12} \, 12^9 \, 8^8 \, 7^2 }&  \\  
                 &&& & && 
                 {\qquad 6^{13} \, 4^5 \, 3^2 \, 2^5 \, 1^2 } &  \\    
                                      \hline                  
 \end{array}
 \]     
 \caption{\label{g2} Invariants of some covers with $G = G'_2(2)$ or $G_2(2)$}
\end{table}    

The last two genus zero covers on Table~\ref{g2}  come only from $T = PSU_3(3)$.  
Equations are

\begin{eqnarray}
\nonumber \lefteqn{f_{8,7,6,3}(j,x)=} \\
\nonumber  &&  4  \left(4 x^7+22 x^6-60 x^5-166 x^4  +236 x^3+858 x^2-3626 x+2401\right)^3   \\
\label{cover8763} &&  \qquad (2 x-1) \left(2 x^2+16 x-49\right)  \\
\nonumber && + 3^{18} j x^7 (x-2)^6   (x+4)^3, \\
\nonumber & \\
\nonumber \lefteqn{f_{12,8^2,7,3,2}(j,x) =} \\
\nonumber  &  &  \left(64 x^{12}-576 x^{11}+2400
    x^{10}-5696 x^9  +7344 x^8-3168 x^7  -4080 x^6 \right. \\ 
    \nonumber && \qquad  \left. +8640 x^5-7380 x^4 -1508 x^3+8982
    x^2-7644 x+2401\right)^3  \\
\label{cover1288732}    && \qquad    \left(4 x^4-20 x^3+78 x^2-92 x+49\right) \\
\nonumber &&  -2^8 3^{12}  j  \left(2 x^2-4 x+3\right)^8 x^7 (x-2)^3  (x+1)^2 .
\end{eqnarray}

\section{Some $5$- and $6$-point Hurwitz-Belyi maps }
\label{fivepoint}
      All the explicit Hurwitz-Belyi maps presented in the paper so far 
  have had ramification number $r=4$.   \S\ref{fivepointBP}-\ref{fivepoint6}
  present some Hurwitz-Belyi maps with $r=5$ and \S\ref{sixpoint}
  presents one with $r=6$. 
      
  \subsection{A Belyi pencil for $\nu = (4,1)$ yielding  $3$-$2$-$\infty$ maps}
  \label{fivepointBP}
  Sections~\ref{23p}-\ref{237} built many Hurwitz-Belyi maps from the single 
  Belyi pencil $u_{3,1}$ into $\AConf_{3,1}$.    This pencil
  has the remarkable property that it produces braid 
  permutation triples $(b_0,b_1,b_\infty)$ in $S_m$ with $b_0$ and $b_1$ of order
  $3$ and $2$ respectively.  This property kept genera
  very low in \S\ref{23p}-\ref{237}.  
  
  Abbreviating $k=j-1$, let 
  \begin{equation}
  \label{belyipencil41}
 s(j,t) =  k^2 t^4 -6 j k t^2-8 j k t -3 j^2.
  \end{equation}
  Define  $u : \AP^1_j - \{0,1,\infty\} \rightarrow \AConf_{4,1}$ by
  $j \mapsto (D_1(j),\{\infty\})$, with $D_1(j) \subset \AP^1_t$ 
  the roots of $s(j,t)$.   
   Let 
  \begin{align}
  \label{braid41}
  B_0 & = \sigma_1 \sigma_2 \sigma_3^2 , &
  B_1 & = (\sigma_1 \sigma_2 \sigma_3)^2   .
  \end{align}
 A braid calculation says that the abstract braid triple 
  of the Belyi pencil $u$ is  $(B_0,B_1,B_1^{-1} B_0^{-1})$, and that
    $B_0$ and $B_1$  likewise have orders $3$ and $2$ in $\br_{4,1}$ respectively. 
  
  Two Hurwitz-Belyi maps built from $u$ are considered in \cite{RobHNF}. 
  First, for $h = (S_5,(2111,5),(4,1))$ the Hurwitz-Belyi map $\pi_{h,u}$ is
  full and an equation is given in \S4.1 there. This Hurwitz-Belyi
  map reappears in Table~\ref{fivepointtab} here.  For $h = (SL_3(2),(22111,421),(4,1))$
  the degree is $192$.  After quotienting by the natural action of $\mbox{Out}(SL_3(2))$,
  one gets a full degree $96$ map with equation given in \cite[\S8.2]{RobHNF}.

  \subsection{A table of $3$-$2$-$\infty$ maps from $T=A_5$} 
  \label{fivepoint5}
  We begin with the smallest nonabelian simple group $T = A_5$ 
  and build our Hurwitz parameters from $G \in \{A_5,S_5\}$.
  Table~\ref{fivepointtab} gives all Hurwitz-Belyi maps $\pi^*_{h,u} : \AX^*_{h,u} \rightarrow \AP^1_j$ with $h = (G,(C_1,C_2),(4,1))$
  and $u$ the Belyi pencil \eqref{belyipencil41}.    The complications described at the end of \S\ref{BP2} arising in the 
  passage from $(3,1)$ to $(4)$ do not arise when one passes from $(4,1)$ to $(5)$.  Accordingly, Table~\ref{fivepointtab} also includes
  cases of the form $h = (G,(C_1),(5))$,
  written on the table as $h = (G,(C_1,C_1),(4,1))$.   Otherwise,
  Table~\ref{fivepointtab} has a format very similar to the first two tables
  in each of Sections~\ref{235} and \ref{237}.  
     \begin{table}[bht]
   \[
  {\renewcommand{\arraycolsep}{3pt}
  \begin{array}{|c|c|c|c|}
  \hline
\mbox{\small{$\! C_1 \!\! \setminus \!\! C_2 \! $}} & 221 & 311 & 5b \\
  \hline
  221 & 96 & 135 & 150 \\
  \hline
  \multirow{2}{*}{311} &    \multirow{2}{*}{288}  & 126 & 225 \\
  & & 216 & 300 \\
  \hline
   \multirow{2}{*}{$5a$} &  \multirow{2}{*}{64} & 9b & 45 \\
   & & 108 & 3b,9a \\
   \hline
    \multirow{2}{*}{$5b$}  &  &  &  0 \\
    &  &  & 96 \\
    \hline
    \multicolumn{4}{c}{\;} \\
    \multicolumn{4}{c}{\;} \\
    \multicolumn{4}{c}{\;} \\
    \hline
    \mbox{\small{$\! C_1 \!\! \setminus \!\! C_2 \! $}} & 221 & 311 & 5 \\
    \hline
    21111 & 0 & 0 & 25 \\
    \hline
     \multirow{2}{*}{41} &  \multirow{2}{*}{1440} & 1107 & 1275 \\
     &  & 1089 & 1225 \\
     \hline
     \multirow{2}{*}{32} &  \multirow{2}{*}{336} & 255 & 225 \\
     & & 192 & 300 \\
     \hline
     \end{array}
     \;\;\;\;
  \begin{array}{| lclll  |}
  \hline
                 M& g_X & \beta_0 & \beta_1 & \beta_\infty \\
  \hline
                      S_{3b} & 0 &3 & 2 \; 1 & 2 \; 1  \\ 
                     H_{9a} & 0 & 3^3 & 2^3 \, 1^3  & 6 \; 3  \\ 
                     H_{9b} & 0 & 3^3 & 2^4 \,1 & 6 \; 2 \; 1 \\
                    A_{45} & 0 & 3^{15} & 2^{22} \, 1 & 9^3 \, 6 \; 3^3 \, 2 \; 1^2  \\
                 A_{64} & 1 & 3^{21} 1 & 2^{32}  & 15 \; 9^3 \, 6^2 \,  5 \; 2^2 \, 1\\
                 A_{96} &  0  & 3^{32} & 2^{44} \, 1^8 & 15^3 \, 9^3 \, 5 \; 3^6 \, 1\\
                   A_{108} & 1 & 3^{36} & 2^{48} \, 1^{12} & 15^4 \, 9^3 \, 6 \; 5^2 \, 3 \;  2  \\
                    A_{126} & 0 & 3^{42} & 2^{58} \, 1^{10} & 15^3 \, 9^2 \, 6^8 \,  5 \; 3^3 \, 1 \\
                  A_{135} & 0  & 3^{45} & 2^{60} \, 1^{15} & 15^3 \, 9^6 \, 6^4 \, 5 \; 3^2 \, 1\\
                  A_{150} & 1 & 3^{50} & 2^{70} \, 1^{10} & 15^3 \, 9^6 \,  6^6 \, 5 \; 3^2 \, 2^2\\
                 A_{216} & 3 & 3^{72} & 2^{102} \, 1^{12} & 15^5 \, 9^9 \, 6^8 \, 3^4 \\
                 A_{225} & 7&  3^{75} & 2^{112} \, 1^1 & 15^8 \, 9^4 \, 6^{10} \, 5 \; 2^2  \\
                 A_{288} & 7 & 3^{96} & 2^{140} \, 1^8 & 15^8 \, 9^{11} \, 6^{10} \, 5 \; 2^2 \\
                  A_{300} & 4&  3^{100} & 2^{140} \, 1^{20} & 15^8 \, 9^{12} \, 6^{10} \, 5 \; 3 \; 2^2  \\
                   \hline
                   A_{25} & 0 & 3^8 1 & 2^{10} \, 1^5 & 12 \; 9 \;  4 \\  
                    A_{192} & 1 & 3^{64} & 2^{88} \, 1^{16} & 18^2  \,12^{10} \, 4^2 \, 3^9 \, 1 \\
                     A_{225} & 7 & 3^{75} & 2^{112} \,  1^{1} & 18^2 12^{10} 6^{10} 4^2  1^1 \\
                    A_{255} & 7 & 3^{85} & 2^{122} \, 1^{11} & 18^6 \, 12^{10} \,  4^2 \, 3^6 \, 1 \\
                   A_{300} & 6 & 3^{100} & 2^{140} \, 1^{20} & 18^6  \,12^{10} \, 6^{10} \, 4^2 \, 3 \; 1 \\
                  A_{336} & 9 & 3^{112} & 2^{160} \, 1^{16} & 18^8 \, 12^{10} \,  6^{10} \,  4^2 \,  3 \; 1  \\
                 A_{1089} & 34 & 3^{363} & 2^{528} \, 1^{33} & 18^{24} \, 12^{37} \, 6^{34} \, 4 \; 2^2 \, 1\\
                  A_{1107} & 36 & 3^{369} & 2^{540} \, 1^{27} & 18^{24} \, 12^{38} \,  6^{34} \, 4^2 \, 3 \; 2^2\\
                   A_{1225} & 46 & 3^{408} 1& 2^{612} \, 1 & 18^{24} \, 12^{47} \, 6^{35} \, 4^3 \, 2^3 \, 1\\ 
                   A_{1275} & 40 & 3^{425} & 2^{620} \, 1^{35} & 18^{24} \,12^{51} \, 6^{35} \, 4^3 \, 3 \; 2^3 \\ 
                   A_{1440} & 40 & 3^{480} & 2^{704} \, 1^{32} & 18^{24} \, 12^{60} \, 6^{34} \, 4^4 \, 3^{21} \, 2^2 \, 1 \\
                   \hline
    \end{array}
    }
   \]
\caption{\label{fivepointtab} Invariants of $\pi_{h,u}$ for $h = (G,(C_1,C_2),(4,1))$ and $u$ the five-point pencil \eqref{belyipencil41}.
Top: $G=A_5$  Bottom: $G=S_5$}
\end{table}

 There is one instance of cross-parameter agreement:  the Belyi map for $(A_5,(5a),(5))$ and
$(A_5,(221),(5))$ are isomorphic; this Belyi maps occurs for a third time in the next section, where we get
an equation for it. 
Spin separation is near generic as follows.  If 
$(C_1,C_2)$ contains either $221$ or $2111$, then the Belyi cover $\AX^*_{h,u}$ is always connected.  Otherwise
both $C_1$ and $C_2$ split in the Schur double cover and one has the spin separation $\AX_{h,u} = \AX_{h,u}^{*+} \coprod \AX_{h,u}^{*-}$.
In all cases $\AX_{h,u}^{*\epsilon}$ has one component except that $\AX_{h,u}^{*+}$ is empty for  $(C_1,C_2) = (5a,5a)$ and $\AX_{h,u}^{*-}$ has two components for $(C_1,C_2) = (5a,5b)$.

Several patterns in Table~\ref{fivepointtab} merit comments.  First, the first two monodromy groups under the header $M$ are odd, being
$S_3$ and $H_{9a}=9T13$.   However, from the left half of Table~\ref{fivepointtab}, they arise together as an even intransitive 
dodecic group.  With this packaging, all monodromy groups are even, including $H_{9b} = 9T11$.    Second, just like in all the tables
in the previous two sections, the exponent on $1$ in $\beta_0$ is always very small; however, in contrast to these previous tables,
the exponent on $1$ in $\beta_1$ is not always small.   Finally, a phenomenon present in the tables of the previous 
two sections is more visible here because of the different organization: the general nature of $\beta_\infty$ depends
on whether $G$ is $A_5$ or $S_5$.

\subsubsection*{A subtle three-way agreement} The rest of this subsection describes the three smallest degree covers in Table~\ref{fivepointtab}
and related unexpected tight interconnections.  A cubic cover and then these three covers are given by 
\begin{align*}
f_{3a}(j,y) & = y^3 - 2 j, & f_{9a}(j,x) & =   4 \left(x^3+6 x^2+3 x-1\right)^3 + 3^6 j  x^3, \\
f_{3b}(j,y) & = y^3 - (3y-2) j, & f_{9b}(j,x) & = 4 \left(x^3-3 x^2+1\right)^3-27 j x^2 (x-3).  
\end{align*}
Like the covers \eqref{common31} and \eqref{common432}, these for covers have bad reduction only at $2$ and $3$.  
We are using notation adapted to our current situation, but all four covers also fit into our standard notational framework via
$f_{3a} = f_{3}$, $f_{3b}=f_{2,1}$, $f_{9a}=f_{6,3}$ and $f_{9b}= f_{6,2,1}$.  

Our notation is chosen because both nonic covers $\AX_{9a}$ and $\AX_{9b}$ of $\AP^1_j$ are imprimitive, the
cubic covers $\AY_{3a}$ and $\AY_{3b}$ being intermediate.    The formulas
\begin{align*}
y & = -\frac{2 \left(x^3+6 x^2+3 x-1\right)}{9 x}, & y &= \frac{2}{3} \left(x^3-3 x^2+1\right)
\end{align*}
give the maps $\AX_{9a} \rightarrow \AY_{3a}$ and $\AX_{9b} \rightarrow \AY_{3b}$.  

Nonic extensions of a given ground field with Galois group $9T18$ of order $108$ come in 
twin pairs.   The function fields $\Q(\SX_{9a})$ and $\Q(\SX_{9b})$ are such a pair
over $\Q(j)$.  This means that they are nonisomorphic as nonic fields, but
their Galois closures are isomorphic.   One difference between the two nonic
fields is seen after base change to $\C(t)$.  Then the common Galois group
drops to the distinct monodromy groups $9T13$ and $9T11$ mentioned above of order $54$.  This
difference is seen already at the level of cubic subfields, where the monodromy
groups are $3T1 = A_3$ and $3T2 = S_3$.  

The agreement we have been describing is after specialization to the Belyi pencil $u$.
But there is agreement of a different nature before specialization as well.  Let $h_{9b} = (\tilde{A}_5,(5b+,311+),(4,1))$ 
and $h_{12} = (\tilde{A}_5,(5a+,5b-),(4,1))$.  Consider covering surfaces of the two-dimensional
base $\AU_{4,1}$ of \cite{RobHNF}.   Then the nonic cover
$\AX_{h_{9b}}$ and the dodecic cover $\AX_{h_{12}}$ are resolvents of each other,
with Galois groups $9T26 =\F_3^2.GL_2(\F_3)$ and the isomorphic group 
$12T157$ respectively.   The monodromy groups have index two and are 
$9T23 = \F_3^2.SL_2(\F_3) \cong 12T122$.  

The Hurwitz parameter $h_{9c} = (S_3 \wr S_2,(21111,222,33),(3,1,1))$ studied in \cite[\S5]{RobHNF}
now enters our considerations as follows.   Via the natural involution $(u,v) \mapsto (v/u^2,v^2/u^3)$ 
of $U_{3,1,1}$ \cite[\S7.5]{RobHNF} and the exceptional identification $U_{3,1,1} \cong U_{4,1}$ of
\cite[\S3.4]{RobHNF}, the covering surfaces $\AX_{h_{9b}}$ and $\AX_{h_{9c}}$ 
are isomorphic.   This cross-parameter agreement in the setting of surfaces 
complements the three presented in \cite[\S3.6]{RobHNF}.  On an explicit level, consider the specialization
of $f_{9c}(u,v,x)$ from \cite[\S5.2]{RobHNF} obtained by the substitutions
$u = (j-1)/3j$ and $v=(j-1)^2/27j^2$.  This specialization 
defines the same nonic cover as $f_{9b}(j,x)$.

  \subsection{Two unexpectedly similar $3$-$2$-$\infty$ maps built from $T=A_6$}  
  \label{fivepoint6}
  Consider the two Hurwitz parameters on the left:
  \begin{align*}
  h_{96} & = (A_6,(3111),(5)), & (\beta_0,\beta_1,\beta_\infty) & = (3^{36}, \; 2^{44} \, 1^8, \; 15^3 \, 9^3 \, 5 \; 3^6 \, 1 ), \\
   h_{192} & =  (A_6,(3111, \; 2211),(4,1)), &  (\beta_0,\beta_1,\beta_\infty) & = (3^{64}, \; 2^{84} \, 1^{24}, \; 15^3 \, 12^5  \, 9^3 \,  6^8 \,  5 \; 4 \; 3).
  \end{align*}
  These cases are amenable to a standard calculation because the 
  five-point covers $\AY_x \rightarrow \AP^1$ all have genus zero.
  A mass formula calculation says that the two parameters have
  their indicated degrees.  
  
  Since the $\AY_x$ have genus zero, Serre's theorem \eqref{serretheorem} applies and the degree $96$ cover
  $\AX^*_{96} := \AX^*_{h_{96},u}$ 
  does not exhibit spin separation, as $\AX_{96}^{*} = \AX_{96}^{*-}$.   The braid monodromy computation using \eqref{braid41}
  shows that in fact the monodromy group  is $A_{96}$.  The braid partition 
  triple is as indicated above, and so $\AX_{96}^*$ also has genus zero.  
  The standard computation eventually yields  
  $f_{15^3,9^3,5,3^6,1}(j,x) =$
   \begin{align*}
   &\left(3 x^8-6
    x^7-60 x^6+202 x^5-110 x^4-74 x^3-52 x^2-10 x-1\right)^3   \\ &
     \left(729
    x^{24}-10206 x^{23}+15552 x^{22}-2045790 x^{21}+52397442 x^{20}-543319218
    x^{19} \right.  \\ & \left. \qquad +3209261832 x^{18}-12210163074 x^{17}+31525143435 x^{16}-55955395164
    x^{15} \right.  \\ & \left.  \qquad+66094935696 x^{14}-43882703964 x^{13}-2654708692 x^{12}+42096515820
    x^{11} \right.  \\ & \left. \qquad -51857004992 x^{10}+37353393228 x^9-17942013057 x^8+5711207034
    x^7 \right.  \\ & \left.  \qquad-1071984720 x^6+65222394 x^5+12734514 x^4-1277306 x^3-182088 x^2 \right.  \\ & \left.  \qquad -3850
    x-3\right)^3 \\
  & + 2^{10} j   \left[3 x^3-7 x^2+11 x-1\right]^{15}   \left[3 x^3-9 x^2+3 x+1\right]^9 \left[x-3\right]^5 
  \\ & \qquad
  \left(x^4+8 x^3-36
    x^2+17 x+1\right)^3 \left[x-1\right]^3 \left[x\right].
  \end{align*}
  
  \begin{figure}[thb]
\includegraphics[width=4.5in]{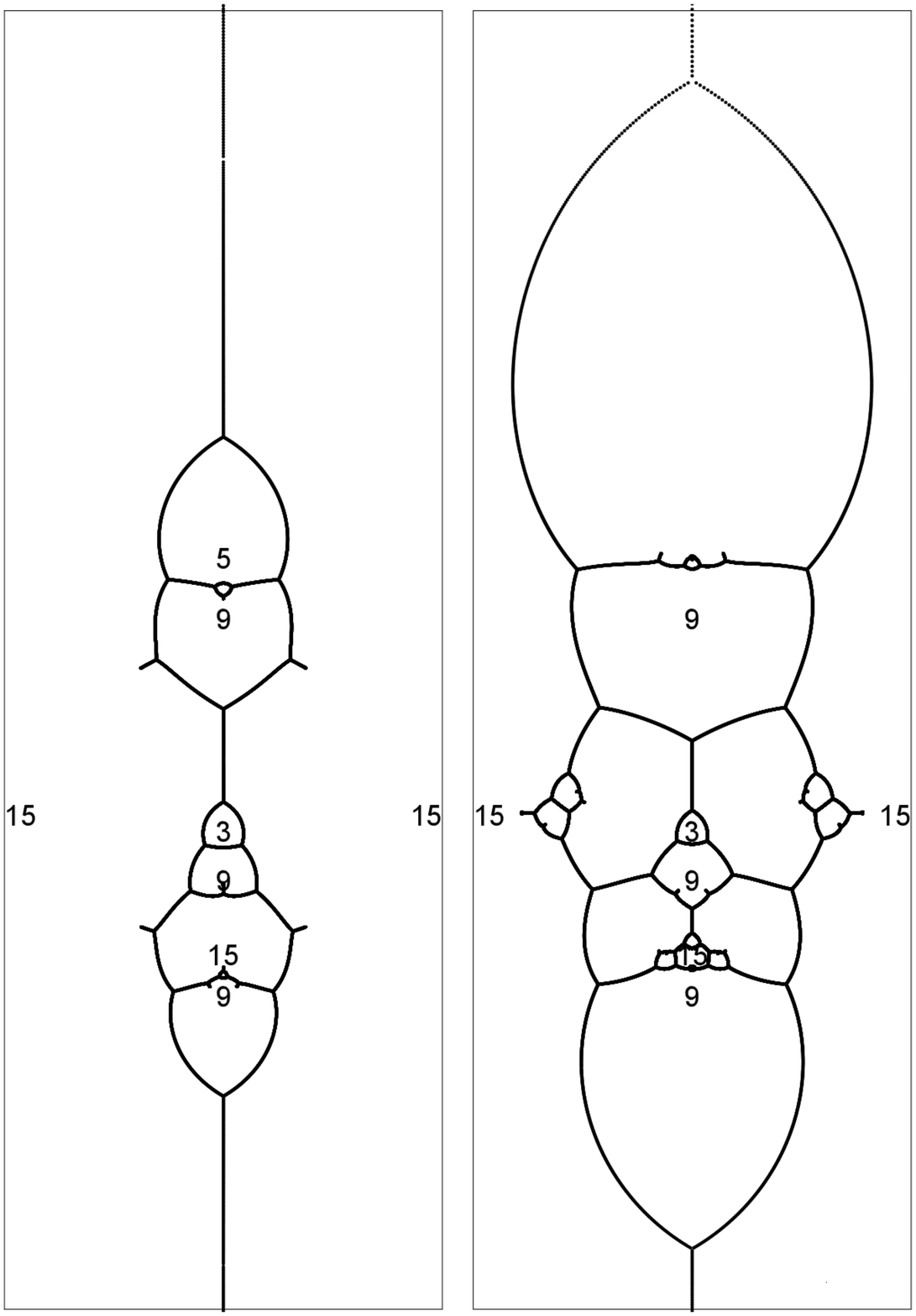}
\caption{\label{twodessins}  Dessins in $\AX^*_{h,u}$ for $h = (A_6,(3111),(5))$ on the
left and $h = (A_6,(3111, 2211),(4,1))$ on the right, illustrating the common
locations of the three 9's and the three 15's
 }
\end{figure}

  The case $h = h_{192}$ has monodromy group $A_{192}$, braid partition
  triple as above, and genus zero.    
  There is a remarkable and unexpected similarity between
  the coefficients of $j$ in the two defining equations, symbolized 
  by 
    \[  
 {15}_A^3 \, {9}_B^3 \, {5}_3 \,  3^4 \, 3 \; 3_1 \, {1}_0  \sim
  {15}_A^3 \, 12^5  \, {9}_B^3 \,  6^8 \,  {5}_3 \, {4}_0 \, {3}_1.  
  \]
  Here the subscripts $3$, $1$, and $0$ indicate that we are normalizing
  so that the coordinates induced  on $\AX^*_{96}$ and  $\AX^*_{192}$ have
  some similarity.   The unexpected similarity is that 
  the cubic polynomials corresponding to the
  two $A$'s coincide and likewise the cubic polynomials corresponding
  to the two $B$'s coincide.     All these agreeing factors are
  bracketed in the two displayed polynomials.  The second polynomial is too 
  large to print, but an excerpt containing the part relevant for the current discussion is 
        \begingroup
\allowdisplaybreaks
    $f_{15^3 ,12^5 ,9^3 ,6^8 ,5, 4, 3}(j,x) = $
   {
  \begin{align*}
&   \left(14659268544 x^{64}-1012884030720 x^{63}  +33879848424192
    x^{62} + \cdots  \right. \\ 
    & \left.  \qquad +40857490944 x^5    -1245316608
    x^4+28200960 x^3 -569088 x^2  +11008 x      -64\right)^3 \\ 
 &  -2^4 3^6  j \left[3 x^3-7 x^2+11 x-1\right]^{15} 
 \left(6 x^5-36 x^4+72 x^3-64 x^2+23  x-4\right)^{12}  \\ & \qquad 
  \left[3   x^3-9 x^2+3 x+1\right]^9  \\ & \qquad 
   \left(9 x^8-72 x^7+240 x^6-444 x^5+474 x^4-280 x^3+72 x^2-12
    x+1\right)^6   \\ & \qquad
  \left[x-3\right]^5 \left[x \right]^4 \left[x-1 \right]^3.
   \end{align*}}%
   \endgroup

Our situation presents many challenges.  For example, we have not worked out
equations for the four covers of largest degree on Table~\ref{fivepointtab} with $g_X=0$.  From 
the degrees given in the table, 45, 96, 126, and 135, the last three are certainly beyond
current implementations of the braid-triple method.  However, if some part of these
equations could be determined ahead of time, perhaps by understanding
better how parts of $f_{96}(j,x)$ repeat in $f_{192}(j,x)$ as just discussed, 
these computations might be brought into the range of feasibility.  

As a second example of a challenge, it would be interesting to build analogs of Table~\ref{fivepointtab}
both for other simple groups $T$ and other Belyi pencils $u$.   The braid monodromy programs described
in  \cite{MSV03} would allow one to go quite far.  For example, consider the Hurwitz parameter
$h = (S_6,(6,51),(4,1))$.   Both classes split in the double cover $\tilde{S}_6$, so one 
has a decomposition $\AX_{h,u} = \AX_{h,u}^+ \coprod \AX_{h,u}^-$.  The mass formula
applied to the group $\tilde{S}_6$ says that the degrees are 49275 and 65400 respectively.  Magaard 
has verified that indeed both $\AX_{h,u}^{\epsilon}$ are full over $\AP^1_j$, with
monodromy groups $A_{49275}$ and $A_{65400}$.  

\subsection{A Hurwitz-Belyi map with $r=6$: wildness at a prime not dividing $|T|$}
\label{sixpoint}
In all the Hurwitz-Belyi maps $\pi_{h,u}$ so far in this paper, 
the bad reduction set $\cP_h$ contains
the bad reduction set $\cP_u$.   
This is not at all the case in general,
and we present an explicit example with $\cP_h = \{2,3\}$ but
$\cP_u = \{2,3,5\}$.  

To get $\cP_h=\{2,3\}$ we leave our main context of almost simple 
$G$ and take $G=S_3$.   The Hurwitz parameter $h = (S_3,(21),(6))$ 
has degree $m=40$ and monodromy group 
$PSp_4(3) \subset A_{40}$. We have computed a polynomial for the three-dimensional slice 
of the Hurwitz cover $\AHur_h \rightarrow \AConf_6$ corresponding
to $D_f \cup \{\infty\}$ with $D_f$ the roots of 
\[
s(t) = t^5 + a t^3 + b t^2 + c t + d.
\]
It has  $1673$ terms as an expanded element of $\Q[a,b,c,d,t]$.
As covering curves $\AY_x \rightarrow \AP^1$ of type $h$ have
genus one, the standard method can
only be used after substantial modification.   Our modification involved that 
points on the cover corresponds classically to three-torsion
subgroups of the Jacobian of the genus two curve $y^2 = s(t)$.

By specializing at suitable Belyi pencils, one can produce many 
degree 40 Hurwitz-Belyi maps $\pi_{h,u}$ with monodromy
group $PSp_4(3)$ and bad reduction set only $\{2,3\}$.  
Here, however, we are trying to explicitly illustrate that
bad reduction at $5$ can be imposed.  So as a Belyi pencil 
we take $u : \AP^1_w \rightarrow \AU_{6} : w \mapsto (D_f(w)\cup \{\infty\})$ 
with $D_f(w) \subset \C_t$ the roots of 
 \begin{equation}
 s(w,t) = 3 t^5 (w-1)^2-10 t^3 (w-1) w+15 t w^2+8 w^2.  
 \end{equation}
 The discriminant is $\disc_t(s(w,t)) = 2^{12} 3^4 5^2$.
 
The specialized polynomial in simplified form is then
 \begin{eqnarray*}
 \lefteqn{f_{40a}(w,x)  = } \\
 \!\!\!& \!\!\! & \left(x^4-30 x^3-240 x^2-450 x-225\right)^5 \left(x^4+30 x^3+240
    x^2+450 x-225\right)^5 \\ 
 \!\!\!    &\!\!\!& -2^4 3^3 5^4 w   \left(x^4+10 x^3+60 x^2+150 x+75\right)^6   \\ &&\;\;\; \left(x^4+30
    x^3+300 x^2+1050 x+675\right)^2 
    \left(x^4-60 x^3-510 x^2-1200 x-675\right) x.
 \end{eqnarray*}
 There are many other Belyi pencils $u$ in $\AConf_6$ which also have $5 \in \cP_u$.
We have chosen the relatively complicated $s(w,t)$ 
 because $f_{40a}(w,x)$ has the very unusual property that the resolvents $f_{27}(w,x)$, $f_{36}(w,x)$, $f_{40b}(w,x)$, $f_{45}(w,x)$ 
 corresponding to other maximal subgroups of $PSp_{4}(3).2$ also have genus zero.
 For example the resolvent
 \[
f_{27}(w,x) = 3^6 x^5  \! \left(x^2-5\right)^5 \! 
    \left(x^2+5 x+10\right)^5 \! \left(2 x^2-5 x+5\right) +  5^4 w \!  \left(3 x^4+10 x^3+25\right)^6
\]
 simplifies the analysis of ramification.   The braid partition triples $(\beta_0, \beta_1, \beta_\infty)$ corresponding to $f_{40a}(w,x)$
 and $f_{27}(w,x)$ are respectively $(5^8, \; 2^{20}, \; 6^4 \; 3 \; 2^4 \; 1^5)$ and
 $(5^5 \, 1^2, \; 2^{10} \, 1^7, \; 6^4 \, 3)$.

\section{Expectations in large degree }
\label{conj}  
In \cite{RV15} with Venkatesh and then in the sequel \cite{RobHNF},
we formulated and supported an unboundedness conjecture 
for number fields.  This final section transposes these
considerations from number fields to 
Belyi maps, with emphasis on phenomena
particular to the Belyi map setting.

\subsection{Full Belyi maps with at most two bad primes}
\label{full2}
Consider Belyi maps defined over $\Q$ with bad reduction
within a given set of primes $\cP$.  
For any prime $p$ and any exponent $k$, it is elementary 
to get $3^k$ different degree $p^k$ such Belyi maps 
$\AP^1 \rightarrow \AP^1$ with monodromy group a $p$-group and
 bad reduction set 
$\{p\}$ \cite{Rob07FCP}.    For any two distinct primes $p$, $\ell$ and certain
$k$, mod $\ell$-reductions of hypergeometric 
monodromy representations give degree $(\ell^{2k}-1)/(\ell-1)$ 
Belyi maps with primitive monodromy group 
$PSp_{2k}(\ell)$ and bad reduction set 
$\{p,\ell\}$.  Fixing $\{p,\ell\}$, the number of such Belyi maps 
for a given $k$ can be arbitrarily large.  

In contrast, it seems very difficult to construct full
Belyi maps defined over $\Q$ with bad reduction
within a two-element set $\cP$.   Returning
to the inverse problem of \S\ref{inverseproblem}, write 
$B_\cP(m)$ for the number of isomorphism classes
of full Belyi maps defined over $\Q$ with bad reduction
within $\cP$.    If $\pi$ contributes to $B_{\cP}(m)$,
then typically the compositions $\sigma \circ \pi$ 
for $\sigma \in \langle t \mapsto 1-t, t \mapsto 1/t \rangle =  \mbox{Sym}(\{0,1,\infty\})$
all contribute separately, so in a sense the numbers 
$B_{\cP}(m)$ are inflated by a factor of six.  
However the $B_{\cP}(m)$  enter the unboundedness conjecture below
only in a qualitative way, and so this 
duplication is not important to us. 
 
To provide
context for the unboundedness conjecture
and support the discussion afterwards, we 
summarize here what we know about the 
numbers $B_{\cP}(m)$ for $|\cP| \leq 2$.   
The trinomial equation $y^k - k t y + (k-1) t = 0$ gives
a Belyi map ramified exactly at the set $\cP_k$ of prime divisors
of $k (k-1)$.   Thus, as an interesting example, $\cP_9 = \{2,3\}$.  
Otherwise one has only the possibilities involving Mersenne primes $M_r = 2^r-1$
and the Fermat primes $F_r = 2^{2^r}+1$, namely 
$(k-1,k) = (M_r,2^r)$ and $(k-1,k) = (2^{2^r},F_r)$.
In \cite{RobCheb}, we are giving two more sequences of 
covers $T_{k-1,k}$ and $U_{k-1,k}$, also ramified
exactly at $\cP_k$.   Degrees are now larger, being 
$k(k-1)/2$ and $(k-1)^2$ respectively.   Our initial 
degree 64 example \eqref{U89} is $U_{8,9}$.

From \cite{MR05} we know also that $B_{\{2,3\}}(m)$
is positive for $m \in \{28,33\}$.  Otherwise we
do not currently know of any instances with $|\cP| \leq 2$ and
$m \geq 20$ with $B_\cP(m)$ positive beyond 
the three sequences just described.

\subsection{An unboundedness conjecture} 
The following conjecture is a direct analog of Conjecture
1.1 of \cite{RobHNF}:
 
 \begin{Conjecture} \label{mc}  Let $B_{\cP}(m)$ be the number of 
 full degree $m$ Belyi maps defined over $\Q$ with bad reduction
 within $\cP$.  Suppose that $\cP$ contains the set of
 primes dividing the order of a finite nonabelian simple group. 
 Then the numbers $B_\cP(m)$ can be arbitrarily large.  
 \end{Conjecture} 
 
 Our heuristic argument for Conjecture~\eqref{mc}  is essentially
 the same as the argument made in \cite{RV15} and \cite{RobHNF}
 for its number field analog.  Namely we expect that Hurwitz-Belyi maps 
 $\pi_{h,u}$ already
 give enough maps to make $B_\cP(m)$ arbitrarily
 large.
 
 In more detail, given $\cP$ as in the conjecture, there
 is at least one nonabelian finite simple group $T$ with
 $\cP_{T} \subseteq \cP$.    From Hurwitz parameters
 $h = (G,C,\nu)$, with $G$ of the form $T^k.A$ as in
 \cite[\S5.1]{RV15}, supplemented if necessary by rational lifting invariants $\ell$,
 there are infinitely many full covers 
 $\AHur^{*\ell}_{h} \rightarrow \AConf_\nu$ defined over 
 $\Q$ with bad reduction within $\cP$.    
 From \cite[\S8]{Rob15} or \cite{Rob07FCP}, there 
 are infinitely many appropriately matching 
 rational Belyi pencils, even with bad reduction set
 consisting of a single prime.  For Conjecture~\ref{mc} to 
 be false, there would be have to be a systematic
 drop from fullness when one specializes from the 
 full family to the Belyi pencil.  We have seen occasional
 drops from fullness in \cite[\S6]{RobHNF} and on
 some of the tables in \S\ref{235}-\ref{fivepoint}
 here.  However these seem to represent a low
 degree phenomenon, and there is no evidence
 of systematic drops in asymptotically
 large degrees.   
 
We have already noted an important difference between Hurwitz number fields and
Hurwitz-Belyi maps in \S\ref{BP1}.   Namely for the former, the specialization
step is arithmetic, as the ground field becomes $\Q$, but for the latter,
the specialization step stays within geometry, as the ground field
becomes only $\C(v)$.   In particular, it seems
to us that Conjecture~\ref{mc} is more within reach than
its analog, as it may be possible to prove it 
using braid groups.  

\subsection{Complements}
  To conclude very speculatively, say that $\cP$ is {\em anabelian} if it contains 
 the set of primes dividing the order of a finite nonabelian simple group,
 and {\em abelian} otherwise.    This terminology seems appropriate
 to us because we suspect that there are connections between
 the material in this paper and investigations into anabelian geometry
 as defined in \cite{Gro97}.
 
 Conjecture~\ref{mc} gives a partial qualitative response to the inverse problem 
 set up in \S\ref{inverseproblem}.  One could ask for a more complete qualitative 
 response.   A guess we find attractive is 
 \begin{itemize}
 \item If $\cP$ is abelian, then $B_{\cP}(m)$ is eventually zero. 
 \item If $\cP$ is anabelian, then $B_{\cP}(m)$ is unbounded because of Hurwitz-Belyi maps, 
 but still zero for $m$ in a set of density one.  
 \end{itemize}
 We put forward the analogous guess for number fields in \cite[\S4.6]{RobHNF}.
 
 The first bulleted statement is supported by the extreme paucity of known Belyi maps
 contributing to $B_{\{p,\ell\}}(m)$, as reported in \S\ref{full2}.  The
 second part of the second bulleted statement is motivated by
 the exponential dependence of the asymptotic mass 
 formula \cite[(3-7)]{RV15} on the multiplicities $\nu_i$.  
  Evidence either supporting or opposing this vision would
 be most welcome.    
 
  \bibliographystyle{amsplain}
\bibliography{mast4}

\end{document}